\documentclass[trsc,nonblindrev]{informs3} 

\OneAndAHalfSpacedXI 

\usepackage{natbib, hyperref}
\usepackage{booktabs, subfig}

\usepackage{ulem}

\usepackage{longtable}

 \usepackage{multirow}
\usepackage[ruled,vlined]{algorithm2e}
\usepackage[dvipsnames]{xcolor}
 \bibpunct[, ]{(}{)}{,}{a}{}{,}%
 \def\bibfont{\small}%

\DeclareMathOperator*{\maximize}{maximize}
\DeclareMathOperator*{\minimize}{minimize}

\def\blue#1{{\color{black}#1}}  
 
\newcommand{\x}{\pmb{x}}
\newcommand{\z}{\pmb{z}}
\newcommand{\uu}{\pmb{u}}
\newcommand{\vv}{\pmb{v}}
\newcommand{\ignore}[1]{}
\newcommand{\vnote}[1]{$\ll${\bf Vish:}
  {\sf \color{red} #1}$\gg$\marginpar{\tiny\bf VN~}}

\TheoremsNumberedThrough     

\EquationsNumberedThrough    

\MANUSCRIPTNO{TS-2021-0349} 

\begin{document}


\RUNAUTHOR{Luo, et al.}

\RUNTITLE{Efficient Algorithms for Stochastic Ride-pooling Assignment  }

\TITLE{Efficient Algorithms for Stochastic Ride-pooling Assignment with Mixed Fleets }

\ARTICLEAUTHORS{%
\AUTHOR{Qi Luo\thanks{Corresponding author: Qi Luo, qluo2@clemson.edu } }
\AFF{Department of Industrial Engineering, Clemson University, Clemson, SC 29634}   
\AUTHOR{Viswanath Nagarajan}
\AFF{Department of Industrial and Operations Engineering, University of Michigan, Ann Arbor, MI 48109} 
\AUTHOR{Alexander Sundt, Yafeng Yin}
\AFF{Department of Civil and Environmental Engineering, University of Michigan, Ann Arbor, MI 48109} 
\AUTHOR{John Vincent, Mehrdad Shahabi}
\AFF{Ford Motor Company, Dearborn,  MI 48120 }
} 

\ABSTRACT{\blue{Ride-pooling, which accommodates multiple passenger requests in a single trip, has the potential to significantly increase fleet utilization in shared mobility platforms.
The ride-pooling assignment problem finds optimal co-riders to maximize the total utility or profit on a shareability graph, a hypergraph representing the matching compatibility between available vehicles and pending requests. 
With mixed fleets due to the introduction of automated or premium vehicles, fleet sizing and relocation decisions should be made before the requests are revealed.}
Due to the immense size of the underlying shareability graph and demand uncertainty, it is impractical to use exact methods to calculate the optimal trip assignments.
Two approximation algorithms for mid-capacity and high-capacity vehicles are proposed in this paper; the respective approximation ratios are $\frac1{p^2}$ and $\frac{e-1}{(2e+o(1)) p \ln p}$, where $p$ is the maximum vehicle capacity plus one.  
\blue{The performance of these algorithms is validated using a mixed autonomy on-demand mobility simulator.}
These efficient algorithms serve as a stepping stone for a variety of multimodal and multiclass on-demand mobility applications.}%




\maketitle

%

\section{Introduction} \label{sec:1}

\blue{Shared mobility fosters an ecosystem of connected travelers, disruptive vehicle technologies, and intelligent transportation infrastructure \citep{shaheen2017mobility, ke2020ride}. 
The platform aims to maximize vehicle fleet utilization by combining multiple requests into a single trip while minimizing the quality of service degradation. 
We term the process of assigning trip requests to available vehicles the \textit{ride-pooling assignment problem}, a notion generalized from the ride-pooling option in ride-hailing services \citep{santi2014quantifying}, whereas the underlying business model can vary from microtransit to shared autonomous vehicles.
Efficient ride-pooling assignment algorithms must balance between maximizing the platform's profits, improving the utility of passengers and drivers, and reducing total energy usage and emissions.
Developing innovative ride-pooling assignment algorithms that meet the vast scales of customers and vehicles has been a significant research task in the burgeoning shared mobility industry.}

\blue{When shared mobility applications migrate from the private to public sectors, the computational challenge of the ride-pooling assignment problem rises inevitably with the increasing vehicle capacity.
Most heuristic methods focus on operating ride-pooling with at most two or three groups of customers \citep{santi2014quantifying, ke2021data, erdmann2021combining, sundt2020heuristics}.
Emerging shared mobility applications, such as microtransit, will operate high-capacity vehicle fleets (between 8 and 14 passengers per vehicle) and will necessitate more scalable algorithms \citep{markov2021simulation,hasan2021benefits,tafreshian2021proactive}.
Per this study, a platform can benefit from operating mixed vehicle fleets to accommodate various customer preferences. 
However, including more types of services complicate fleet operations because of this \textit{uncertain} customer preference.}
\blue{Now we introduce the \textit{Stochastic Ride-pooling Assignment with Mixed Fleets} (SRAMF) problem that arises with the diversification of vehicle fleets in shared mobility platforms.}

\subsection{Overview of the SRAMF Problem} \label{sec1-1}
\blue{The SRAMF problem extends the scalable framework developed for the deterministic ride-pooling assignment for homogeneous fleets in \cite{santi2014quantifying, alonso2017demand}.
The main idea here is to decouple  routing and trip-to-vehicle assignment into two sequential tasks based on the notion of ``shareability graphs''.
Given a batch of trip requests and available vehicles in each time interval, the procedure guarantees the anytime optimality \citep{alonso2017demand} such that the resulting ride-pooling assignments attain the exact optimal solution to the joint problem (see Appendix \ref{app:b1}).
The procedure is summarized as follows:}
\blue{
\begin{enumerate}
    \item First, it constructs a shareability graph that describes the matchable relationship between all trip requests  (demand) and available vehicles (supply) in each batch (see Figure \ref{fig:1}) and compute the associated values of each route. 
    \item Next, the platform solves a general assignment problem (GAP) on this graph that maximizes the total matching value.
    \item Finally, the shareability graph is updated by removing occupied vehicles and assigned demand and inserting arriving requests and available vehicles.
\end{enumerate}}


\blue{
Our research focuses on the stochastic extension of  the  assignment problem (step 2 above). We defer the discussion on constructing shareability graphs (step 1 above) to Appendix~\ref{sec:app2}. There, 
we describe a sequence of matching rules that can construct a hierarchical tree of matchable requests and significantly reduce the computational time.
We show that the framework is computationally efficient and can be deployed with a rolling horizon with demand forecast \citep{yang2020optimizing} as well as varying time intervals \citep{qin2021optimal}.}

\blue{While using the shareability graph can reduce the computational burden of} solving dial-a-ride problems, deploying this framework with mixed fleets is nontrivial. 
Blending two or more types of vehicle fleets leads to new operational challenges due to uncertainty. 
\blue{SRAMF} is a natural yet significant generalization of the ride-pooling assignment problem, which is motivated by the following \blue{applications}:

\noindent \textbf{Example 1:} 
A platform provides several classes of mobility services for various market segments. 
For example, Uber operates a standard service (UberX) alongside a luxury service (Uber Black). 
In most instances, the platform tends to match \blue{customers with the requested vehicle class}.
However, if there are excessive demands in the standard class, it may be more advantageous for the platform to dispatch vehicles from the luxury class to the standard class to \blue{avoid} reneging. 

\noindent \textbf{Example 2:}
A ride-hailing platform hires drivers as either permanent employees or freelancers \citep{dong2021optimal}. 
The platform pays a fixed hourly rate to its permanent employees and pays freelancers per finished trip. 
As a result, the platform needs to determine the priority of drivers in the ride-pooling assignment to facilitate the long-term hiring strategy for permanent employees. 

\begin{figure}[!htb]
    \centering
    \includegraphics[width = 0.6\textwidth]{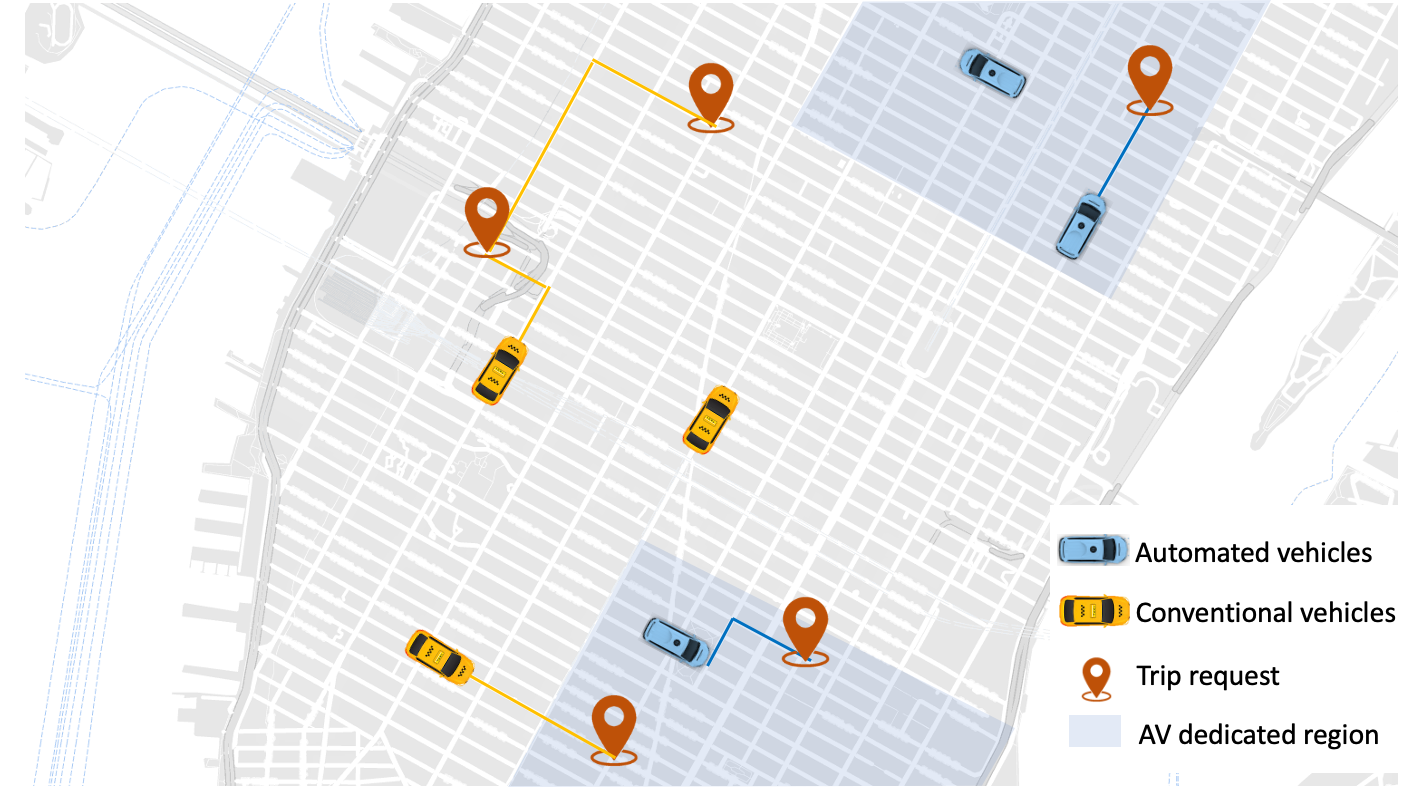}
    \caption{Example of ride-pooling with AVs and CVs. The first-stage decision is positioning vehicles in dedicated regions; The second-stage decisions are solving the assignment problem. }
    \label{fig:0}
\end{figure}

\noindent \textbf{Example 3:} 
A shared mobility platform operates both fully automated vehicles (AVs) and conventional human-driven vehicles (CVs) for on-demand mobility services \citep{wei2020mixed}. 
Operating a mixed autonomy platform faces two potential obstacles (see Figure \ref{fig:0}). 
First, customers may prefer or trust one type of vehicle more than the other type \citep{lavieri2019modeling}. 
Second, AVs requiring dedicated road infrastructure may induce different pickup times and restrict the serviceable areas \citep{shladover2018connected, chen2017towards}.

These examples address the importance of investigating the \blue{mixed-fleet operations problems. 
The platform's decision in SRAMF is twofold: determining each supply source's fleet size and location, which is termed the \textit{vehicle selection decision}, and identifying potential ride-pooling assignments.
The introduction of mixed fleets creates an endogenous source of stochasticity.
The platform must make the vehicle selection decision before the actual demands are revealed (the detailed procedure is described in Figure \ref{fig:1})}.

\blue{It is critical to quantify how these sources of uncertainty affect the efficacy of ride-pooling assignment. 
For example, the platform with mixed autonomy fleets can only direct AVs to specific regions} with roadside units or other road infrastructure \citep{chen2017optimal}.
If the number of AVs in each region \blue{does} not match the estimated demand in the first stage, vehicle repositioning in the second stage \blue{is time-consuming} and detrimental to the platform's profitability.

\subsection{Main Results and Contribution}


\blue{The study's primary objectives are to explain why solving the SRAMF problem is computationally intractable and to develop approximation algorithms for different vehicle fleets.}
\blue{
\begin{enumerate}
    \item We prove that that SRAMF is NP-hard with any finite number of scenarios. Moreover,   its objective does not have favorable submodular properties, which necessitates the development of new approximation algorithms.
    \item Let $p$ denote the largest vehicular capacity plus one.
    For mid-capacity vehicles, we develop a local-search-based approximation algorithm with approximation ratio of $\frac{1}{p^2}$.
    \item For high-capacity vehicles, we develop a primal-dual approximation algorithm with approximation ratio of $\frac{e-1}{(2e + o(1)) p \ln p}$. 
\end{enumerate}
}

\blue{
These polynomial-time algorithms leverage the special structure of shareability graphs to reduce the iteration load of evaluating different positioning policies in the first stage. 
Mid-capacity vehicles can carry less or equal to four requests, suitable for applications in Example~1 and Example~2. 
High-capacity vehicles can carry more than four requests, suitable for on-demand shuttles in Example~3. 
These approximation ratios are close to the best possible: no polynomial-time algorithm can achieve a ratio better than $O(\frac{\ln p}p)$, under standard complexity assumptions. }





\blue{
This work focuses on the profit-maximization setting not only because of the practical concerns \citep{ashlagi2018maximum, simonetto2019real} but also because developing approximation algorithms for maximizing GAPs is generally a more challenging task than minimizing GAPs \citep{fleischer2006tight}.
This study contributes to the growing literature on the operations of ridesharing and other mobility-on-demand platforms as follows:}
\blue{
\begin{enumerate}
    \item Design algorithms with worst-case performance guarantees. \quad  
    This work develops approximation algorithms to facilitate ride-pooling with mixed fleets with provable bounds in the stochastic setting.
    These algorithms are easy-to-implement and resemble the best possible bounds for GAP.
    \item Evaluate the marginal value of increase fleet size in matching. \quad  
    SRAMF necessitates the algorithm to assess the marginal value of including additional vehicles in the existing fleet. 
    This proof technique is of independent interest to the broader fleet sizing studies \citep{benjaafar2021dimensioning}.
    \item Demonstrate algorithms' performance in a mixed-autonomy case study. \quad  
    We demonstrate these algorithms' computational efficiency and optimality gaps in a mixed autonomy traffic case study. 
    The performance of these algorithms is competitive in real-data numerical experiments, and the derived worst-case approximation ratios are conservative.
\end{enumerate}
}

This work focuses on \blue{the rolling-horizon approach with look-ahead demand estimates}.
\blue{Note that the exact method (i.e., MIP-based algorithms) in \cite{alonso2017demand} and \cite{ke2021data} can also be replaced with these approximation algorithms to accelerate vehicle dispatching processes.}
\blue{The} performance guarantee of any approximation algorithm built for the dual-source setting also holds for the single-source \citep{mori2020request} and can be extended to mixed fleets of more than two vehicle types.


\subsection{Organization and General Notation }
The remainder of the paper is organized as follows. We first review the related literature in Section \ref{sec:2}. 
Section \ref{sec:3} \blue{ formulates   the SRAMF} problem and shows its hardness.
Two approximation algorithms are proposed in Section \ref{sec:4} with nearly-tight approximation ratios. 
We test the effectiveness of these approximation algorithms using real-world and simulated data in Section \ref{sec:5} and draw conclusions in Section \ref{sec:6}.

The following notation conventions are followed throughout this work.
The notation $:=$ stands for ``defined as''.  
For any integer $n$, we let $[n]:= \{ 1,2, \dots, n\}$.
We use $v(\cdot)$ as the real value function and $\hat{v}(\cdot)$ as the approximate or estimated value function. P  stands for the class of questions for which some algorithm can provide an answer in polynomial time, and NP stands for those with nondeterministic polynomial time algorithms.
For any set $S$, $|S|$ is its cardinality, i.e., the number of elements in the set. 
Given two sets $A$ and $B$, we let $A+B$ or $A\cup B$ represent the union of $A$ and $B$; we let $A-B$ or $A\backslash B$ represent modifying $A$ by removing the elements belonging to $B$.
We let $A\sim B$ denote that set $A$ intersects with $B$, \blue{i.e.}, $A\cap B \neq \varnothing$.
\blue{i.i.d. stands for ``independent and identically distributed'';
w.r.t. stands for ``with respect to''.
}
\blue{Other notation} and acronyms used in this paper are summarized in Table \ref{table1} in Appendix \ref{sec:app1}.

\section{Literature Review} \label{sec:2}
\blue{We refer to ride-pooling (also called ride-splitting/carsharing rides) in the broad context and focus on operations-level decisions.
The following review covers the recent development of computational methods for ride-pooling applications with different objectives of maximizing the utilization of vehicles or reducing the negative externalities related to deadhead miles.}

\noindent \textbf{\blue{Overview of ride-pooling algorithms.} } \quad  
Solving the optimal ride-pooling assignment is challenging in essence because \blue{a) the number of possible shared trips grows exponentially regarding the number of trip requests and the vehicle capacity; b) trip requests can be inserted during the trip.}
\blue{As a result, the exact approaches for solving the dynamic vehicle routing problem (DVRP) \citep{pillac2013review} are not suitable for platforms that operates thousands of vehicles and expect to compute the assignment solutions in real-time.}

Compared to the substantial body of literature for matching supply and demand without the ride-pooling option \citep{wang2019ridesourcing}, there exist only several attempts to solve the \blue{ride-pooling assignment problem} by heuristic or decomposition methods \citep{yu2019integrated,herminghaus2019mean,sundt2020heuristics}.
\blue{Although these methods achieved satisfying performance metrics in experiments, they cannot balance computational efficiency and accuracy with theoretical guarantees.
This problem is tractable with fixed travel patterns such as providing services for daily commuting.
\cite{hasan2020commute} proposed a commute trip-sharing algorithm that maximized total shared rides for a set of commute trips satisfying various time-window, capacity, pairing, ride duration, and driver constraints.
}

\blue{Since the supply and demand processes are determined by previous decisions, the design of non-myopic policies that considers the future effects of assignments is attractive to platforms. 
\cite{shah2020neural} developed an approximate dynamic programming method that can learn from the IP-based assignment and approximate the value function by neural networks.
\cite{qin2021reinforcement} provided a comprehensive review of the current practice of reinforcement learning methods for assignment and other sequential decision-making problems in the ridesharing industry.  }

\noindent \textbf{Scalable ride-pooling assignment algorithms \blue{for} shareability graphs.} \quad  
To tackle those unprecedented technical challenges in shared mobility platforms, \cite{santi2014quantifying} \blue{quantified} the trade-off between social benefits and passenger discomfort from ride-pooling by introducing the concept of ``shareability networks''. 
They found that the total empty-car travel time was reduced 40\% in the offline setting \blue{(i.e., with ex post demand profiles) or 32\% when demands are revealed en route}. 
This work suffers a limitation in capacity as the matching-based algorithm can only handle up to three-passenger shared rides. 

\cite{alonso2017demand} expanded the framework to up to ten riders per vehicle. 
The high-capacity ride-pooling trip assignment is solved by decomposing the shareability graph into trip sets and vehicle sets and then solving the optimal assignments by a large-scale integer program (IP). 
As the vehicle capacities increase, the moderate size of shared vehicle fleet ($2,000$ vehicles with capacities of four rides in their case studies) can serve most travel demand with short waiting time and trip delay. 
\cite{simonetto2019real} improved this approach's computational efficiency by formulating the master problem as a linear assignment problem. 
The resulting large-scale assignment on shareability networks is calculated in a distributed manner. 
However, despite the easy implementation of these methods, they lack theoretical performance guarantees.

\noindent \textbf{Approximation algorithms for maximization GAP.} \quad 
Approximation algorithms can find near-optimal assignments with provable guarantees on the quality of returned solutions.  
Since the ride-pooling assignment problem is a variant  of GAP \citep{oncan2007survey}, we list the significant results here.
\cite{shmoys1993approximation} and \cite{chekuri2005polynomial} transferred GAP to a scheduling problem and obtained polynomial-time $\frac{1}{2}$-approximation algorithms. 
\cite{fleischer2006tight} obtained an LP-rounding based $(1-\frac{1}{e})$-approximation algorithm and a local-search based $\frac{1}{2}$-approximation algorithm. 

Previous studies have explored GAP algorithms for both instant dispatching and batched dispatching settings. 
Instant dispatching assigns requests to available vehicles upon arrival.
\cite{lowalekar2020competitive} developed approximation algorithms for online vehicle dispatch systems. 
Their setting with i.i.d. demand assumptions is markedly different from the current work. 
Batched dispatching utilizes GAP on a hypergraph to search for locally optimal assignments. 
\cite{mori2020request} developed \blue{$\frac{1}{e}$-approximation LP-rounding algorithms} for the deterministic request-trip-vehicle assignment problem. 
In contrast, the current work considers two-stage decisions of fleet sizing and trip assignment with stochastic demand in \blue{SRAMF}.
As a batched dispatching algorithm, this \blue{stochastic formulation can be applied to} arbitrary demand distributions.



\noindent \textbf{Shared mobility with mixed fleet\blue{s}.} \quad  
\blue{Mixed-fleet mobility systems are emerging research topics in} the ridesharing literature. 
\blue{The first stream of  researches is motivated by the change of workforce structure at present.}
\cite{dong2020managing} investigated the staffing problem in ride-hailing platforms with a blended workforce of permanent employees and freelance workers.
The platform needs to determine the number of hired drivers considering its impact on disclosing a flexible workforce.
\cite{dong2021optimal} justified the dual-source strategy for mitigating the demand uncertainty in ride-hailing systems and designed optimal contracts to coordinate the mixed workforce. 
\cite{castro2020matching} modeled the ridesharing market as matching queues where strategic drivers had different flexibility levels. 
They proposed a throughput-maximizing capacity reservation policy that is robust against drivers' strategies. 

The transition from traditional ridesharing services to \blue{AVs in the foreseeable future motivates a second stream of mixed-fleet researches.}
\cite{lokhandwala2018dynamic} used an agent-based model to evaluate the impact of heterogeneous preferences and revealed that the transition to a mixed fleet would reduce the total number of vehicles, focus on areas of dense demands, and lower the overall service levels in the suburban regions. 
\cite{wei2020mixed} studied the equilibrium of mixed autonomy network in which AVs are fully controlled by the platform and CVs are operated by individual drivers.
The optimal pricing for the mixed service is formulated as a convex program. 
\cite{lazar2020routing} proposed a network model for mixed autonomous traffic and showed how the price of anarchy in routing games was affected. 
In contrast, this work is one of the first attempts to develop algorithms for the operations \blue{of mixed fleets}. 

\section{Problem Description} \label{sec:3}
\subsection{Basic Setting} \label{sec:3-1}
This section introduces the formulation of the SRAMF problem as a two-stage stochastic program and shows its NP-hardness. 
These technical challenges motivate the design of new approximation algorithms in the remainder of this work.

\subsubsection{Preliminaries: construction of shareability graph.}
\blue{ride-pooling assignment is conducted on the  shareability graph, which is represented by a \textit{hypergraph} $G = \{S, D, E\}$. Here, $S$ denotes  the supply/vehicles and $D$ denotes the demand/ride-requests; $S\cup D$ represents the vertices of the hypergraph. 
Each hyperedge $e \in E$ consists of one vehicle and a subset of ride-requests. 
In contrast to  conventional ridesharing where each vehicle can only serve one ride-request per time (which corresponds to assignment on a bipartite graph), we consider the hypergraph  generalization, where each hyperedge $e \in E$ can contain any number of ride-requests (within the vehicle's capacity). 
Other constraints such as detour times and preferences for co-riders are also   considered in the construction of the shareability graph (see Appendix \ref{sec:app2}). 
We also refer to  hyperedges as  \textit{cliques}.
}

The \blue{mixed-fleet} supply contains two separate sets of available vehicles $S_A$ and $S_B$ such that $S = S_A \cup S_B$, $|S_A| = n_A$, and $|S_B| = n_B$. 
The setup costs of vehicles in $S_A$ are significantly higher than those in $S_B$. 
We denote $p = 1+\max_{i\in S_A\cup S_B} \{C_i\}$ where $C_i$ is the capacity of vehicle $i$.
Without loss of generality, we let the cost of using vehicles in $S_B$ be $0$ and let the cost of each vehicle in $S_A$ be normalized to $1$.
Depending on the specific applications, the setup cost can include the extra salary paid to full-time drivers (Example 1 and Example 2 in Section \ref{sec1-1}) or the cost of positioning AVs in advance (Example 3).
The set $S_B$ is called the \textit{basis} supply (which is always available) and the set $S_A$ is called the \textit{augmented} supply (from which we need to select a limited number of vehicles). Section \ref{sec4.3} discusses an extension to more than two vehicle types.


This setting is a generic model for shared mobility applications described in Section \ref{sec1-1}. 
\blue{Each hyperedge $e = \{i, J\}_{i \in S, J\subseteq D }$ corresponds to a potential trip where vehicle $i$ serves requests $J$. }
In the market segmentation example, the augmented supply includes available luxury-service vehicles. 
In the mixed autonomy traffic example, vacant CVs are the basis supply and potential locations to prearrange AVs are the augmented supply.

Solving the GAP problem on a hypergraph is the final step in the dispatching process \citep{alonso2017demand}.
The main analysis for  SRAMF  is conditional on having access to a hyperedge value oracle $\mathcal{O}(v_e)$ that queries the expected profit obtained from any collection of requests in polynomial time. 
The hyperedge values $v_e$ include the total travel time of a Hamiltonian path $t$ picking up all   requests $j \in e$. 
For completeness of deploying these proposed algorithms, Appendix \ref{sec:app2} summarizes the preprocessing procedure.


\subsubsection{\blue{Formulation of SRAMF.}}
Given a budget $K$, the platform chooses a subset $S_R \subset S_A$  of at most $K$ vehicles from the augmented supply. This decision is made before requests are revealed. 
After requests are revealed, the platform can only assign requests to these chosen vehicles $S_R\cup S_B$ and collect profits from finished trips immediately.
Using a hypergraph representation, the second-stage assignment decision is equivalent to choosing a set of hyperedges in which every pair of hyperedges is disjoint. This condition guarantees that each vehicle and each request can be included no more than once in the final assignment.

\begin{figure}[!htb]
    \centering
    \vspace{-0.5in}
    \subfloat[ \label{fig101}]{\includegraphics[angle=270, width = 0.45\textwidth]{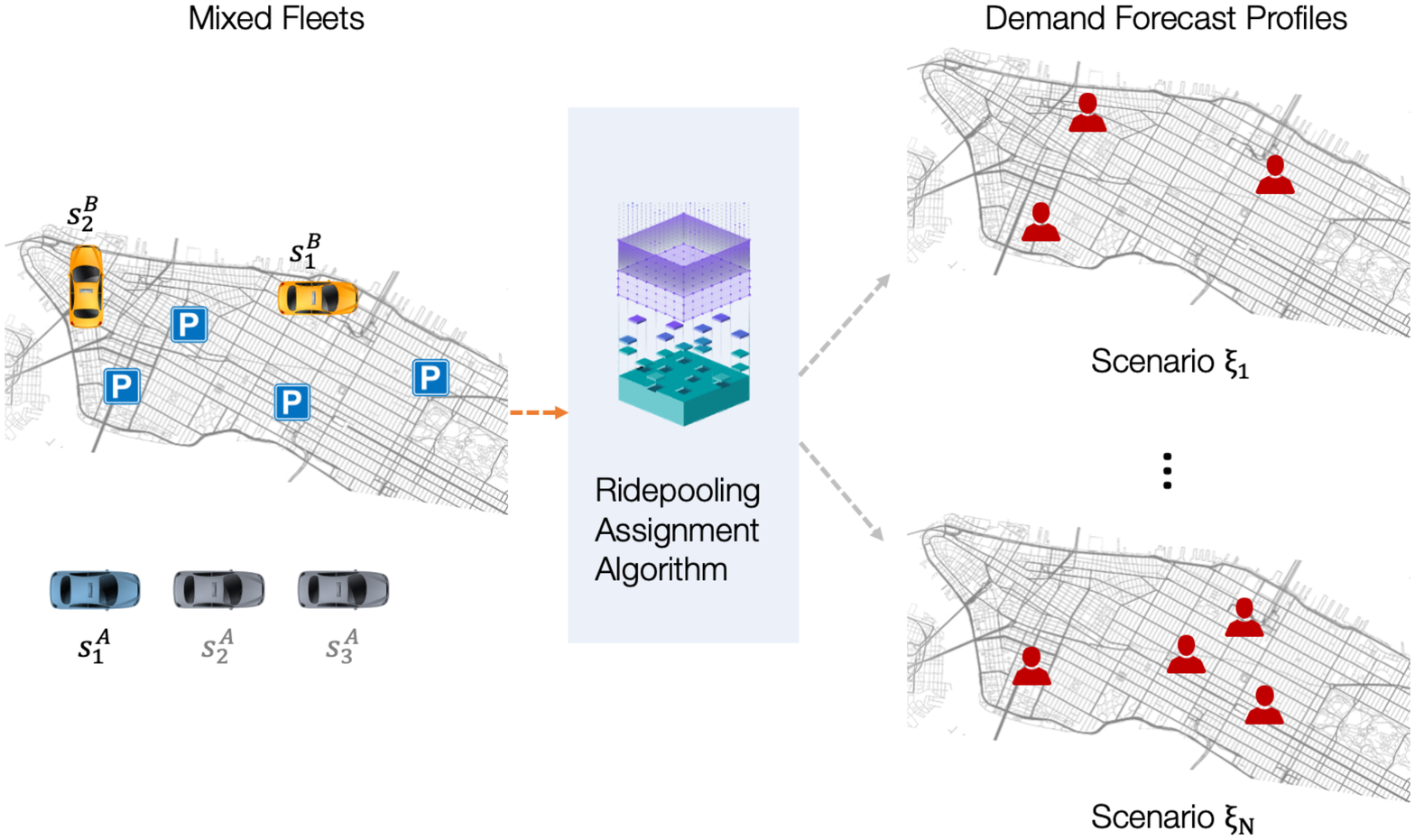}}
    \qquad 
    \subfloat[ \label{fig102}]{\includegraphics[angle=270, width = 0.45\textwidth]{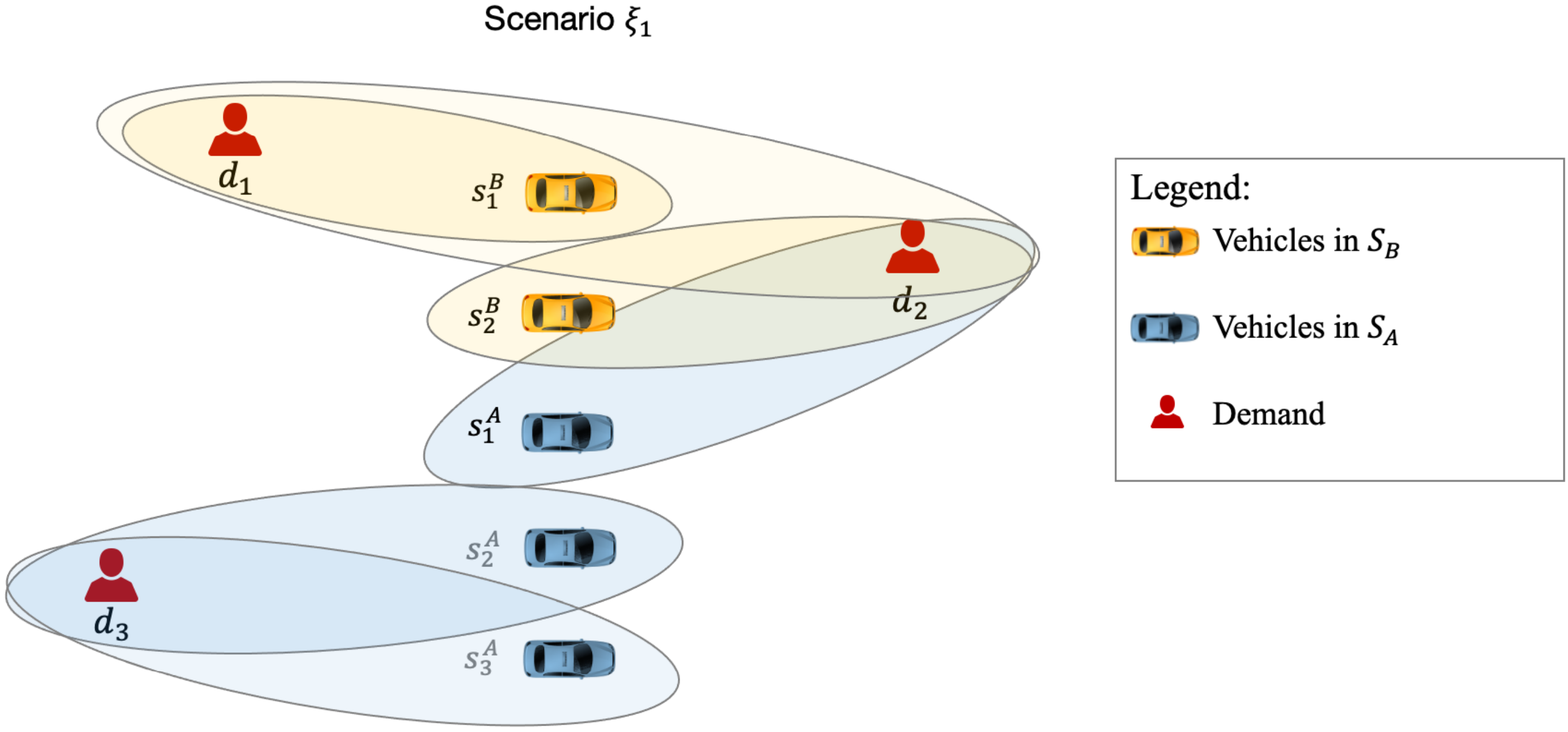}} \\
    \subfloat[\label{fig103}]{\includegraphics[angle=270,width = 0.45\textwidth]{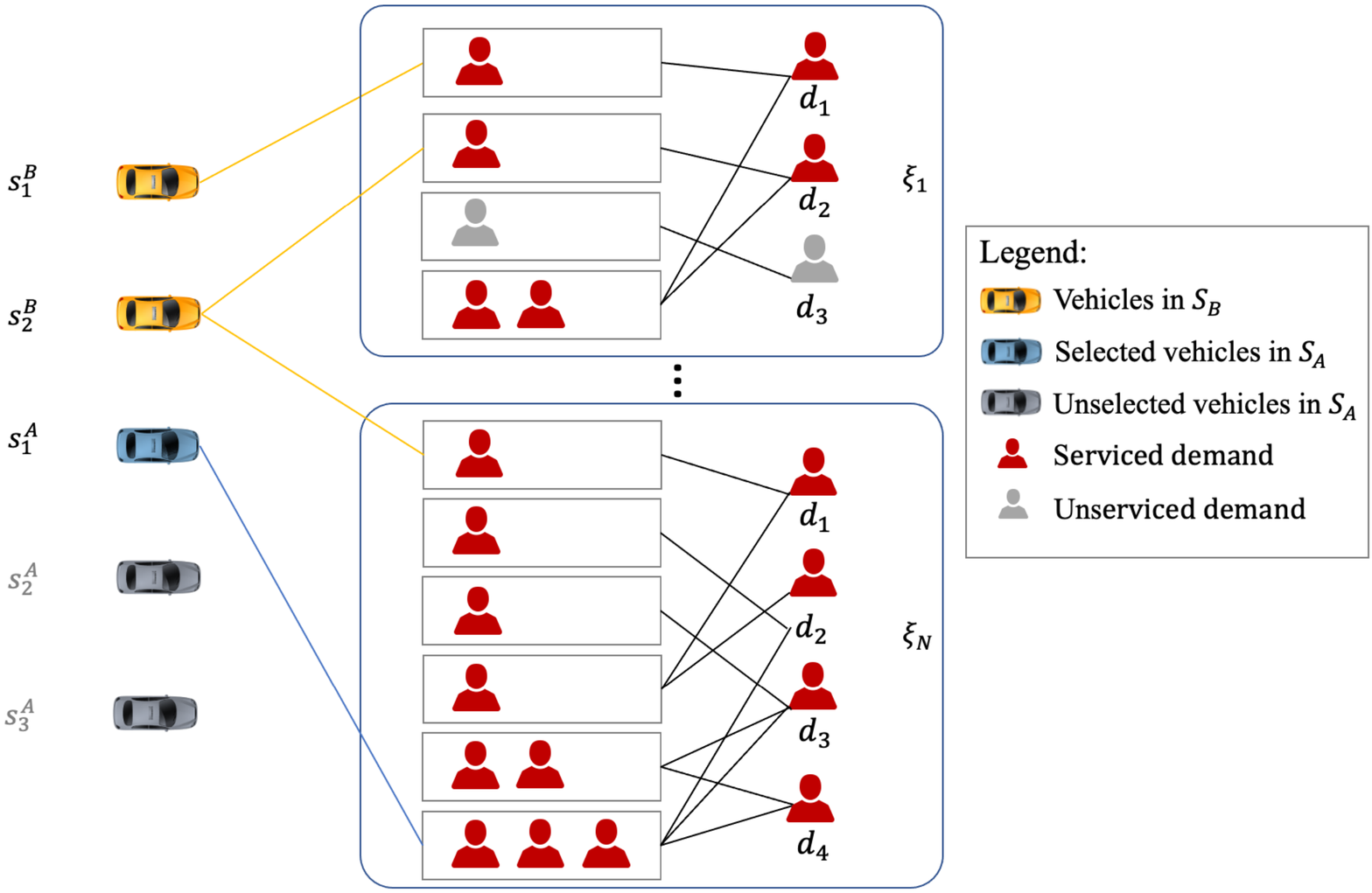}}
    \qquad 
    \subfloat[ \label{fig104}]{\includegraphics[angle=270,width = 0.45\textwidth]{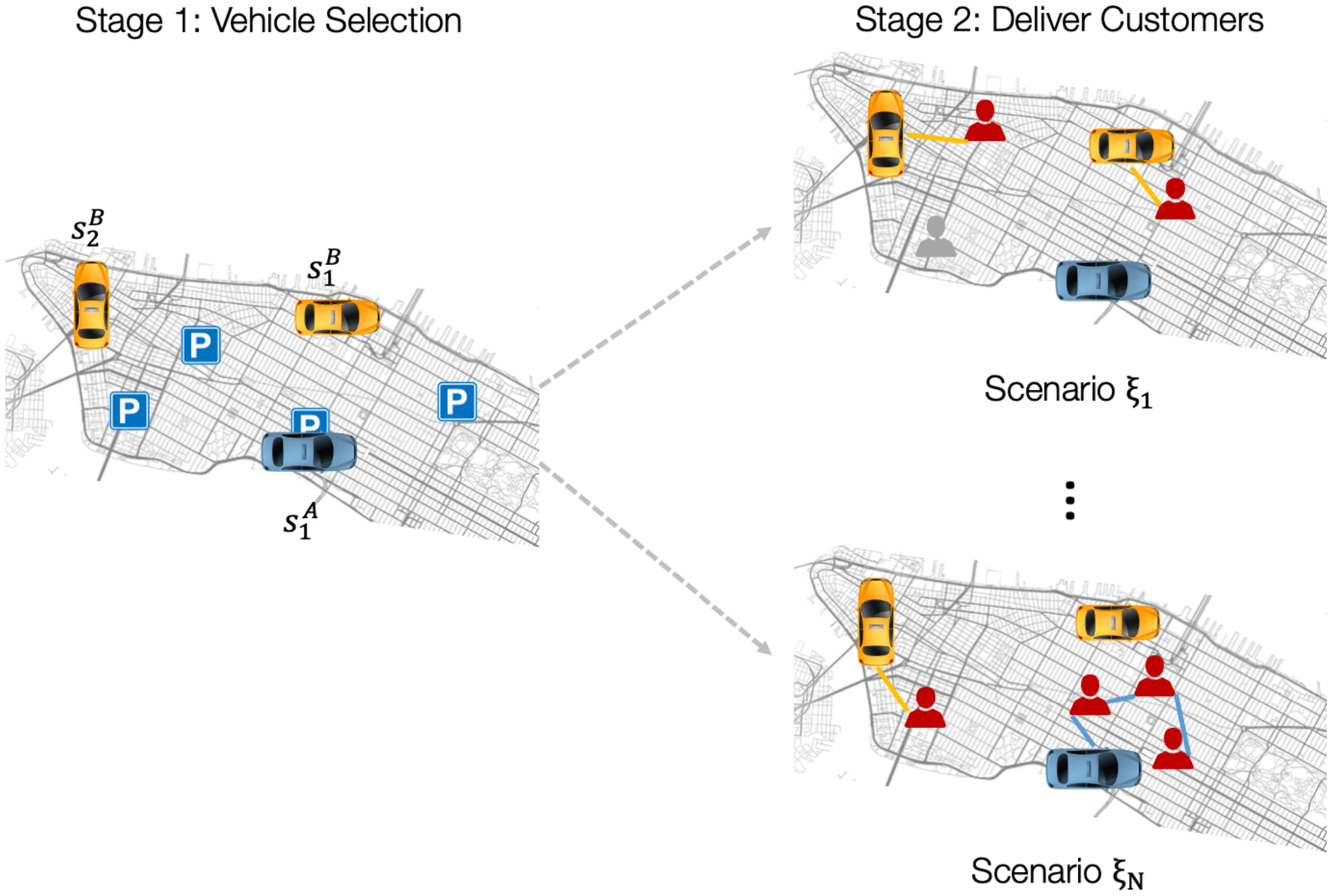}}
    \caption{
    \blue{The illustration of SRAMF procedure per step. 
     $S_B = \{s_1^B, s_2^B\}$ is the basis set (e.g., CVs) and $S_A = \{s_1^A, s_2^A\}$ is the augmented set  (e.g., AVs). 
    In the first step in Figure \ref{fig101}, the algorithm observes the current locations of $S_A$ and $S_B$ and obtain demand forecast.  
    In the second step in Figure \ref{fig102}, the algorithm constructs the shareability graph for each scenario, where each trip is a clique containing one vehicle and multiple matchable requests. 
    In the third step in Figure \ref{fig103}, the SRAMF problem is solved by the approximation algorithm.
    In the final step in Figure \ref{fig104}, the two-stage decisions are implemented and the system state is updated. 
    In each scenario $\xi$, one or more requests are linked by the assigned vehicle in a single ride.}}
    \label{fig:1}
\end{figure}


To be more specific in implementation, the sequence of decision-making in \blue{SRAMF} is as follows:
\begin{enumerate}
	\item For each vehicle $i \in S_A$, we let $y_i = 1$ denote if we put the vehicle in service and $y_i=0$ if not.   
	\blue{These augmented vehicles or potential locations to preposition vehicles are spatially distant from each other and behave differently in the assignment stage.}
We use
	$S_R:= \{i\in [n_A]: y_{i} = 1\} \subseteq S_A$ to denote all \textit{selected} augmented vehicles. 
	All vehicles in the basis supply are included in the first-stage decision as they are no cost to use; i.e., we set $y_i=1$ for all $i\in S_B$.  So, the chosen supply is $S_R\cup S_B$. 
	The first-stage decision space is $Y \in \{0,1\}^{n_A+n_B}$. 
	
	\item   In the second-stage, the platform observes the realized scenario $\xi$, which corresponds to a set of ride-requests $D(\xi)$ and  hyperedges $E(\xi)$; we also observe the   values of all hyperedges. 
	The scenario $\xi$ follows a random distribution $F(\xi)$ with support on $\Xi$, which \blue{is incorporated into} a demand forecast model.
	
	\blue{Each hyperedge $e\in E(\xi)$ includes some vehicle $i$ and \blue{a subset of} requests $J\subset D(\xi)$.  The total number of passengers in  the ride-requests $J$ must be at most the capacity $C_i$ of vehicle $i$, i.e., $\sum_{j\in J} w_j \leq C_i$ where   $\{w_j\}_{j\in J}$ denote the  numbers of passengers in each ride-request.  
	The hyperedge's value} considers following elements:
	\begin{enumerate}
		\item The expected profit $u_j$ \blue{gained from serving the request $j$}.
		\item Trip $t = \{j_1, j_2, \dots:j_k\in J\}$ as a sequence of requests and the associated traveling cost $c(i, t)$ for vehicle $i$ to pick up all requests.
		\item Each \blue{request} $j$ gains additional utility $\tilde{u}_{ij}$ if matched with their preferred vehicle type. 
	\end{enumerate}
	
	The hyperedge value for $e\in E(\xi)$ collected from a potential assignment \blue{is given by}
	\begin{equation}
	v_{e} = \sum_{j\in J} u_j  +  \sum_{j\in J} \tilde{u}_{ij}  - c(i, t) \geq 0.  \label{eq0}
	\end{equation}
	
	\blue{This value  captures many stochastic factors at play between the 1st and 2nd stages, i.e., vehicle selection   and  trip assignment. 
	$u_j$ considers the uncertain number of trip requests and their origin and destination;
	$w_j$ and the set $J$ considers the unknown number of passengers in each trip request; 
	$\tilde{u}_{ij}$ considers the customers' uncertain preference for vehicle types.  
	Finally, due to fluctuating traffic conditions and different vehicle technology (e.g., CVs and AVs), $c(i, t)$ represents that pickup times are uncertain. However, in the 2nd stage (after the scenario $\xi$ is observed), all hyperedge values are known precisely. 
	} 
	
	
	In addition, the batched dispatch \blue{can be expanded to more general policies by using more advanced} value function approximation. 
	\blue{For example,} \cite{tang2019deep} calculated the associated hyperedge value as a reward signal derived from a reinforcement-learning-based estimator.

	\item  The platform assigns  ride-requests to each  available vehicle by determining $x_{e} \in \{0,1\}$ for all $e \in E(\xi)$.   
	The assignment is only available between the chosen supply $S_R\cup S_B$ (denoted as $e\sim S_R \cup S_B$) and realized demand $D(\xi)$ in each scenario. 
	
	

\end{enumerate}

The optimal value of assignments in scenario $\xi$ is denoted by $Q(y, \xi)$ supported on $Q: Y\times \Xi \to \mathbb{R}$. 
Given a scenario, the second-stage decisions are trip assignments denoted by $x = \{ x_e \}_{e\in E(\xi)}$. 
Our objective is to maximize the \blue{\textit{expected total  value}.}

The SRAMF problem can be formulated as a two-stage stochastic program:
\begin{subequations} \label{eq1}
    \begin{align}
        &\maximize_{y}  \mathbf{E} [Q(y, \xi)]   \tag{\ref{eq1}} \\
        s.t. & \sum_{i\in S_A} y_i \leq K   &\textit{(budget) }    \\
        & y_i \in \{ 0,1 \} &\forall i \in S_A \cup S_B,   &
    \end{align}
\end{subequations}
and \blue{the second-stage problem is given by}
\begin{subequations}\label{eq2}
    \begin{align}
         Q(y, \xi) & = \maximize_{x} \sum_{e \in E(\xi)} v_{e} x_{e}  \tag{\ref{eq2}}  \\
        s.t. & \sum_{e\in E(\xi): j\in e } x_{e} \leq 1     & \forall j \in D(\xi)  \quad \textit{(assignment I ) }  \label{eq2d} \\
        & \blue{\sum_{e\in E(\xi): i\in e } x_{e} \leq y_i} & \blue{\forall i \in S_B \cup S_A \quad \textit{(assignment II )}}  \label{eq2e} \\
        & x_{e} \in \{0,1\}       & \forall e \in E(\xi). \label{eq2g}
    \end{align}
\end{subequations}

In the first-stage problem \eqref{eq1}, $K$ is the maximum number of chosen vehicles from the augmented supply. 
In the second-stage problem \eqref{eq2}, the constraints \eqref{eq2d} and \eqref{eq2e} guarantee that each supply and demand are matched at most once \blue{and} only the vehicles selected in the first stage are used.  We allow both vehicles and requests to remain unassigned in the  hypermatching $x$.  
\blue{The second-stage problem is} also known as the \blue{$p$-set packing problem, where ($p-1$) denotes the maximum number of requests per hyperedge} (e.g., ~\cite{furedi1993fractional} and \cite{chan2012linear}). This $p$-set packing problem is already NP-hard. 


\subsubsection{\blue{Road-map for proving SRAMF approximation algorithms.}}
\begin{figure}[!htb]
    \centering
    \includegraphics[width = 0.9\textwidth]{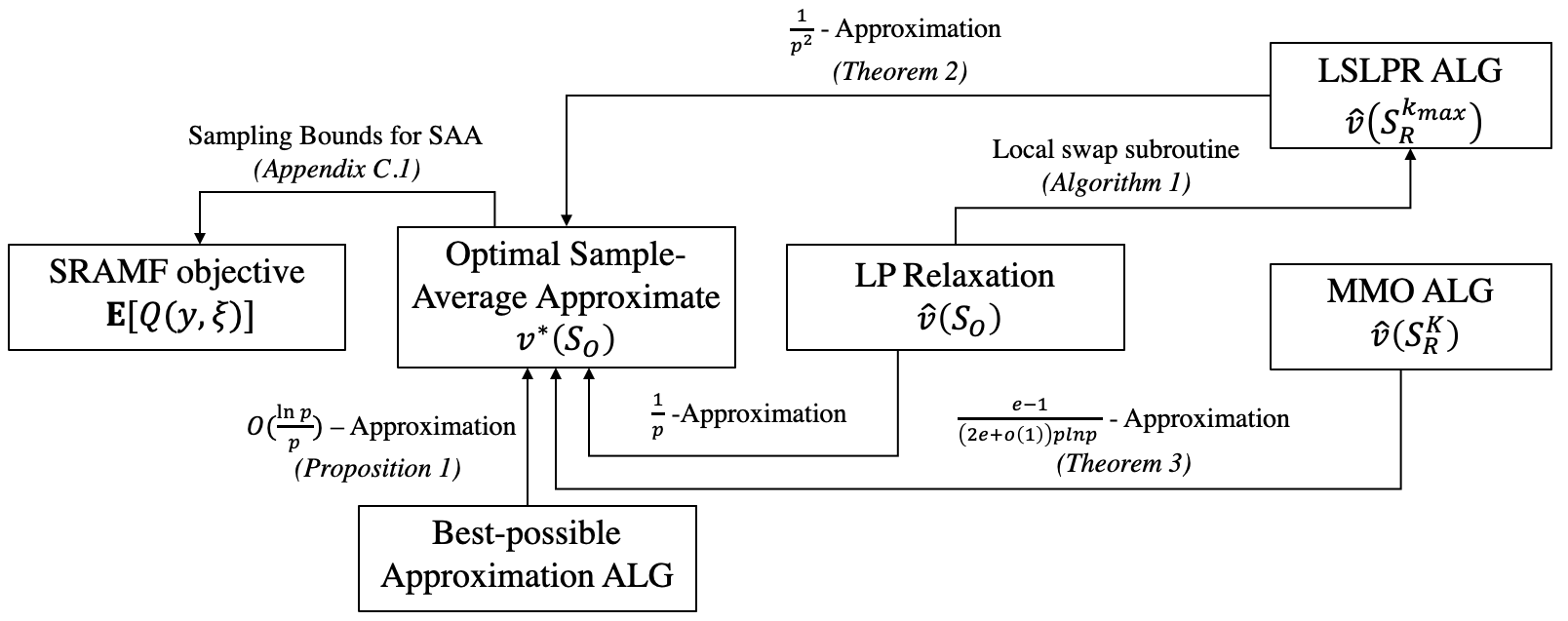}
    \caption{\blue{An illustration of the performance analysis on SRAMF algorithms; the approximation ratios on arrows refer to the results in this paper; $S_O$ is the optimal selection of vehicles and $S_R$ is the section of vehicles generated by approximation algorithms. }}
    \label{fig:r11}
\end{figure}

\blue{Figure \ref{fig:r11} provides an overview of the performance analysis on two proposed approximation algorithms and their approximation ratios, respectively. 
We start with reducing the objective of \eqref{eq1} to the sample-average estimate in Section \ref{sec:3.2}.  
We then show the hardness of the SRAMF problem in Section \ref{sec:3.3}. A key challenge is that  the 2nd stage problem \eqref{eq2} is itself NP-hard. So, our approximation algorithms rely heavily on ``fractional assignments'' that relax the integrality constraints in \eqref{eq2} and can be solved in polynomial-time via linear programming. 
We provide two different approximations algorithms, LSLPR and MMO,  in Section \ref{sec:4.1} and Section \ref{sec:4.2}, respectively. 
}

\subsection{Reduction to Sample-Average Estimate} \label{sec:3.2}

The sample-average approximation (SAA) method is commonly used to solve such two-stage stochastic programs. 
It draws $N$ scenarios $\{\xi_\ell \}_{\ell=1}^N$ from a scenario-generating oracle \blue{(e.g., demand forecasting and vehicle simulation models)} and approximates the expected objective function by a sample-average estimate $\mathbf{E}[Q(y, \xi)] \approx \frac{1}{N} \sum_{\ell=1}^N Q(y, \xi_\ell)$.

To simplify the analysis with regard to $\mathbf{E}(Q(y,\xi))$, we reduce the objective function to finite-sample proximity.
The main analysis is conditional on $N$ mutually disjoint sets of hyperedges as $E(\xi)$. 
\blue{Since the second-stage assignment ensures unique matching per scenario, we can make $n'$ copies when a trip consisting of one vehicle and several requests duplicates across $n'$ scenarios.}
The consistency and shrinking bias of the sample-average estimate are well-studied in literature.
The optimal value of any approximation algorithm converges to $\mathbf{E}[Q(y,\xi)]$ as the number of scenarios $N\to \infty$.
We include the standard SAA proof in \blue{Appendix \ref{sec:app1}} for completeness of the results, which helps determine the required sample size $N$ for any confidence level. 


The current paper's focus is developing algorithms to solve the SRAMF  problem in eq.\eqref{eq1} \blue{with the sample-average estimate.} As mentioned earlier, we will work with an LP relaxation of~\eqref{eq2} as the original problem is NP-hard.  
For any subset $S_R\subseteq S_A$ and scenario $\xi$, define $\hat{v}(S_R,\xi)$ to be the optimal value of the following LP:
\begin{subequations}\label{eq2.5}
    \begin{align}
         \maximize_{x}& \sum_{e \in E(\xi)} v_{e} x_{e}  \tag{\ref{eq2.5}}  \\
        s.t. & \sum_{e\in E(\xi): j\in e } x_{e} \leq 1     & \forall j \in D(\xi)  \label{eq2.5d} \\
        & \sum_{e\in E(\xi): i\in e } x_{e} \leq 1     & \forall i \in S_A \cup S_B    \label{eq2.5e}  \\
        &  x_{e} = 0    & \forall e \sim S_A\setminus S_R \label{eq2.5f}  \\
        &  x_{e} \geq 0       & \forall e \in E(\xi). \label{eq2.5g}
    \end{align}
\end{subequations}

We call the LP solution ``fractional assignments''. Also, define ${v}(S_R,\xi)$ to be the optimal value of the \blue{IP} corresponding to the formulation above; so ${v}(S_R,\xi)$ equals the optimal value in \eqref{eq2}. These are used to define two useful objective functions w.r.t. $S_R$: 
\begin{itemize}
    \item The objective value using the exact GAP in eq.\eqref{eq2} for each scenario is given by:
    \begin{align}
        v^*(S_R) =   \frac{1}{N} \sum_{\ell\in[N]} v(S_R, \xi_{\ell}).
    \end{align}
    \item The objective value using the LP-relaxation \eqref{eq2.5} is given by:
    \begin{align}\label{eq:vhat}
        \hat{v}(S_R) = \frac{1}{N} \sum_{\ell\in[N]} \hat{v}(S_R, \xi_{\ell}).
    \end{align}
\end{itemize}

\blue{Fractional assignments   enjoy many well-known  properties.} 
\blue{First, the integrality gap of the above LP relaxation for $p$-set packing is at most $p$ \blue{~\citep{arkin1998local}}. Second,} a greedy algorithm that selects hyperedges $e$ in decreasing order of their values $v_e$ (while maintaining feasibility) also achieves a $\frac1p$-approximation to the LP value. 
\blue{We restate them in the following theorem:}
\begin{theorem}\label{thm:greedy}
  For any $S_R\subseteq S_A$, we have $    v^*(S_R) \leq \hat{v}(S_R) \leq p  \cdot  v^*(S_R)$; furthermore, the greedy algorithm  obtains a solution of value at least $\frac1p\cdot \hat{v}(S_R)$. 
\end{theorem}

These reductions narrow down the main task to bounding the approximation ratios with regard to $\hat{v}(S_R)$. 
In particular, we will focus on the \blue{SRAMF} problem with fractional assignments:
\begin{equation}\label{eq:frac-SRAMF }
    \max_{S_R\subseteq S_A : |S_R|\le K}\quad \hat{v}(S_R).
\end{equation}

 \blue{ If we obtain an $\alpha$-approximation algorithm for \eqref{eq:frac-SRAMF }, then combined with Theorem~\ref{thm:greedy}, we would obtain an $\frac{\alpha}{p}$-approximation algorithm for SRAMF (with integral assignments). }

Before jumping into the design of approximation algorithms, the following subsection elaborates on some technical challenges.

\subsection{\blue{Hardness and Properties of SRAMF}} \label{sec:3.3}
We show that solving SRAMF  is computationally challenging due to the following reasons.  
First, the second-stage SRAMF problem is NP-hard in general, as demonstrated in Proposition \ref{prop0}. 
Second, Proposition \ref{prop2} shows that $\hat{v}(S_R)$ is   {\em not} submodular, which prevents the use of classic algorithms such as (\cite{nemhauser1978analysis}).     
These facts motivate the development of new approximation algorithms in Section \ref{sec:4}.

\ignore{
\blue{We begin the discussion about the properties of SRAMF with following definitions.}
\blue{
\begin{definition}
The hypergraph $G$ is
\begin{enumerate}
    \item $p$-uniform if all its hyperedges have size $p$.
    \item  $p$-strongly colorable if there is a partition of vertices $V$ into $p$ sets such that each part is an independent set. 
\end{enumerate}
\end{definition}
}

\blue{With vehicles with capacity $(p-1)$, the ride-pooling assignment problem is equivalent to} the \blue{following} $p$-dimensional matching problem:
\begin{definition}\label{def:pdim mathcing}
A $p$-dimensional matching finds a hypermatching of maximum size on a $p$-uniform $p$-strongly colorable hypergraph $G$.
\end{definition}

We can now show the hardness of SRAMF.
}

\begin{proposition} \label{prop0}
There is no algorithm for SRAMF  \blue{(even with $N=1$ scenario)} with an approximation ratio better than $O(\frac{\ln p}{p})$, unless $P=NP$.
\end{proposition}

\proof{Proof for the hardness of SRAMF :} 
\blue{We reduce from the \textit{$p$-dimensional matching} problem, defined as follows. There is a hypergraph $H$ with vertices $V$ partitioned into $p$ parts $\{V_r\}_{r=1}^p$, and hyperedges $E$. Each hyperedge  contains exactly one vertex from each part (so each hyperedge has size exactly $p$). The goal is to find a collection $F$ of disjoint hyperedges that has maximum cardinality $|F|$. }

\blue{Given any instance of $p$-dimensional matching (as above), we generate the following SRAMF instance. The augmented vehicles are $S_A=V_1$  and the basis vehicles are  $S_B = \varnothing$. There is $N=1$ scenario with  ride-requests   
$V_2 \cup \dots  V_p$ and hyperedges  $E$ (each of value $1$).  
Each vehicle has capacity $p-1$ and each ride-request has $1$ passenger. Note that each hyperedge contains exactly one vehicle, as required in SRAMF. 
The bound $K=|S_A|$: so the optimal 1st stage solution is clearly $S_R=S_A$ (select \textit{all} augmented vehicles). 
Now, the SRAMF problem instance reduces to its 2nd stage problem \eqref{eq2}, which involves selecting a maximum cardinality subset of disjoint hyperedges. 
This is precisely the 
$p$-dimensional matching problem. }

\blue{It follows that if there is any $\alpha$-approximation algorithm for SRAMF with $N=1$ scenario  then there is an $\alpha$-approximation algorithm for $p$-dimensional matching. 
Finally, \cite{hazan2006complexity} proved that it is NP-hard to approximate $p$-dimensional matching better than an  $O(\frac{\ln p}{p})$ factor (unless $P = NP$). The proposition now follows. }
\hfill $\square$
\endproof

\blue{
This intractability is the reason that we work with the \textit{fractional} assignment problem  \eqref{eq:frac-SRAMF }. A natural approach for budgeted maximization problems such as \eqref{eq:frac-SRAMF } is to prove that the objective function is \textit{submodular}, in  which case one can directly use the  
$(1-\frac1e)$-approximation  algorithm by \citep{nemhauser1978analysis}. 
However, we show a negative result about the submodularity of $v^*(S_R)$ as well as  $\hat{v}(S_R)$, which precludes the use of such an approach.  Recall that a set function $f:2^\Omega\rightarrow \mathbb{R}_+$ on groundset $\Omega$ is submodular if $f(U\cup \{i\}) - f(U) \ge f(W\cup \{i\}) - f(W) $ for all $U\subseteq W\subseteq \Omega$ and $i\in \Omega\setminus W$.
}

\ignore{
\blue{
\begin{definition}
The set function  $\hat{v}(\cdot)$ is submodular if for every every $S_1$, $S_2 \subseteq \Omega$ with $S_1\subseteq S_2$ and every $s\in \Omega \backslash S_2$, we have that $\hat{v}(S_1 \cup \{s\}) -\hat{v}(S_1) \geq \hat{v}(S_2 \cup \{s\}) -\hat{v}(S_2)$.
\end{definition}
}


\blue{If the set function in eq.\eqref{eq:frac-SRAMF } is submodular, SRAMF could be solved by the classic $(1-\frac1e)$-approximation submodular maximization algorithm \citep{nemhauser1978analysis}.}
However, we show a negative result about the submodularity of $v^*(S_R)$ as well as  $\hat{v}(S_R)$.  
}

\begin{proposition} \label{prop2}
 $v^*(S_R)$ and $\hat{v}(S_R)$ are \blue{{\bf not}} submodular functions.
\end{proposition}

\begin{figure}[!htb]
	\centering
	\includegraphics[width=0.5\textwidth]{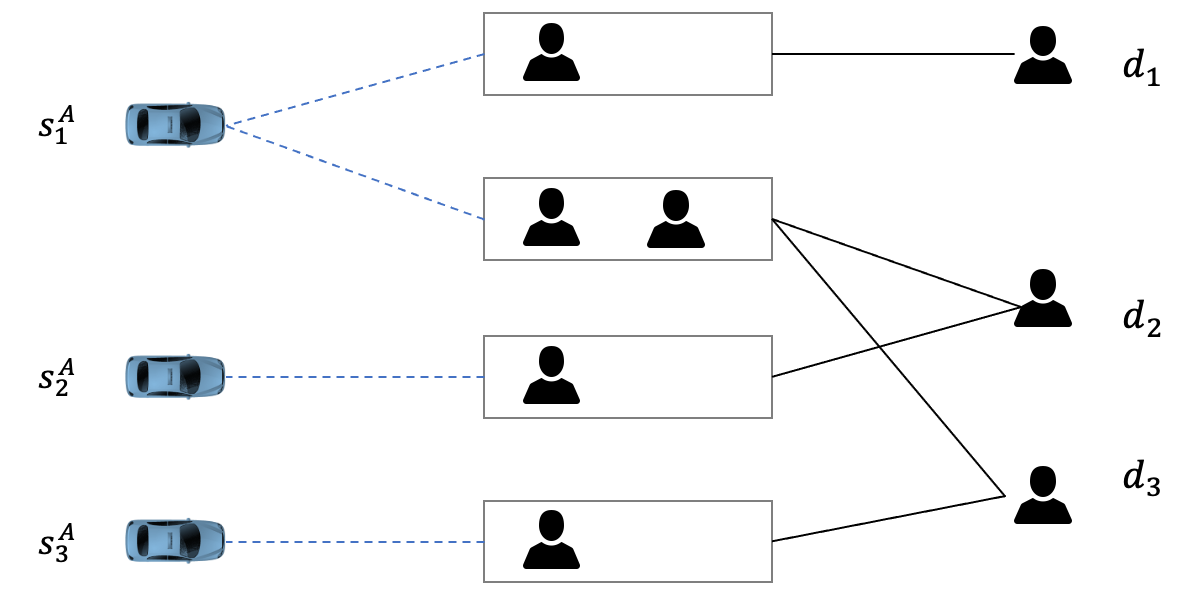}
	\caption{An example for non-submodularity of function $v^*(S_R)$.}
	\label{fig:2}
\end{figure}

\proof{Proof:}
\blue{ Recall that the groundset for both functions $v^*$ and $\hat{v}$ is $\Omega:= S_A$ the set of augmented vehicles. We provide an SRAMF instance with $N=1$ scenario where these functions are not submodular.  
Consider a shareability graph with  $|S_A|=3$,  $S_B = \varnothing$ and three ride-requests $\{d_1, d_2, d_3\}$. 
Let $p=3$, i.e., each vehicle can carry at most two requests. 
The set of hyperedges  is }
\blue{
$$\{(s_1^A, d_1) , (s_1^A, d_2,d_3) , (s_2^A, d_2) , (s_3^A, d_3)\}.$$ 
See also  Figure \ref{fig:2}. The value of each hyperedge is the number of ride-requests covered by it. 
}

\blue{
Let subsets $U= \{s_1^A\}$ and $W=\{s_1^A, s_2^A\}$. Also let $i=s_3^A$. Clearly,
$v^*(U) = 2$ (serving $d_2, d_3$), $v^*(W) = 2$  (serving $d_1, d_2$  or $d_2, d_3$), $v^*(U\cup\{i\}) = 2$  (serving $d_1, d_3$  or $d_2, d_3$), and $v^*(W\cup\{i\}) = 3$. Therefore, we have:
\begin{align*}
    v^*(W\cup \{i\}) - v^*(W) = 1 > 0 = v^*(U\cup \{i\}) - v^*(U),
\end{align*}
which implies the set function $v^*$ is not submodular. It is easy to check that the  LP value function  $\hat{v} = v^*$ for this instance: so function $\hat{v}$ is also not submodular. }
\hfill $\square$
\endproof

\section{Approximation \blue{Algorithms for SRAMF} } \label{sec:4}
This section provides two different approximation algorithms for SRAMF. \blue{Both the algorithms focus on solving the fractional assignment problem \eqref{eq:frac-SRAMF }, and achieve approximation ratios $\frac1p$ and $\approx\frac{e-1}{2e\cdot \ln p}$ respectively.} 
Combined with Theorem~\ref{thm:greedy}, these imply approximation algorithms for SRAMF with an additional factor of $\frac1p$. 

\subsection{Local Search Algorithm for \blue{Mid-Capacity SRAMF} } \label{sec:4.1}
\blue{The Mid-Capacity SRAMF models the current ride-hailing market where the limit for palletizing requests is two or three.}
In this section, we propose a Local-Search LP-Relaxation (LSLPR) algorithm that obtains $\frac{1}{p}$-approximation for the fractional problem~\eqref{eq:frac-SRAMF }. 

\subsubsection{Overview of the LSLPR algorithm. }
Let $\epsilon>0$ be an arbitrarily small parameter that will be  used in the  stopping criterion.  
The outline of the LSLPR algorithm is as follows:
\begin{enumerate}
	\item Start from any solution $S_R\subseteq S_A$   \blue{with $|S_R| = K$}.
	\item \blue{Consider all solutions  \blue{$S_{R'} = S_R- \{i\} + \{i'\} $ where $i \in S_R$ and $i'\notin S_R$ (i.e., swapping one vehicle)} and evaluate the corresponding LP value $\hat{v}(S_{R'})$.  } 
	\item \blue{Change the current solution $S_R$ to $S_{R'}$ if the objective value improves significantly, i.e.,   $\hat{v}(S_{R'}) > (1+\epsilon)\cdot \hat{v}(S_{R}) $.}
	\item \blue{Stop when the current solution does not change.}
\end{enumerate}

 
\blue{Formally, let $k$   index  the iterations, where  the current solution changes in each iteration.  Let 
$S_R^k$ denote the current solution at the start of iteration $k$. The following subroutine implements a single iteration.}

\begin{center}
	\begin{minipage}{\linewidth}
		\begin{algorithm}[H]
			\SetAlgoLined
			\For{$i\in S_R^k$ and $i' \in S_A\backslash S_R^k$} {
			obtain $\hat{v}(S^{k}_{R} -i +i' )$ by solving an LP \;
			}
			let $(c,c')$ be the pair that maximizes $\hat{v}(S^{k}_{R} -i +i')$ over $i\in S_R^k$ and $i' \in S_A\backslash S_R^k$\;
			\uIf{$\hat{v}^*(S^{k}_{R} -c +c') > (1+\epsilon)\cdot \hat{v}(S^{k}_R)$} {
	\blue{set  $S_R^{k+1}\leftarrow S^{k}_{R} -c +c'$ and continue with $k\gets k+1$}  \;						}
			\Else{halt local search and output $S^{k}_{R}$\; }
			
			\caption{Local swap subroutine. \label{algorithm1}}
		\end{algorithm}
	\end{minipage}
\end{center}

In a broad sense, the local swap subroutine does not necessarily \blue{enumerate all pairs $(i, i')$ to search for the optimal $(c, c')$.} 
\blue{A more efficient alternative is terminating each iteration} at the first pair of \blue{$i\in S_R$ and $i'\in S_A\backslash S_R$} that increases the objective by more than  $\epsilon\cdot \hat{v}(S_R^k)$. 
The complete LSLPR algorithm is now as follows: 
\begin{center}
    \begin{minipage}{\linewidth}
    \begin{algorithm}[H]
        \SetAlgoLined
        \KwData{Augmented supply $S_A$, basis supply $S_B$, scenarios $\{\xi_\ell\}_{\ell=1}^N$  and $\epsilon > 0$.}
        \KwResult{Near-optimal $S_R \subset S_A$ and the corresponding trip assignment.}
        
        \textbf{Initialization: } \blue{Set $k = 1$ and
        randomly select $K$ vehicles from $S_A$ as $S^1_{R}$}\;
           \While{$k \leq k_{\max} $} 
       {Run the local swap subroutine in Algorithm 1\;
       }
       Obtain the final trip assignment with  \blue{$S_R=S^{k_{\max}}_R$} using  the greedy algorithm (Theorem~\ref{thm:greedy}).
     \caption{Local search LP-relaxation algorithm for SRAMF \label{alg2}}
    \end{algorithm}
    \end{minipage}
\end{center}

\blue{Algorithm \ref{alg2} obtains the final selection of vehicles $S^{k_{\max}}_R$ where the maximal number of iterations will be derived below.}
In the final step, the algorithm obtains an \textit{integral}  assignment for each scenario \blue{instead of the fractional assignments in $\hat{v}(S_R)$}. To this end, we use the greedy algorithm (Theorem~\ref{thm:greedy}) to select the assignment  for each scenario, which is guaranteed to have value at  least $\frac1p$ times the fractional assignment. 
\blue{In Section \ref{sec:4.1.2}, we first analyze the approximation ratio and then the computational complexity  of LSLPR.}

\subsubsection{Analysis of the LSLPR algorithm.} 
\label{sec:4.1.2}

\def\nx{\bar{\pmb{x}}}

Recall that $S_R$ is the locally optimal solution obtained by our algorithm. Let $S_O$ denote the optimal solution; 
\blue{note that  $S_O$ is a fixed subset that is  only  used in the analysis.} 
Also, let $\x=\langle \x^\xi\rangle$ and $\z=\langle \z^\xi\rangle$ denote the optimal LP solutions to $\hat{v}(S_R)$ and $\hat{v}(S_O)$, respectively. 
 
It will be convenient to consider the overall hypergraph on vertices $S_A\cup S_B \cup \left(\cup_\xi D(\xi)\right)$ and hyperedges $\cup_\xi E(\xi)$. By duplicating hyperedges (if necessary), we may assume  that $E(\xi)$ are disjoint across scenarios $\xi$. 
Recall that $\x^\xi$ (and $\z^\xi$) has a decision variable corresponding to each hyperedge in $E(\xi)$.  
For each demand $d\in \cup_\xi D(\xi)$, let $H_d$ denote the hyperedges incident to it. 
For each vehicle $i\in S_A\cup S_B$ and scenario $\xi$, let $E_{i,\xi}$ denote the hyperedges in $E(\xi)$  containing $i$. So,  $F_i := \cup_{\xi} E_{i,\xi}$ is the set of  hyperedges incident to vehicle $i$.

For any demand $d$, the following lemma sets up a mapping between the hyperedges (incident to $d$) used in the solutions $\x$ and $\z$. 
For the analysis, we add a dummy hyperedge $\bot$ incident to $d$ so that the  assignment constraints in the LP solutions $\x$ and $\z$ are binding at $d$. 
So, $\sum_{e\in H_d} x_e + x_{\bot}= 1$ and $\sum_{f\in H_d} z_f  + z_{\bot}= 1$. Let $H'_d := H_d \cup \{\bot\}$ denote the hyperedges incident to $d$. 

\begin{lemma}\label{lemma4-1}
For any demand $d$, there exists a decomposition mapping $\Delta_d: H'_d \times H'_d  \to \mathbb{R}$ satisfying the following conditions:
    \begin{enumerate}
    \item $\Delta_d(e,f) \geq 0$ for all $e,f \in H_d'$;
    \item $	\sum_{e \in H'_d} \Delta_{d}(e, f) = z_{f}$ for all $f\in H'_d$;
    \item $	\sum_{f \in H'_d} \Delta_{d}(e, f) = x_{e}$ for all  $e\in H'_d$.
\end{enumerate}
\end{lemma}

\begin{figure}[!htb]
    \centering
    \includegraphics[width = 0.5\textwidth]{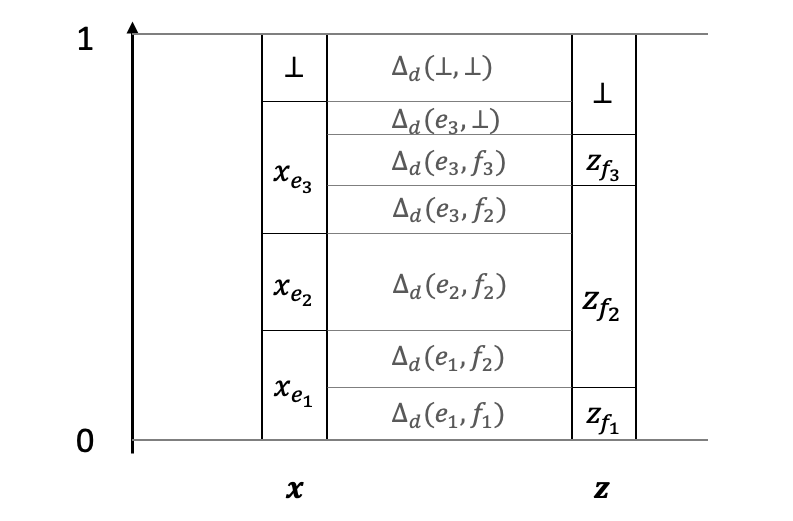}
    \caption{Illustration of mapping $\Delta_d(e,f)$. }
    \label{fig:5}
\end{figure}

Figure \ref{fig:5} illustrates this mapping. 
Appendix \ref{sec:app3} includes the definition of $\Delta_d(e,f)$ and the proof for Lemma \ref{lemma4-1}.
Note that $\sum_{e\in H'_d} \sum_{f\in H'_d} \Delta_d(e,f) = 1$ for any demand $d$. 
For any subset $F\in H_d'$, we use the shorthand $\Delta_d(e,F) := \sum_{f\in F} \Delta_d(e, f)$.

Here is an outline of the remaining analysis. 
 Let ${\cal L}$ denote any bijection between $S_R$ (algorithm's solution) and $S_O$ (optimal solution), consisting of  pairs $(i_1,i_2)$ where $i_1\in S_R$ and $i_2\in S_O$. 
 \blue{We first consider a swap $S_R-i_1+i_2$ where $(i_1,i_2)\in {\cal L}$, and  lower bound  the objective increase.}
Note that the local optimality of $S_R$ implies that the objective increase is at most $\epsilon\cdot \hat{v}(S_R)$. 
Then, we  add the inequalities corresponding to the objective increase for the swaps in ${\cal L}$ and obtain the approximation ratio. 

\noindent \textbf{Analysis of a single swap $(i_1,i_2)$.} \quad Consider any $i_1\in S_R$ and $i_2 \in S_O$. We now lower bound $\hat{v}(S_R - \{i_1\} + \{i_2\}) - \hat{v}(S_R)$. 
Recall that for any subset $S$, $\hat{v}(S) =   \frac{1}{N} \sum_{\xi} \hat{v}(S, \xi)$ where  $\hat{v}(S, \xi)$ is the LP value for scenario $\xi$. So, we have
\begin{align*}
    \hat{v}(S_R - \{i_1\} + \{i_2\}) - \hat{v}(S_R)= \frac1N \sum_\xi \left( \hat{v}(S_R - \{i_1\} + \{i_2\}, \xi) - \hat{v}(S_R, \xi)\right).
\end{align*}

We now focus on a single scenario $\xi$ and lower bound 
$\hat{v}(S_R - \{i_1\} + \{i_2\}, \xi) - \hat{v}(S_R, \xi)$. To this end,  we will define a feasible solution  $\nx^\xi$ for the LP  $\hat{v}(S_R - \{i_1\} + \{i_2\}, \xi)$. 
Recall that $\x^\xi$ denotes the optimal solution for LP $\hat{v}(S_R,\xi)$. 
So, we can then bound: 
\begin{equation}
\hat{v}(S_R - \{i_1\} + \{i_2\}, \xi) - \hat{v}(S_R, \xi)\quad \ge\quad \vv^{\intercal}\nx^\xi - \vv^{\intercal} \x^\xi,
\end{equation}
where $\vv$ is the vector of hyperedge values for $E(\xi)$. 
As we focus on a single scenario $\xi$, we   drop $\xi$ from the notation whenever it is clear.


\ignore{

\begin{figure}[!htb]
	\centering
	\includegraphics[width=0.55\textwidth]{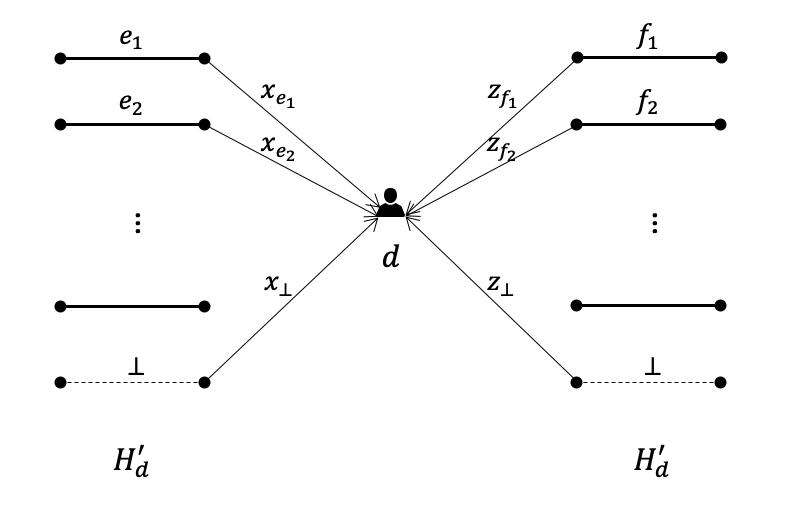}
	\caption{Illustration of single swap $i_1, i_2$ with regard to demand $d$. }
	\label{fig:4}
\end{figure}
}

We are now ready to construct the new fractional assignment $\bar{\pmb{x}}$. Define:
\begin{enumerate}
    \item $\bar{x}_e = 0$ for all $e \in F_{i_1}$. This corresponds to dropping vehicle $i_1$ from $S_R$. 
    \item $\bar{x}_e = z_e$ for all $e \in F_{i_2}$. This corresponds to adding vehicle $i_2$ to $S_R$. 
    \item $\bar{x}_e = x_e - \max_{d\in e} \Delta_d(e, F_{i_2}\cap H_d )$ for all $e \in E(\xi) \setminus F_{i_1}\setminus F_{i_2} $.
\end{enumerate}
If $i_1=i_2$ then we simply drop case~1 above.  
The third case above is needed to make space for the hyperedges incident to the new vehicle $i_2$ (which are increased in case 2). 
The following two lemmas prove the feasibility of this solution $\nx$ and bound its objective value. Below, we assume that $i_1\ne i_2$ (the proof for $i_1=i_2$ is nearly the same, in fact even simpler). 

\begin{lemma}\label{lemma2}
     $\bar{\x}$ is a feasible solution for $\hat{v}(S_R - \{i_1\} + \{i_2\})$. 
\end{lemma}

\proof{Proof for Lemma \ref{lemma2}:} 

We show the feasibility by checking all constraints in eq.\eqref{eq2.5}. Note that $\bar{x}_e=0$ for all  hyperedges $e$ incident to a vehicle in $S_A\setminus \left(S_R - \{i_1\} + \{i_2\} \right)$.

\noindent \textbf{Constraint} $\bar{\x} \geq 0$. It suffices to check this for hyperedges $e\in E \setminus F_{i_1}\setminus F_{i_2}$. Note that
$$\bar{x}_e = x_e - \max_{d\in e} \Delta_d(e, F_{i_2} \cap H_d ) = \min_{d\in e} \left( x_e - \Delta_d(e,  F_{i_2} \cap H_d )\right) \ge 0,$$ where the inequality uses Lemma~\ref{lemma4-1}, i.e., $x_e = \Delta_d(e, H'_d) \geq \Delta_d(e, F_{i_2}\cap H_d )$.

\noindent \textbf{Constraint} \eqref{eq2.5d}: 
By definition of $\bar{\x}$, for any   demand $d$, we have:
\begin{align*}
\sum_{e\in H_d} \bar{x}_e & = \sum_{e\in H_d\cap F_{i_2}} z_e + \sum_{e\in H_d \setminus F_{i_1}\setminus F_{i_2}} \left[x_e - \Delta_d(e, F_{i_2}\cap H_d ) \right]  \\
& \leq \sum_{e\in H_d\cap F_{i_2}} z_e + \sum_{e\in H'_d } \left[x_e - \Delta_d(e, F_{i_2}\cap H_d ) \right]  \\
& = \sum_{e\in H_d\cap F_{i_2}} z_e + \sum_{e\in H'_d } x_e - \sum_{f\in  F_{i_2}\cap H_d} \Delta_d( H'_d , f)\\
& =  \sum_{e\in H'_d} x_e \,\, =\,\, 1.
\end{align*}

\noindent \textbf{Constraint} \eqref{eq2.5e}: 
The augmented vehicle set can be divided into three groups.
\begin{enumerate}
    \item Vehicle $i_1$: $\sum_{e\in F_{i_1}} \bar{x}_e = 0$.
    \item Vehicle $i_2$:  $\sum_{e\in F_{i_2}} \bar{x}_e = \sum_{e\in F_{i_2}} z_e \leq 1$ by definition.
    \item Vehicles $j\neq i_1, i_2$: 
    \blue{$\sum_{e\in F_j} \bar{x}_e \leq  \sum_{e\in F_j} x_e\leq 1 $.}
\end{enumerate}

Therefore, $\bar{\x}$ is a feasible fractional assignment solution. 
\hfill $\square$
\endproof

\begin{lemma} \label{lemma3}
The increase in objective is:
    \begin{align*}
        \sum_{e\in E(\xi)} v_e (\bar{x}^\xi_e - x^\xi_e) \geq 
        \sum_{e\in F_{i_2}\cap E(\xi)} v_e z^\xi_e - \sum_{f\in F_{i_1}\cap E(\xi)} v_f x^\xi_f - \sum_{e\in E(\xi)} v_e  \sum_{d\in e} \Delta_d(e, F_{i_2}\cap H_d ).
    \end{align*}
\end{lemma}

\proof{Proof for Lemma \ref{lemma3}:} 
By definition of $\bar{\x}$, 
\begin{align*}
    \bar{x}_e - x_e =
    \begin{cases}
        z_e   & \text{ if } e \in F_{i_2} \\
         - x_e   & \text{ if } e \in  F_{i_1} \\
        -\max_{d\in e} \Delta_d(e, F_{i_2}\cap H_d )   & \text{ otherwise}
    \end{cases}.
\end{align*}
Note that $x_e = 0$  for all $e\sim S_A\backslash S_R$ in $\hat{v}(S_R)$.  
So we have 
\blue{
\begin{align*}
    \sum_{e\in E(\xi)} v_e (\bar{x}_e - x_e) & \geq  \sum_{e\in F_{i_2}\cap E(\xi)}  v_e z_e -\sum_{f\in F_{i_1}\cap E(\xi)} v_f x_f - \sum_{e \in E(\xi)}v_e \max_{d\in e} \Delta_d(e, F_{i_2} \cap H_d) \\
    & \geq \sum_{e\in F_{i_2}\cap E(\xi)} v_e z_e - \sum_{f\in F_{i_1}\cap E(\xi)} v_f x_f - \sum_{e\in E(\xi)} v_e  \sum_{d\in e} \Delta_d(e, F_{i_2} \cap H_d).
\end{align*}
}
\hfill $\square$
\endproof

Combining   Lemmas~\ref{lemma2} and \ref{lemma3}, and adding over all scenarios $\xi$, we obtain:
\begin{lemma}\label{lem:1-swap}
For any $i_1\in S_R$ and $i_2\in S_O$, we have 
$$ \hat{v}(S_R - \{i_1\} + \{i_2\}) - \hat{v}(S_R) \geq 
        \sum_{e\in F_{i_2}} v_e z_e - \sum_{f\in F_{i_1}} v_f x_f - \sum_{e\in E} v_e  \sum_{d\in e} \Delta_d(e, F_{i_2}\cap H_d).$$
\end{lemma}

\noindent\textbf{Combining all the swaps. } \quad
Recall that  ${\cal L}$ is a  bijection between $S_R$ and $S_O$. Using the local optimality of $S_R$,
\begin{align}
K\epsilon\cdot \hat{v}(S_R) \geq &  \sum_{(i_1, i_2)\in \mathcal{L}  } \left[ \hat{v}(S_R - \{i_1\} + \{i_2\}) - \hat{v}(S_R)\right] \notag \\ 
\ge & \sum_{(i_1, i_2)\in \mathcal{L}  } \left[ \sum_{e\in F_{i_2}} v_e z_e - \sum_{f\in F_{i_1}} v_f x_f - \sum_{e\in E} v_e  \sum_{d\in e} \Delta_d(e, F_{i_2}\cap H_d ) \right] \label{eq:ls-1}\\
    = & \sum_{i_2 \in S_O} \sum_{e\in F_{i_2}} v_e z_e  - \sum_{i_1 \in S_R} \sum_{f\in F_{i_1}} v_f x_f   - \sum_{i_2 \in S_O}\sum_{e\in E} v_e  \sum_{d\in e} \Delta_d(e, F_{i_2}\cap H_d )\notag \\
    \ge & \sum_{i_2 \in S_O} \sum_{e\in F_{i_2}} v_e z_e  - \sum_{i_1 \in S_R} \sum_{f\in F_{i_1}} v_f x_f   -  \sum_{e\in E} v_e  \sum_{d\in e} \Delta_d(e, H_d) \label{eq:ls-2}\\
    \ge & \sum_{i_2 \in S_O} \sum_{e\in F_{i_2}} v_e z_e  - \sum_{i_1 \in S_R} \sum_{f\in F_{i_1}} v_f x_f   -  \sum_{e\in E} v_e  \sum_{d\in e} x_e \label{eq:ls-3}\\
    =& v^T z  - v^T x -    \sum_{e\in E} |\{d\in e\}| v_e x_e \notag \\
    \ge &  v^T z  - v^T x -   (p-1) v^T x \,\,=\,\, v^T z -   p\cdot  v^T x \,\,=\,\, \hat{v}(S_O) - p \cdot \hat{v}(S_R). \label{eq:ls-4}
    \end{align}
Above, \eqref{eq:ls-1} is by Lemma~\ref{lem:1-swap}, \eqref{eq:ls-2} uses that $\{F_{i_2}\}$ are disjoint, \eqref{eq:ls-3} uses Lemma~\ref{lemma4-1}, and the inequality in \eqref{eq:ls-4} uses that each hyperedge has at most $p-1$ demands.

\blue{Setting $\epsilon=\frac{1}{pK^2}$, it follows that $\hat{v}(S_R)\ge \frac{1}{p+o(1)} \cdot \hat{v}(S_O)$. Combined with Theorem~\ref{thm:greedy}, we obtain $v^*(S_R)\ge \frac1p \cdot \hat{v}(S_R)\ge \frac{1}{p^2+o(p)} \cdot \hat{v}(S_O)$.} Hence,
\begin{theorem} \label{theo1}
	The LSLPR algorithm for SRAMF    is a $\frac{1}{p^2}$-approximation algorithm. 
\end{theorem}

\ignore{
\vnote{The additive case seems like a digression- drop it?} \blue{A relevant local-search algorithm that obtains a  tighter approximation ratio of $(1/2 - \epsilon)$ when the hyperedge costs are additive.
For example, for any driver $i \in S_R \cup S_B$, the associated hyperedge $e$'s  value satisfies  $v_e = \sum_{j\in e} \beta_{ij} v_{ij} $. 
The detailed algorithm and the proof of its approximation ratio are included in Appendix \ref{sec:appendc3}. 
}
}

\noindent \textbf{Time complexity of local search. } \quad Note that each iteration (i.e., Algorithm~\ref{algorithm1}) involves considering $K(n_A-K)$ potential swaps; recall that $n_A=|S_A|$. For each swap, we need to evaluate $\hat{v}$, which can be done using any polynomial time LP algorithm such as the ellipsoid method.  So, the time taken in each iteration is polynomial. 

We now bound the  number of local search iterations. 
In each iteration, the objective value increases by a factor at least $1+\epsilon$. So, after $k$ iterations, 
\begin{equation*}
	 \hat{v}(S_{R}^{k+1})  \geq (1+\epsilon)^k
	 \hat{v}(S_R^1).
\end{equation*}

Clearly,  the assignment associated with the initial selected $S_R^0$ has a lower bound $\hat{v}(S_R^0) \geq \frac1N \cdot v_{\min}$ where $v_{\min}=\min_{e : v_e>0} v_e$ is the minimum value over all hyperedges. 
\blue{Recall that $|S_A| = n_A$ and $|S_B|=n_B$. }
The maximum objective of any solution is at most $(n_A+n_B)\cdot v_{\max}$ where $v_{\max}=\max_{e } v_e$ is the maximum value over all hyperedges. Hence,
$$ (n_A+n_B)\cdot v_{\max} \ge \hat{v}(S_{R}^{k+1})  \geq (1+\epsilon)^k\cdot \frac1N v_{\min},$$
	 which implies 
that the maximum number of iterations 
$$k_{\max}\le \log_{1+\epsilon}\left( \frac{N(n_A+n_B) v_{\max}}{v_{\min}} \right) = O\left( \frac1\epsilon \,\log\frac{N(n_A+n_B) v_{\max}}{v_{\min}} \right).$$
Using $\epsilon=\frac{1}{p K^2}$, it follows that the number of iterations is polynomial. 

Finally, the last step of Algorithm~\ref{alg2} just implements the greedy $p$-set packing algorithm (for each scenario), which also takes polynomial time. It follows that 
LSLPR  solves the SRAMF  problem in polynomial time with regard to parameters $p$, $K$, $N$, $|E|$, $n_A$, and $n_B$.

\ignore{
\begin{remark}
With access to the oracle $\mathcal{O}(v_e)$, LSLPR-based algorithm can solve the SRAMF  problem in polynomial time with regard to parameters $p$, $K$, $N$, $|E(\xi)|$, $n_A$, and $n_B$.
\end{remark}

The polynomial run time of the proposed LSLPR algorithm in the end-to-end \blue{stochastic ride-pooling procedure} is guaranteed as follows.
First, in the local swap subroutine, the LP-relaxation of the ride-pooling assignment can be solved by any polynomial-time LP algorithm such as the ellipsoid method. 
Second, the size of mapping for $(i_2, i_1) \in \mathcal{M}$ is upper-bounded by $K^2$.   
Finally, the computation of the exact assignment in the final iteration is equivalent to a $p$-set packing problem with the selected vehicles $S_R$. 
}

\subsection{\blue{Max-Min Online Algorithm for High-Capacity SRAMF}   }  \label{sec:4.2}
The LSLPR algorithm is capable of assigning rides in shared mobility applications with midsize vehicles. 
When the maximal vehicle capacity is large (e.g., the maximum capacity is ten in \citep{alonso2017demand}), $\frac{1}{p^2}$-approximation ratio becomes unacceptable in the worst case. 
We propose an alternative method for high-capacity SRAMF. 
\blue{The main idea of the Max-Min Online (MMO) algorithm is to use LP-duality to reformulate $\hat{v}$ as a \textit{covering} linear program. 
Then, we use an existing framework for max-min optimization from~\cite{feige2007robust}.
This framework requires two technical properties (monotonicity and online competitiveness), which we show are satisfied for SRAMF.} 
We will prove that  the MMO algorithm obtains an approximation ratio of  $(1-\frac{1}{e}) \frac{1}{ 2p \ln p }$.

Using LP-duality and the definition of $\hat{v}(S_R)$ (see the derivation in Appendix~\ref{sec:app4}), we can reformulate: 
\begin{align} \label{eq8}
    \hat{v}(S_R) \,\,= &\,\,\minimize_u \sum_{\xi} \sum_{ g \in G  }    u_{g, \xi } \\ 
 &      s.t.\,\, \sum_{g\in e}  u_{g, \xi }  \geq \frac{v_e}{N} , & \forall e \in F_{i,\xi},\quad \forall \xi, \quad  \forall i\in S_R \cup S_B  \nonumber \\
    & \qquad \pmb{u} \geq 0. &  \nonumber
\end{align}
\blue{Here, $G=S_A\cup S_B\cup \left( \cup_\xi D(\xi)\right)$ is a combined groundset consisting of all vehicles and demands from all scenarios. 
For any vehicle $i$ and scenario $\xi$, set $F_{i,\xi}\subseteq E(\xi)$ denotes all the hyperedges incident to $i$ in scenario $\xi$.  }

We can scale the covering constraints to normalize the right-hand-side to $1$ and rewrite the constraints as $\sum_{g\in e} \frac{N}{v_e} u_{g, \xi }\geq 1$. \blue{Note that the \textit{row sparsity} of this constraint matrix (i.e., the maximum number of non-zero entries in any constraint) is $\max_{e\in E} |e| =p$ and $v_e >0$ for all hyperedges.}
Let $\pmb{c}_e$ be the row of constraint coefficients for any hyperedge $e\in E=\cup_\xi E(\xi)$, i.e., 
\begin{align*}
    c_e(g,\xi) = \begin{cases}
        \frac{N}{v_e} & \text{ if } g\in e \mbox{ and } e\in E(\xi) \\
        0   & \text{ otherwise}
    \end{cases}.
\end{align*}

Then, the SRAMF  problem with fractional assignments $\max_{S_R\subseteq S_A : |S_R|\le K} \hat{v}(S_R)$ can be treated as the following max-min problem:
\begin{align} \label{eq:maxmin}
    \max_{S_R\subseteq S_A : |S_R|\le K} \quad \min_{u} \{ \mathbf{1}^{\intercal} \uu \,\, | \,\,
    \pmb{c}_e^{\intercal} \uu  \geq 1, \,\, \forall e\in F_i, \,\, \forall i \in  S_R \cup S_B;\,\, \uu \geq 0 \},
\end{align}
where $F_i=\cup_\xi F_{i,\xi}$ for each vehicle $i$. 

The main result is:
\begin{theorem} \label{thm3}
 	 There is a $\frac{e-1}{(2e+o(1)) \ln p}$-approximation algorithm for \eqref{eq:maxmin}.
 \end{theorem}
 
 Before proving this result, we introduce two important properties.

\begin{definition}{(Competitive online property)}\label{defn:online}
An $\alpha$-competitive online algorithm for the covering problem eq.\eqref{eq8} takes as input 
any sequence $(i_1, i_2, \dots, i_t, \dots)$ of vehicles from $S_A$ and maintains a non-decreasing solution $\uu$ such that the following hold for all steps $t$. 
\begin{itemize}
\item $\uu$ satisfies constraints $\pmb{c}_e^{\intercal} \uu \geq 1$ for $e\in F_i$, for all vehicles $i\in \{ i_1, i_2, \dots, i_t\}$, and  
\item $\uu$  is an $\alpha$-approximate solution, i.e.,  the objective $\pmb{1}^{\intercal} \uu \,\leq \, \alpha \cdot \hat{v}(\{ i_1, i_2, \dots, i_t\})$.  
\end{itemize}
\end{definition}

Note that the online algorithm may only increase  variables $\uu$ in each  step $t$.

\begin{definition}{(Monotone property)}
For any $\uu\ge 0$ and $S\subseteq S_A$, let 
$$Aug^*(S|\uu) :=  \{\min_{\pmb{w} \geq 0} \pmb{1}^{\intercal} \pmb{w}: \pmb{c}_e^{\intercal}(\uu +\pmb{w}) \geq 1,\, \forall e\in F_i, \, \forall i\in S\cup S_B  \}.$$  
The covering problem eq.\eqref{eq8} is said to be monotone if for any $\uu \geq \uu'\ge 0$ (coordinate wise) and   any $S  \subseteq S_A$, $Aug^*(S |\uu) \leq Aug^*(S |\uu')$.
\end{definition}

These properties were used by \cite{feige2007robust} to show the following result. 
\begin{theorem}{\citep{feige2007robust}} \label{theorem4}
If the covering problem \eqref{eq8} satisfies the monotone and $\alpha$-competitive online properties, there is a $\frac{e-1}{e\cdot \alpha}$-approximation for the max-min problem in eq.\eqref{eq:maxmin}.
\end{theorem}

Our max-min problem indeed satisfies both these properties. 
\begin{lemma} \label{lem:online}	The covering problem \eqref{eq8}  has an 
     $\alpha = O(\ln p)$ competitive online algorithm. Moreover, when $p$ is large, the factor $\alpha=(2+o(1))\ln p$.
\end{lemma}
\proof{Proof:} Recall that \eqref{eq8}  is a covering LP with row-sparsity $p$. Moreover, in the online setting, constraints to \eqref{eq8} arrive over time. So, this is an instance of online covering LPs, for which an $O(\ln p)$-competitive algorithm is known \citep{gupta2014approximating}. See also \citep{buchbinder2014online} for a simpler proof. Moreover, one can optimize the constant factor in \citep{buchbinder2014online} to get $\alpha = (2+o(1))\ln p$. We note that these prior papers 
work with the online model where only one covering constraint arrives in each step. Although \blue{Lemma~\ref{lem:online}} involves multiple covering constraints $F_i$ arriving in each step, this can be easily reduced to the previous setting: just introduce the constraints in $F_i$ one-by-one in any order, and then the algorithms in \citep{gupta2014approximating, buchbinder2014online} can be used directly. 
 \hfill $\square$
\endproof 

\begin{lemma} \label{lem:mon}
	The covering problem \eqref{eq8} is monotone. 
\end{lemma}
\proof{Proof:}
Consider any $\uu \geq \uu'\ge 0$ and any $S  \subset S_A$. Let 
$\pmb{w}'\ge 0$ denote an optimal solution to $Aug^*(S_R|\uu')$. As all constraint-coefficients $\pmb{c}_e \geq 0$, it follows that  $\pmb{c}_e^{\intercal}(\uu + \pmb{w}') \geq \pmb{c}_e^{\intercal}(\uu' + \pmb{w}') \geq 1$ for all $e\in F_i$ and $i\in S\cup S_B$. Hence, $\pmb{w}'$ is also a feasible for the constraints in  $Aug^*(S|\uu)$. Therefore, $Aug^*(S|\uu) \leq \pmb{1}^{\intercal} \pmb{w}' = Aug^*(S|\uu')$, which proves the monotonicity.
 \hfill $\square$
\endproof

Combining Lemmas~\ref{lem:online} and \ref{lem:mon} with Theorem~\ref{theorem4},  we obtain Theorem~\ref{thm3}. 
We note that our $\Omega(\frac{1}{\ln p})$ approximation ratio is nearly the best possible for the max-min problem \eqref{eq:maxmin}, as the problem is hard to approximate to a factor better than $O(\frac{\ln\ln p}{\ln p})$; see \cite{feige2007robust}. 


We now describe the complete algorithm below in the context of SRAMF. This is a combination of the max-min algorithm from \cite{feige2007robust} and the online LP algorithm from \cite{buchbinder2014online}. 
\def\von{\hat{v}_{ON}}
For any ordered subset $S$ of vehicles, let $\von(S)$ denote the objective value of the online algorithm for \eqref{eq8} after adding constraints corresponding to the vehicles in $S$ (in that order). Algorithm~\ref{alg3} describes the updates performed by the online algorithm when a vehicle $i$ is added. 


\begin{center}
	\begin{minipage}{\linewidth}
		\begin{algorithm}[H]
			\SetAlgoLined
			For a given $i\in S_A\cup S_B$, perform the following updates\;
			\For{$e \in F_i=\bigcup_\xi F_{i,\xi}$}{
			let $\{u_{g, \xi}^{-}\}_{g\in e}$ be the   values of variables in hyperedge $e$  and $\Gamma_e^{-}= \sum_{g\in e} u^{-}_{g,\xi}$\;
			
				\If{ $\Gamma_e^{-} < \frac{v_e}{N} $ }{ 
    \begin{align*}
			    \mbox{update } \, u_{g,\xi} \leftarrow  \left( u_{g,\xi}^- + \frac{v_e}{N} \delta\right)\cdot \frac{1 + |e|\cdot \delta}{\frac{N}{v_e} \Gamma^-_e + |e|\cdot \delta} - \frac{v_e}{N} \delta, \quad \mbox{ for all }g\in e.  
			\end{align*}
			 }
				}

			\caption{Updating subroutine in Max-Min Online algorithm } \label{alg3}
		\end{algorithm}
	\end{minipage}
\end{center}

\proof{Proof for the updating subroutine in MMO algorithm: }
 Consider the updates when vehicle $i$ is added. Consider any scenario $\xi$ and hyperedge $e\in F_{i,\xi}$: the corresponding  covering constraint is $c_e^T u = \frac{N}{v_e}\sum_{g\in e} u_{g,\xi} \ge 1$. 
 Let $\tau$ be a continuous variable denoting time. 
 The online LP algorithm in \citep{buchbinder2014online} 
 \blue{raises variables $u_{g,\xi}$} in a continuous manner as follows:
  \begin{align}\label{eq:online-update}
 	\frac{\partial u_{g, \xi }}{\partial \tau } = \frac{N}{v_e} u_{g, \xi } + \delta, \qquad \forall g\in e,  
\end{align}	
until the constraint is satisfied.  
Letting $\Gamma_e=\sum_{g\in e}u_{g,\xi}$, we have 
$$ 	\frac{\partial \Gamma_e}{\partial \tau } = \frac{N}{v_e}\sum_{g\in e} u_{g, \xi } + |e|\cdot \delta = \frac{N}{v_e} _e + |e|\cdot \delta.$$

By integrating, it follows that the duration of this update is 
$$T=  \int_{\Gamma = \Gamma_e^-}^{\Gamma_e^+} \frac{ \partial \Gamma_e }{ \frac{N}{v_e} \Gamma_e + |e|\cdot \delta}  = \frac{v_e}{N}\cdot  \ln\left(\frac{\frac{N}{v_e} \Gamma^+_e + |e|\cdot \delta}{\frac{N}{v_e} \Gamma^-_e + |e|\cdot \delta} \right)=\frac{v_e}{N}\cdot  \ln\left(\frac{1 + |e|\cdot \delta}{\frac{N}{v_e} \Gamma^-_e + |e|\cdot \delta} \right).$$

Above $\Gamma_e^-$ and $\Gamma_e^+$ denote the values of $\Gamma_e$ at the start and end of this update step; note that $\Gamma^+_e= v_e / N$ as the updates stop as soon as the constraint is satisfied. 
For each $g\in e$, using \eqref{eq:online-update},
$$T= \int_{\tau=0}^{T} \frac{ \partial u_{g,\xi}  }{ \frac{N}{v_e} u_{g,\xi} + \delta} = \frac{v_e}{N}\cdot  \ln\left(\frac{\frac{N}{v_e} u_{g,\xi}^+ + \delta}{\frac{N}{v_e} u_{g,\xi}^- + \delta} \right).$$

Again, $u_{g,\xi}^-$ and $u_{g,\xi}^+$ denote the values of $u_{g,\xi}$ at the start and end of this update step. 
Combined with the above value for $T$, we get a closed-form expression for the new variable values:
$$ \frac{N}{v_e} u_{g,\xi}^+ + \delta =\left( \frac{N}{v_e} u_{g,\xi}^- + \delta\right)\cdot \frac{1 + |e|\cdot \delta}{\frac{N}{v_e} \Gamma^-_e + |e|\cdot \delta}, \qquad \forall g\in e. $$
 \hfill $\square$
\endproof 

The complete MMO algorithm is described in Algorithm \ref{alg4}:

\begin{center}
	\begin{minipage}{\linewidth}
		\begin{algorithm}[H]
			\SetAlgoLined
			\KwData{Augmented supply $S_A$, basis supply $S_B$,  hypergraph $G$ with $E(\xi)$, and $\epsilon > 0$.}
            \KwResult{Near-optimal $S_R \subset S_A$ and the corresponding trip assignment.}
			Initialization: $S_R \gets  \varnothing $ and dual variables $\mathbf{u}\gets 0$\;
			For each vehicle in $S_B$ (in any order), run   Algorithm \ref{alg3} to obtain $\von(S_B)$

			\For{$k= 1, \dots ,K$ }
			{   
				\For{ $i \in S_A\backslash S_R$  }
				{
					Run the updating subroutine in Algorithm \ref{alg3} and obtain $ \von(S_B + S_R  + \{i\} )$.
				}
				$i^* = \arg \max_{i\in S_A \backslash S_R } \von(S_B +S_R + \{i\} )  $ \;
				$S_R \gets S_R + \{ i^* \}$\;	
			}	
			
			\caption{ Max-Min online algorithm for SRAMF  } \label{alg4}
		\end{algorithm}
	\end{minipage}
\end{center}

\subsection{\blue{Extensions to SRAMF under Partition Constraints}} \label{sec4.3}

We now consider  a more general setting where  the augmented set $S_A$ is partitioned into $M$ subsets $S_A(1),\cdots S_A(m)$ and the platform requires $K_m$ vehicles from each subset. 
For example, there are $M$ types of vehicles so the cardinality constraint is further specified for each type as $\sum_{i\in S_A(m)} y_i \leq K_m$ for all $m\in [M]$. 
Alternatively, in the mixed autonomy example,  there are $M$ subregions of AV zones and the requirement is proportional to the demand density in each subregion. Instead of \eqref{eq1}, we now want to solve:
\begin{subequations}\label{eq4.3}
    \begin{align}
        &\maximize_{y}  \mathbf{E} [Q(y, \xi)]   \tag{\ref{eq4.3}} \\
        s.t. & \sum_{i\in S_{A}(m)} y_i \leq K_m  &  \forall m \in [M]   \\
        & y_i \in \{ 0,1 \} & \forall i \in S_A.   
    \end{align}
\end{subequations}

\blue{We can extend our result to obtain:}
\begin{theorem}
     The MMO algorithm is a $\frac{1}{(4+o(1))p\log p}$-approximation algorithm for SRAMF  with partition constraints. 
\end{theorem}
The proof is identical to that of Theorem~\ref{thm3}. The only difference is the use of the following result for max-min covering under a partition constraint (instead of Theorem~\ref{theorem4}, which only holds for a cardinality constraint). 


\begin{theorem}{\citep{gupta2015robust}} \label{lemma6}
If the covering problem \eqref{eq8} satisfies the monotone and $\alpha$-competitive online properties, there is a $\frac{1}{2 \alpha}$-approximation for the max-min problem with a partition (or matroid) constraint.
\end{theorem}

\section{\blue{Numerical Experiments in Mixed Autonomy Traffic} } \label{sec:5}

\subsection{\blue{Data Description and Experiment Setup} }  \label{sec:5.1}
\blue{We evaluate the performance of the proposed algorithms through two settings of experiments. 
\begin{enumerate}
    \item \textit{Setting 1} represents the high-capacity SRAMF in which the AV fleet belongs to a public transit operator.
    On-demand AVs are treated as a complementary mode to the existing transit system. 
    A relatively small number of high-capacity microtransit vehicles serve the demand in AV zones, most of which are first- or last-mile connection trips.
    \item \textit{Setting 2} represents the mid-capacity SRAMF in which a private platform operates an electric AV fleet alongside conventional gas vehicles.
   We assume that the electric AVs will return to a charging station to recalibrate and top off their battery after completing a trip. However, the exact station should be chosen carefully so that the AV can serve future demand with minimal pickup times.
\end{enumerate}
}

We test the performance of these approximation algorithms in an on-demand mobility simulator.
The simulator operates a mixed autonomy fleet and  \blue{integrates the batch-to-batch procedure in \cite{alonso2017demand} for ride-pooling and a demand forecast model to evaluate the performance with a rolling horizon. }

\subsubsection{Data description and preprocess. }
The main components of the data input are:
\begin{enumerate}
    \item Road network: 
    \blue{The road network was constructed from OpenStreetMap data in New York City (NYC) and traveling times are computed using the average speed profiles from historical data \citep{sundt2020heuristics}.}
    \blue{In Setting 1, two AV zones are located in the high-density areas (Figure \ref{fig:7-1}).
    AVs only operate within these AV zones.}
    These zones were chosen based on high-demand areas that have been proposed for pedestrianization or could feasibly be closed off to most vehicles other than AV shuttles. 
    \blue{In Setting 2, this zone restriction is lifted and AVs can operate over the entire area of Manhattan.}

    \begin{figure}[!htb]
        \centering
        \subfloat[\blue{CVs' and AVs' initial locations in Setting 1} \label{fig:7-1}]{\includegraphics[ height = 3.5in]{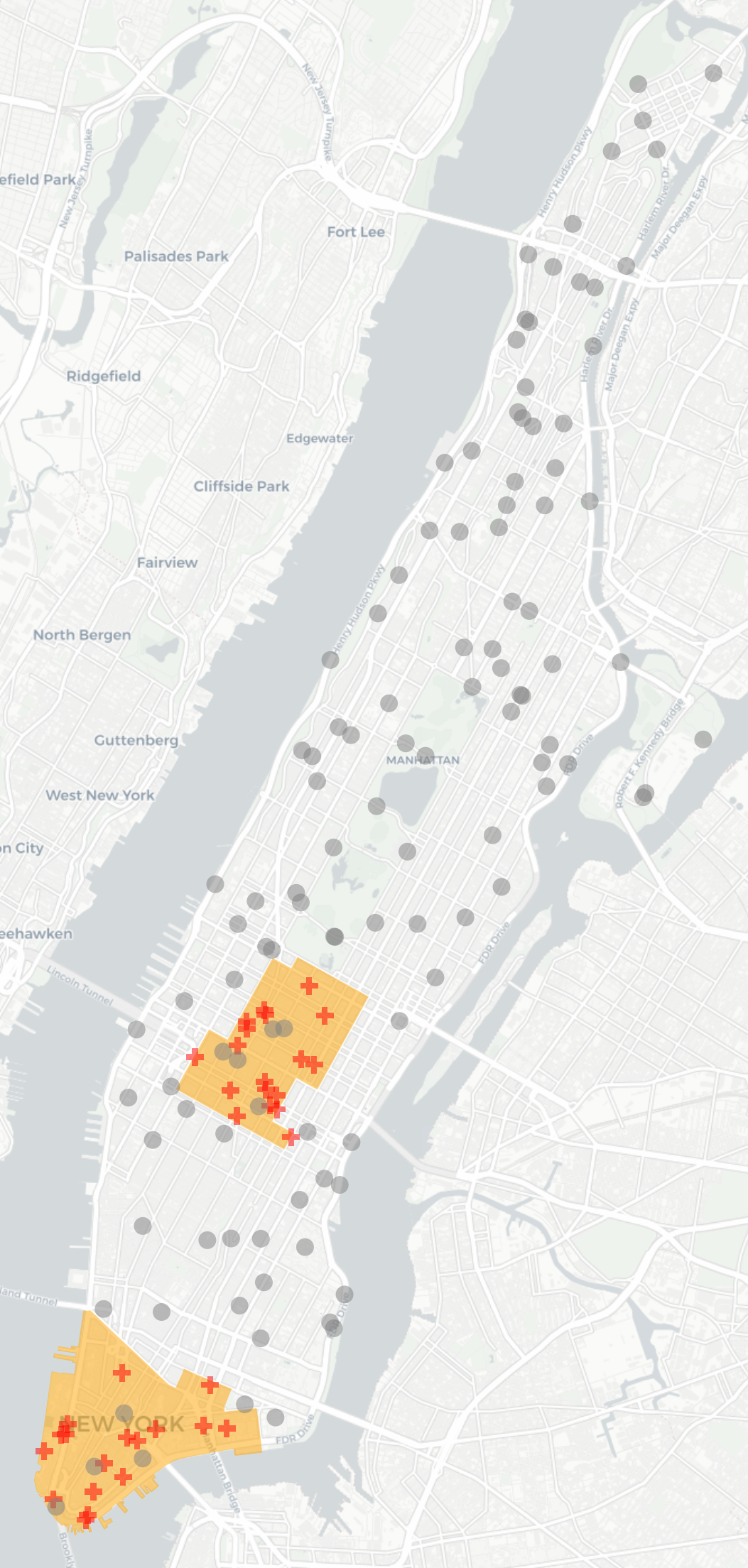}}
        \qquad 
        \subfloat[\blue{Request pickup locations in a sampled experiment}  \label{fig:7-2} ]{\includegraphics[ height = 3.5in]{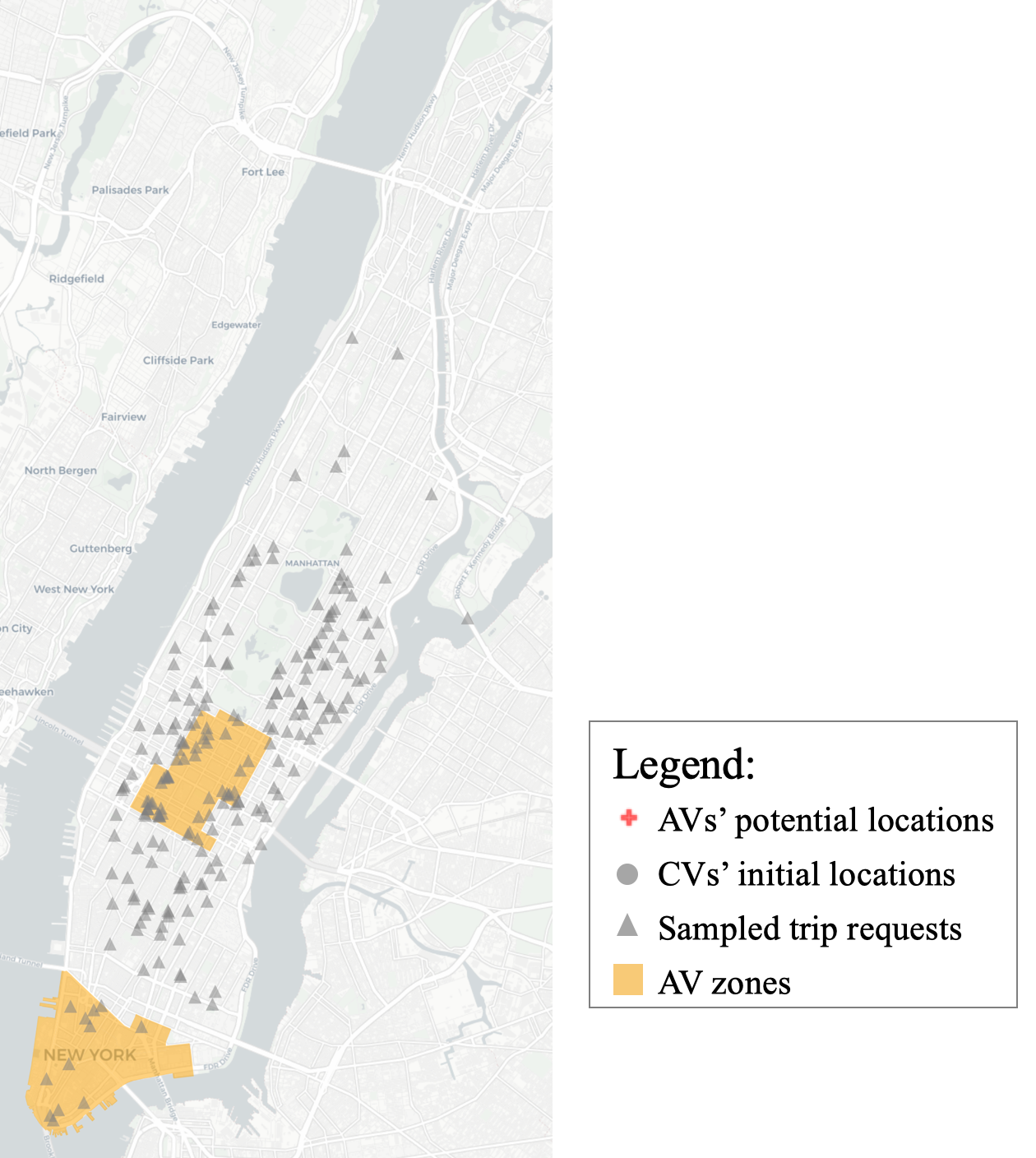}} \\
        \subfloat[Sampled hyperedge value distributions in shareability graphs \label{fig:7-3}  ]{\includegraphics[ width = 0.65\textwidth]{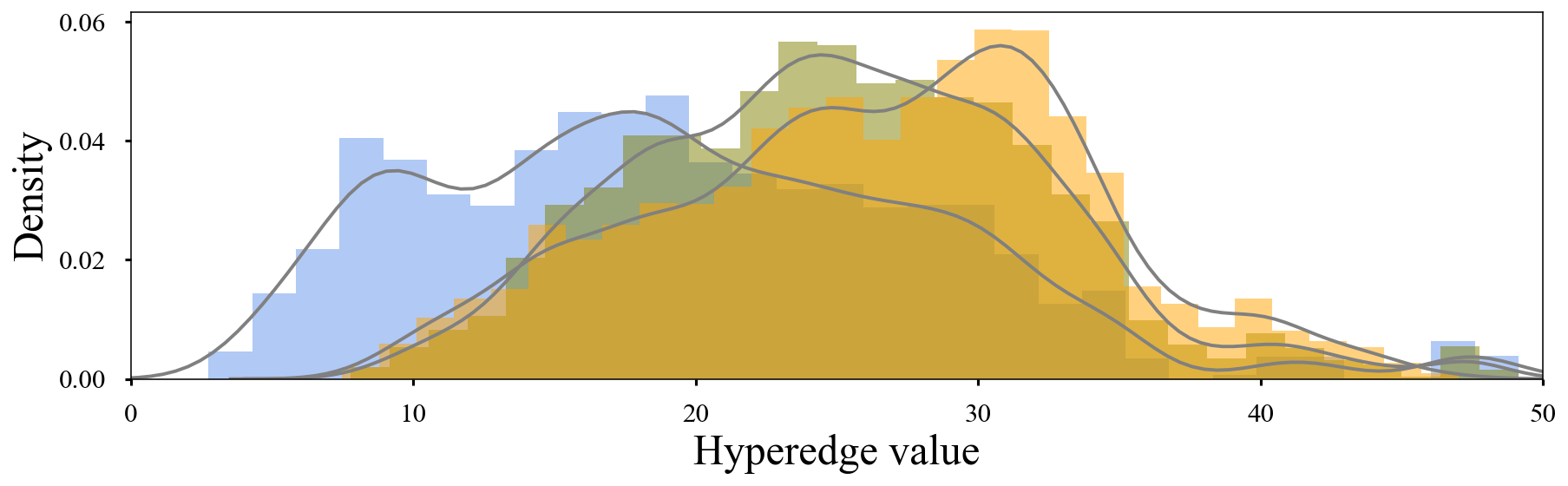}} 
        \caption{\blue{High-capacity mixed autonomy traffic experiment in Manhattan, NYC. } }
        \label{fig:7}
    \end{figure}
    
    \item Supply: 
   \blue{The basis set $S_B$ represents the CVs with a fixed capacity of two that provides the conventional ride-hailing service.}
    
    \blue{
    \begin{itemize}
        \item  
        In the high-capacity SRAMF (Setting 1 in Figure \ref{fig:7}), the augmented set $S_A$ represents the initial locations of automated shuttles fleets that can be redistributed for the regular mobility service.  
         Each shuttle is of capacity up to ten and can only operate within AV zones.
        \item  In the mid-capacity SRAMF (Setting 2 in Figure \ref{fig:r7}), the augmented set $S_A$ represents a set of locations to preposition taxi-like AVs of capacity three.
        These locations are sampled from charging stations in NYC from the NREL Alternative Fueling Stations data \cite{nrel}. The algorithm will choose among $S_A$ to position idle AVs for charging before the next assignment. 
    \end{itemize}
    }

    \begin{figure}[!htb]
    \centering
    \subfloat[ \blue{Public EV charging stations in Manhattan, NYC} \label{fig:r7-1}]{\includegraphics[ height = 3.5in]{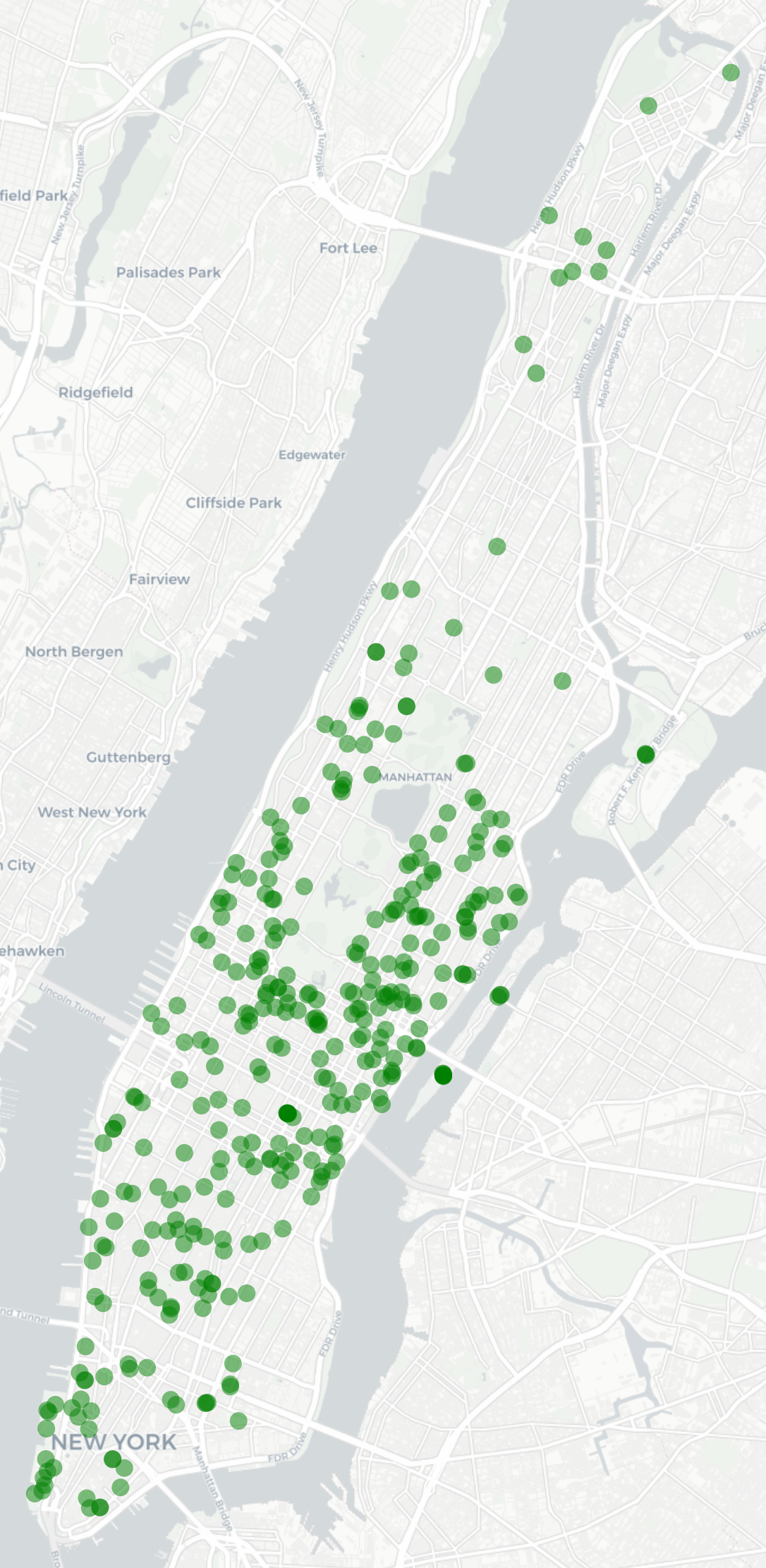}}
    \quad 
    \subfloat[\blue{$S_A$ and $S_B$ for Setting 2 in a given scenario}\label{fig:r7-2} ]{\includegraphics[ height=3.5in]{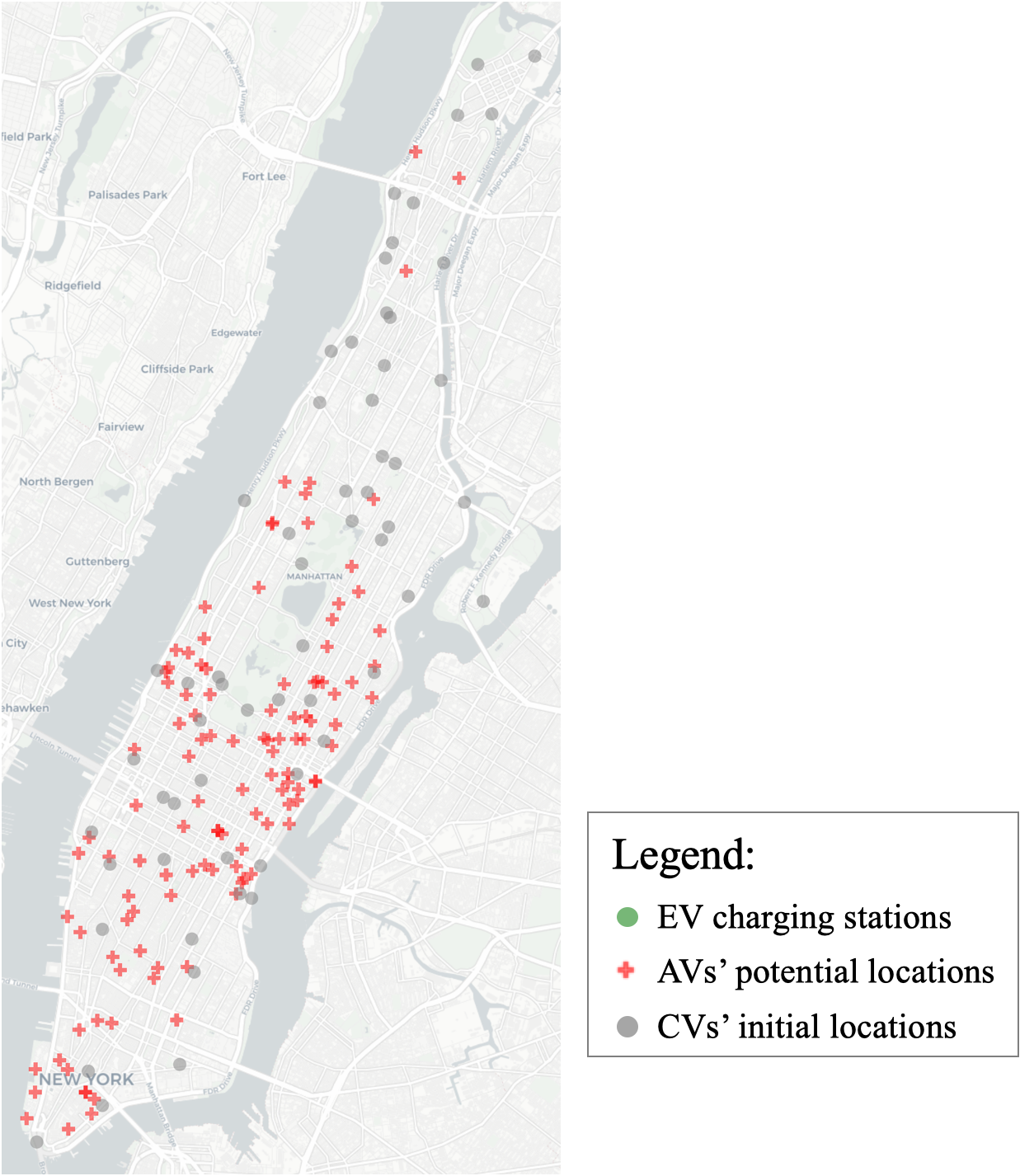}}
    \caption{\blue{Mid-capacity mixed autonomy traffic experiment in Manhattan, NYC. 
    $S_A$ in Figure \ref{fig:r7-2} are randomly chosen from the current EV charging stations \citep{nrel} in Figure \ref{fig:r7-1}. } }\label{fig:r7}
    \end{figure}

    \item Demand: The trip requests are sampled from the NYC Taxi and Limousine Commission trip data \blue{which includes the origin-destination, number of passengers, trip time, and fares} \citep{noauthor_tlc_nodate} (Figure \ref{fig:7-2}). 
    
    \item Hyperedge values: The hyperedge cost follows eq.\eqref{eq0}.
    Each trip's pickup time is computed from the shortest path connecting trip requests on a road network.  
    Customer's preference over CVs and AVs is randomly generated with $v_e > 0$ for all hyperedges $e$.
    $N$ scenarios are independent and identically distributed and used to compare approximation algorithms and the benchmark model (Figure \ref{fig:7-3}).
    \item \blue{Time intervals: 
    We consider the rolling-horizon policy with a fixed interval of 15 minutes. 
    This interval allows for repositioning of selected idle vehicles. In addition, many charging stations in NYC are equipped with high-voltage chargers or superchargers, which can restore up to 200 miles in 15 minutes of charging \citep{xiong2017optimal}. }
\end{enumerate}


\subsubsection{ \blue{Assessment of algorithms in the mixed-autonomy simulator.}}


\blue{In both settings, the objective is to choose a subset of $S_A$ to reposition these vehicles between trips and maximize the total value of assignments through the rolling horizon, including the values of fulfilled demand, trip costs due to the increasing pickup times, and customers' preference over the vehicle type in eq.\eqref{eq0}. 
These candidates in $S_A$ can be the real-time GPS locations of vehicles (Example 1 or Example 2) or transit stops and parking areas for AVs (Example 3).}
The solution algorithms select vehicles in each batch at the beginning of each period and determine the trip assignment after all demands are revealed.
\blue{During preprocessing, the demand forecast model generates $N$ scenarios and constructs the shareability graph using the process in Appendix \ref{sec:app2}.
It generates all hyperedges for the sampled demand as well as the hyperedge values associated with all available vehicles. 
The hyperedge numbers and other results reported in figures and tables for experiments in this paper are averaged out over all samples.}

\blue{All computation times are reported from performance on a server with an 18-core 3.1 GHz processor and 192GB RAM. 
The benchmark model uses a state-of-the-art IP solver (Gurobi 9.1). Setting 1 reports 4-core performance due to the smaller problem size, while Setting 2 is large enough to require the use of all CPU cores.
Solving the SRAMF problem in \eqref{eq1} to its optimality by an exact method is impractical considering the massive size of the potential trips. 
The number of variables in the IP is equal to the product of the number of hyperedges and the sample size, both of which can grow exponentially in real-world applications.
The computational time limit is set as six hours per instance.
Although the proposed approximation algorithms can handle more extensive networks, our numerical experiments downsample from the taxicab data to keep solvable scales for the benchmark IP solver.}

The performance of the proposed approximation algorithm is evaluated under different supply and demand distributions.
As the demand profiles are sampled from the taxicab trip data, these algorithms do not depend on any distributional assumption and 
can be connected to more sophisticated demand forecast models \citep{geng2019spatiotemporal}.
The system is tested in both a relatively balanced demand scenario as well as a massive under-supply scenario, with mean numbers of demand across scenarios ranging from $250$, to $4000$, respectively, which are considerably large in a \textit{stochastic} setting. 





\begin{table}[!htb]
\caption{\blue{Parameters in numerical experiments}} \label{tab1}
\resizebox{\textwidth}{!}{
\blue{
\begin{tabular}{c|c|c|c|c|c|c|c|c}
\toprule
\multirow{2}{*}{Setting} &
\multicolumn{3}{c|}{\begin{tabular}[c]{@{}c@{}}AVs \\ (augmented set)\end{tabular}} & \multicolumn{3}{c|}{\begin{tabular}[c]{@{}c@{}}CVs \\ (basis set)\end{tabular}}                                              & \multirow{2}{*}{\begin{tabular}[c]{@{}c@{}}Ratio of demand to \\ supply of vehicles\end{tabular}} & \multirow{2}{*}{\begin{tabular}[c]{@{}c@{}}Number of \\  sampled scenarios\end{tabular}} \\ \cline{2-7}
& Capacity      & \begin{tabular}[c]{@{}c@{}}Setup \\ cost\end{tabular}      & K     & Capacity & \begin{tabular}[c]{@{}c@{}}Setup \\ cost\end{tabular} & \begin{tabular}[c]{@{}c@{}}Vehicle\\ number\end{tabular} &  &                                                                                           \\ \hline
Setting 1  & 5-10             &   1                                                         & 5     & 3       & 0                                                     & 115                                                      &          1.7-2 & 50                \\ \hline
Setting 2 &3             &   1                                                         & 30-60     & 2       & 0                                                     & 60                                                     &  35-45        & 20-30                                                                                            \\ \bottomrule
\end{tabular}
}
}
\end{table}

\subsection{\blue{Numerical Results for High-Capacity SRAMF}}

\begin{figure}[!htb]
    \centering
    \subfloat[\blue{CV routes}]{\includegraphics[width=0.3\textwidth ]{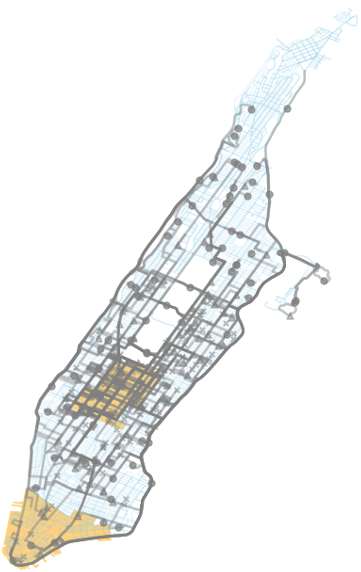} }
    \quad 
    \subfloat[\blue{AV routes}]{\includegraphics[width=0.51\textwidth]{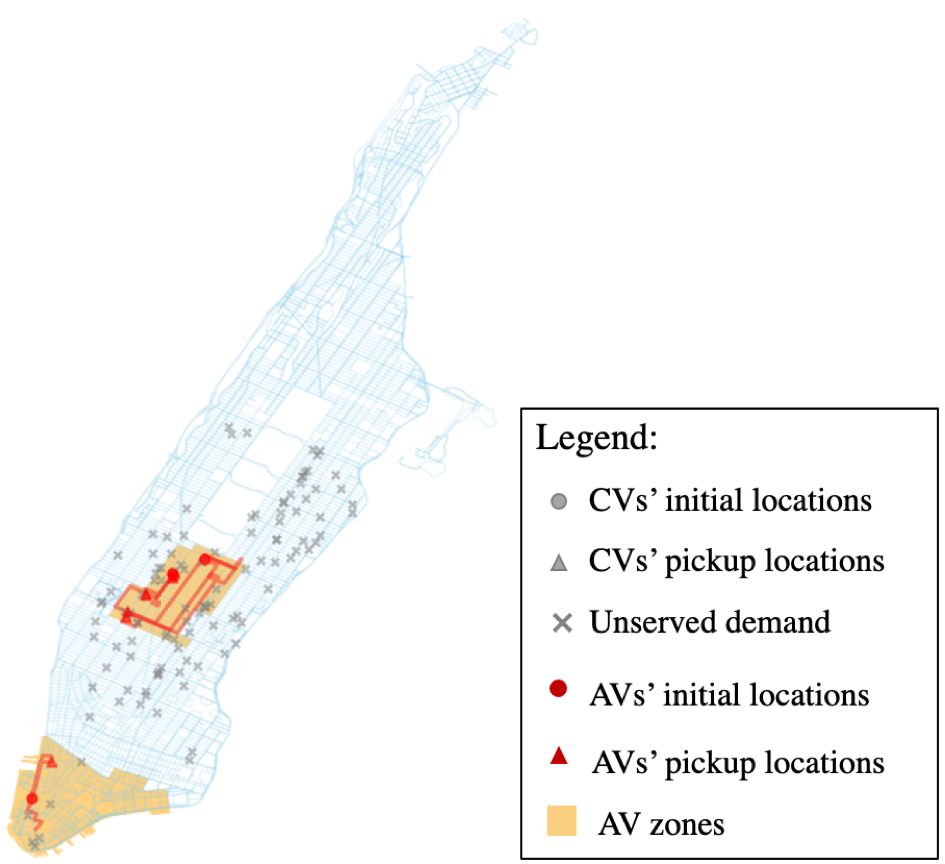} }

    \caption{\blue{Optimal trip assignment and routes in mixed autonomy, high-capacity SRAMF.} }
    \label{fig:map}
\end{figure}

\blue{In Setting 1, the proposed algorithm computes the near-optimal solutions for the mixed-autonomy fleet, including selecting AVs and routing vehicles for each sampled demand profile. }
Figure \ref{fig:map} shows how those chosen AVs (red) and CVs (grey) service trip requests in the face of uncertainty where AVs only operate within the AV zones.
The algorithm allocates one AV in the lower AV zone and four AVs in the upper AV zone on the map, which matches the density of demand in Figure \ref{fig:7-2}.

The remainder of this section compares two approximation algorithms with a benchmark IP method with regard to the total run times and optimality gap. 
The computation time comparison is to validate the polynomial-time reduction in terms of the network size, and the optimality-gap comparison is to examine those proved approximation ratios. 

\subsubsection{Computation times. }
The reported computation time includes solving for the near-optimal selection of vehicles and exact assignment in each scenario. 
\blue{Since the same hypergraph is generated beforehand and used throughout all algorithms, we do not report them in this comparison. }
Two parameters that determine the size of the shareability graph is $|S_A|$ (the size of the augmented set) and the number of hyperedges (the number of decision variables in each scenario). 
\blue{Table \ref{table2}} shows how the total run time grows with the increasing size of the hypergraph.
The computation time of LSLPR and MMO algorithms are shown in \blue{Table \ref{table2}}. 

\blue{Note that a significant difference between the proposed algorithms and the exact solver is that these tailored algorithms add vehicles from the set $S_A$ sequentially. 
A parallel-computing scheme can significantly reduce the total run time of approximation algorithms because the evaluation of each scenario can be carried out simultaneously. 
The exact approach (viz., IP solver) cannot reduce the total run time since it must solve the stochastic program \eqref{eq1} across all samples.
The results report the LSLPR and MMO algorithms under a finite-computing-resource (at most eight threads per time) and an infinite-computing-resource setting. 
The infinite computing resource means that we can evaluate all pairs of vehicles in the active set in LSLPR or the dual variables for all hyperedges in MMO simultaneously. 
The reported run time is the maximal runtime per iteration. In our experiments, we were able to achieve times close to the reported MMO infinite computing setting, because that can be achieved by evaluating all drivers in $S_A \setminus S_R$ in parallel. 
The LSLRP limit is much harder to achieve as the number of potential swaps is combinatorial, however, the computation time will always improve when additional resources are used.
The number of hyperedges is the maximum number of hyperedges per scenario.}

\subsubsection{Optimality gaps. }
The optimality gaps of algorithms are compared with the IP benchmark model. 
Let the objective of the IP solver be $OPT$ and the approximation algorithm's solution be $ALG$. 
The optimality gap is measured by $(OPT-ALG)/OPT$.
\blue{Table \ref{table2}} shows how the optimality gaps grow with the increasing network size under various supply-demand ratios.



\begin{table}[!htb]
\begin{center}
\caption{\blue{Summary of Numerical Results for High-Capacity SRAMF}} \label{table2}
\blue{{\scriptsize 
\begin{tabular}{c|c|c|c|c|c|ccc|ccc}
\toprule
\multicolumn{1}{l|}{\multirow{2}{*}{\begin{tabular}[c]{@{}l@{}}Demand-\\ supply \\ ratio\end{tabular}}} & \multirow{2}{*}{$S_A$} & \multirow{2}{*}{$S_B$} & \multicolumn{1}{l|}{\multirow{2}{*}{\begin{tabular}[c]{@{}l@{}}Total \# IP\\ variables \\ (hyperedges \\ $\times$ samples)\end{tabular}}} & \multicolumn{1}{l|}{\multirow{2}{*}{\begin{tabular}[c]{@{}l@{}}\# of\\ samples\end{tabular}}} & Benchmark                                                                              & \multicolumn{3}{c|}{LSLPR}                                                                                                                                                                                                                                                                      & \multicolumn{3}{c}{MMO}                                                                                                                                                                                                                                                 \\ \cline{6-12} 
\multicolumn{1}{l|}{}                                                                                   &                       &                       & \multicolumn{1}{l|}{}                                                                                & \multicolumn{1}{l|}{}                                                                         & \multicolumn{1}{l|}{\begin{tabular}[c]{@{}l@{}}IP \\ runtime \\ (second)\end{tabular}} & \multicolumn{1}{l|}{\begin{tabular}[c]{@{}l@{}}LSLPR\\ runtime \\ (8-\\ thread)\\ (second)\end{tabular}} & \multicolumn{1}{l|}{\begin{tabular}[c]{@{}l@{}}LSLPR\\ runtime \\ (max-\\ thread)\\ (second)\end{tabular}} & \multicolumn{1}{l|}{\begin{tabular}[c]{@{}l@{}}Opt.\\ Gap\end{tabular}} & \multicolumn{1}{l|}{\begin{tabular}[c]{@{}l@{}}MMO \\ (8-\\ thread)\\ (second)\end{tabular}} & \multicolumn{1}{l|}{\begin{tabular}[c]{@{}l@{}}MMO \\ (max-\\ thread)\\ (second)\end{tabular}} & \multicolumn{1}{l}{\begin{tabular}[c]{@{}l@{}}Opt.\\  Gap\end{tabular}} \\ \hline
1.78                                                                                                     & 10                    & 115                   & 366550                                                                                               & 50                                                                                            & 1177                                                                                   & \multicolumn{1}{c|}{384}                                                                                 & \multicolumn{1}{c|}{20}                                                                                    & 0.2\%                                                                  & \multicolumn{1}{c|}{38}                                                                      & \multicolumn{1}{c|}{16}                                                                        & 0.4\%                                                                   \\ \hline
1.78                                                                                                     & 15                    & 115                   & 372050                                                                                               & 50                                                                                            & 1332                                                                                   & \multicolumn{1}{c|}{383}                                                                                 & \multicolumn{1}{c|}{19}                                                                                    & 0.9\%                                                                  & \multicolumn{1}{c|}{51}                                                                      & \multicolumn{1}{c|}{19}                                                                        & 0.8\%                                                                   \\ \hline
1.78                                                                                                     & 20                    & 115                   & 384650                                                                                               & 50                                                                                            & 1470                                                                                   & \multicolumn{1}{c|}{337}                                                                                 & \multicolumn{1}{c|}{23}                                                                                    & 0.7\%                                                                  & \multicolumn{1}{c|}{58}                                                                      & \multicolumn{1}{c|}{15}                                                                        & 0.8\%                                                                   \\ \hline
1.78                                                                                                     & 25                    & 115                   & 392500                                                                                               & 50                                                                                            & 1354                                                                                   & \multicolumn{1}{c|}{357}                                                                                 & \multicolumn{1}{c|}{24}                                                                                    & 0.9\%                                                                  & \multicolumn{1}{c|}{72}                                                                      & \multicolumn{1}{c|}{26}                                                                        & 0.8\%                                                                   \\ \hline
1.78                                                                                                     & 30                    & 115                   & 403000                                                                                               & 50                                                                                            & 1464                                                                                   & \multicolumn{1}{c|}{548}                                                                                 & \multicolumn{1}{c|}{31}                                                                                    & 0.5\%                                                                  & \multicolumn{1}{c|}{82}                                                                      & \multicolumn{1}{c|}{15}                                                                        & 0.9\%                                                                   \\ \hline
1.78                                                                                                     & 35                    & 115                   & 403250                                                                                               & 50                                                                                            & 1448                                                                                   & \multicolumn{1}{c|}{599}                                                                                 & \multicolumn{1}{c|}{27}                                                                                    & 0.4\%                                                                  & \multicolumn{1}{c|}{97}                                                                      & \multicolumn{1}{c|}{17}                                                                        & 1.0\%                                                                   \\ \hline
1.83                                                                                                    & 10                    & 115                   & 385450                                                                                               & 50                                                                                            & 1425                                                                                   & \multicolumn{1}{c|}{460}                                                                                 & \multicolumn{1}{c|}{29}                                                                                    & 0.3\%                                                                  & \multicolumn{1}{c|}{48}                                                                      & \multicolumn{1}{c|}{15}                                                                        & 0.1\%                                                                   \\ \hline
1.83                                                                                                    & 15                    & 115                   & 407100                                                                                               & 50                                                                                            & 1516                                                                                   & \multicolumn{1}{c|}{509}                                                                                 & \multicolumn{1}{c|}{63}                                                                                    & 0.1\%                                                                  & \multicolumn{1}{c|}{60}                                                                      & \multicolumn{1}{c|}{17}                                                                        & 0.0\%                                                                   \\ \hline
1.83                                                                                                    & 20                    & 115                   & 427400                                                                                               & 50                                                                                            & 1730                                                                                   & \multicolumn{1}{c|}{449}                                                                                 & \multicolumn{1}{c|}{33}                                                                                    & 0.5\%                                                                  & \multicolumn{1}{c|}{77}                                                                      & \multicolumn{1}{c|}{19}                                                                        & 0.3\%                                                                   \\ \hline
1.83                                                                                                  & 25                    & 115                   & 433950                                                                                               & 50                                                                                            & 1817                                                                                   & \multicolumn{1}{c|}{483}                                                                                 & \multicolumn{1}{c|}{35}                                                                                    & 0.3\%                                                                  & \multicolumn{1}{c|}{91}                                                                      & \multicolumn{1}{c|}{22}                                                                        & 0.7\%                                                                   \\ \hline
1.83                                                                                                   & 30                    & 115                   & 458550                                                                                               & 50                                                                                            & 2027                                                                                   & \multicolumn{1}{c|}{566}                                                                                 & \multicolumn{1}{c|}{31}                                                                                    & 0.7\%                                                                  & \multicolumn{1}{c|}{104}                                                                     & \multicolumn{1}{c|}{21}                                                                        & 0.8\%                                                                   \\ \hline
1.83                                                                                                    & 35                    & 115                   & 478000                                                                                               & 50                                                                                            & 2373                                                                                   & \multicolumn{1}{c|}{565}                                                                                 & \multicolumn{1}{c|}{44}                                                                                    & 0.5\%                                                                  & \multicolumn{1}{c|}{137}                                                                     & \multicolumn{1}{c|}{31}                                                                        & 1.2\%                                                                   \\ \hline
1.87                                                                                                   & 10                    & 115                   & 356250                                                                                               & 50                                                                                            & 1180                                                                                   & \multicolumn{1}{c|}{353}                                                                                 & \multicolumn{1}{c|}{26}                                                                                    & 0.7\%                                                                  & \multicolumn{1}{c|}{87}                                                                      & \multicolumn{1}{c|}{65}                                                                        & 0.3\%                                                                   \\ \hline
1.87                                                                                                    & 15                    & 115                   & 375050                                                                                               & 50                                                                                            & 1350                                                                                   & \multicolumn{1}{c|}{398}                                                                                 & \multicolumn{1}{c|}{12}                                                                                    & 0.5\%                                                                  & \multicolumn{1}{c|}{66}                                                                      & \multicolumn{1}{c|}{34}                                                                        & 0.1\%                                                                   \\ \hline
1.87                                                                                                    & 20                    & 115                   & 383100                                                                                               & 50                                                                                            & 1509                                                                                   & \multicolumn{1}{c|}{455}                                                                                 & \multicolumn{1}{c|}{55}                                                                                    & 1.0\%                                                                  & \multicolumn{1}{c|}{78}                                                                      & \multicolumn{1}{c|}{38}                                                                        & 0.1\%                                                                   \\ \hline
1.87                                                                                                    & 25                    & 115                   & 408650                                                                                               & 50                                                                                            & 1650                                                                                   & \multicolumn{1}{c|}{488}                                                                                 & \multicolumn{1}{c|}{18}                                                                                    & 1.0\%                                                                  & \multicolumn{1}{c|}{128}                                                                     & \multicolumn{1}{c|}{70}                                                                        & 0.3\%                                                                   \\ \hline
1.87                                                                                                    & 30                    & 115                   & 427350                                                                                               & 50                                                                                            & 1690                                                                                   & \multicolumn{1}{c|}{536}                                                                                 & \multicolumn{1}{c|}{35}                                                                                    & 1.0\%                                                                  & \multicolumn{1}{c|}{105}                                                                     & \multicolumn{1}{c|}{32}                                                                        & 0.3\%                                                                   \\ \hline
1.87                                                                                                   & 35                    & 115                   & 435950                                                                                               & 50                                                                                            & 1795                                                                                   & \multicolumn{1}{c|}{535}                                                                                 & \multicolumn{1}{c|}{77}                                                                                    & 1.3\%                                                                  & \multicolumn{1}{c|}{141}                                                                     & \multicolumn{1}{c|}{42}                                                                        & 0.3\%                                                                   \\ \hline
1.93                                                                                                    & 10                    & 115                   & 579800                                                                                               & 50                                                                                            & 4468                                                                                   & \multicolumn{1}{c|}{1011}                                                                                & \multicolumn{1}{c|}{196}                                                                                   & 0.4\%                                                                  & \multicolumn{1}{c|}{121}                                                                     & \multicolumn{1}{c|}{63}                                                                        & 0.2\%                                                                   \\ \hline
1.93                                                                                                    & 15                    & 115                   & 619000                                                                                               & 50                                                                                            & 6411                                                                                   & \multicolumn{1}{c|}{1036}                                                                                & \multicolumn{1}{c|}{58}                                                                                    & 0.5\%                                                                  & \multicolumn{1}{c|}{250}                                                                     & \multicolumn{1}{c|}{140}                                                                       & 0.4\%                                                                   \\ \hline
1.93                                                                                                    & 20                    & 115                   & 651850                                                                                               & 50                                                                                            & 11243                                                                                  & \multicolumn{1}{c|}{1237}                                                                                & \multicolumn{1}{c|}{204}                                                                                   & 0.2\%                                                                  & \multicolumn{1}{c|}{224}                                                                     & \multicolumn{1}{c|}{72}                                                                        & 0.3\%                                                                   \\ \hline
1.93                                                                                                    & 25                    & 115                   & 692250                                                                                               & 50                                                                                            & 11820                                                                                  & \multicolumn{1}{c|}{1318}                                                                                & \multicolumn{1}{c|}{28}                                                                                    & 0.4\%                                                                  & \multicolumn{1}{c|}{390}                                                                     & \multicolumn{1}{c|}{167}                                                                       & 0.2\%                                                                   \\ \hline
1.93                                                                                                   & 30                    & 115                   & 734800                                                                                               & 50                                                                                            & 5432                                                                                   & \multicolumn{1}{c|}{1453}                                                                                & \multicolumn{1}{c|}{99}                                                                                    & 0.8\%                                                                  & \multicolumn{1}{c|}{397}                                                                     & \multicolumn{1}{c|}{86}                                                                        & 0.6\%                                                                   \\ \hline
1.93                                                                                                    & 35                    & 115                   & 845950                                                                                               & 50                                                                                            & 7038                                                                                   & \multicolumn{1}{c|}{1830}                                                                                & \multicolumn{1}{c|}{64}                                                                                    & 0.6\%                                                                  & \multicolumn{1}{c|}{516}                                                                     & \multicolumn{1}{c|}{133}                                                                       & 0.4\%                                                                   \\ \hline
2.02                                                                                                    & 10                    & 115                   & 593200                                                                                               & 50                                                                                            & 4348                                                                                   & \multicolumn{1}{c|}{998}                                                                                 & \multicolumn{1}{c|}{50}                                                                                    & 0.6\%                                                                  & \multicolumn{1}{c|}{226}                                                                     & \multicolumn{1}{c|}{182}                                                                       & 1.0\%                                                                   \\ \hline
2.02                                                                                                    & 15                    & 115                   & 662350                                                                                               & 50                                                                                            & 5843                                                                                   & \multicolumn{1}{c|}{1792}                                                                                & \multicolumn{1}{c|}{157}                                                                                   & 0.2\%                                                                  & \multicolumn{1}{c|}{175}                                                                     & \multicolumn{1}{c|}{93}                                                                        & 1.0\%                                                                   \\ \hline
2.02                                                                                                    & 20                    & 115                   & 860950                                                                                               & 50                                                                                            & 9802                                                                                   & \multicolumn{1}{c|}{2177}                                                                                & \multicolumn{1}{c|}{210}                                                                                   & 0.2\%                                                                  & \multicolumn{1}{c|}{442}                                                                     & \multicolumn{1}{c|}{214}                                                                       & 1.0\%                                                                   \\ \hline
2.02                                                                                                    & 25                    & 115                   & 1123050                                                                                              & 50                                                                                            & 16787                                                                                  & \multicolumn{1}{c|}{3508}                                                                                & \multicolumn{1}{c|}{300}                                                                                   & 0.5\%                                                                  & \multicolumn{1}{c|}{609}                                                                     & \multicolumn{1}{c|}{201}                                                                       & 0.9\%                                                                   \\ \hline
2.02                                                                                                    & 30                    & 115                   & 1212400                                                                                              & 50                                                                                            & 22141                                                                                  & \multicolumn{1}{c|}{4143}                                                                                & \multicolumn{1}{c|}{198}                                                                                   & 0.6\%                                                                  & \multicolumn{1}{c|}{1234}                                                                    & \multicolumn{1}{c|}{477}                                                                       & 0.4\%                                                                   \\ \hline
2.02                                                                                                    & 35                    & 115                   & 1243050                                                                                              & 50                                                                                            & 20920                                                                                  & \multicolumn{1}{c|}{6327}                                                                                & \multicolumn{1}{c|}{161}                                                                                   & 0.4\%                                                                  & \multicolumn{1}{c|}{1362}                                                                    & \multicolumn{1}{c|}{487}                                                                       & 0.5\%                                                                   \\ \bottomrule 
\end{tabular} 
}}
\end{center}

\hspace{0.25in} \text{\blue{{\scriptsize 
- The demand-supply ratio is the average trip requests over the total number of vehicles ($K+|S_B|$); $K=5$. }  }}

\hspace{0.3in}\text{\blue{{\scriptsize - 8-thread and max-thread are the number of parallel programs; they are not feasible for the IP benchmark.}}}

\hspace{0.3in}\text{\blue{{\scriptsize - The max-thread runtimes assume enough threads to evaluate all potential swaps or drivers at once: } }}

\hspace{0.35in} \text{\blue{{\scriptsize In LSLRP, max threads = 150; in MMO, max threads = 30.}  }}

\end{table}

The overall performance of the proposed approximation algorithms is surprisingly satisfactory. 
The optimality gap is below 2\% throughout all experiments. 
These results confirm that those approximation ratios proven for the worst-case, $\frac{1}{p^2}$ or $\frac{e-1}{(2e + o(1)) p \ln p}$, are loose with the real-world trip data.  
In other words, the performance degradation of these approximation algorithms is negligible when implementing them in shared mobility systems.

\subsubsection{Sensitivity analysis. }
Three sensitivity analyses conducted in the sensitivity analysis are a) distribution of hyperedge values, b) vehicle number and capacity, and c) sample size.
They test how the performance of these approximation algorithms is affected based on changes in input data and model assumptions.
We discuss them individually in this section.

\begin{figure}[!htb]
    \centering
    \subfloat[Computation time comparison]{ \includegraphics[width = 0.7\textwidth]{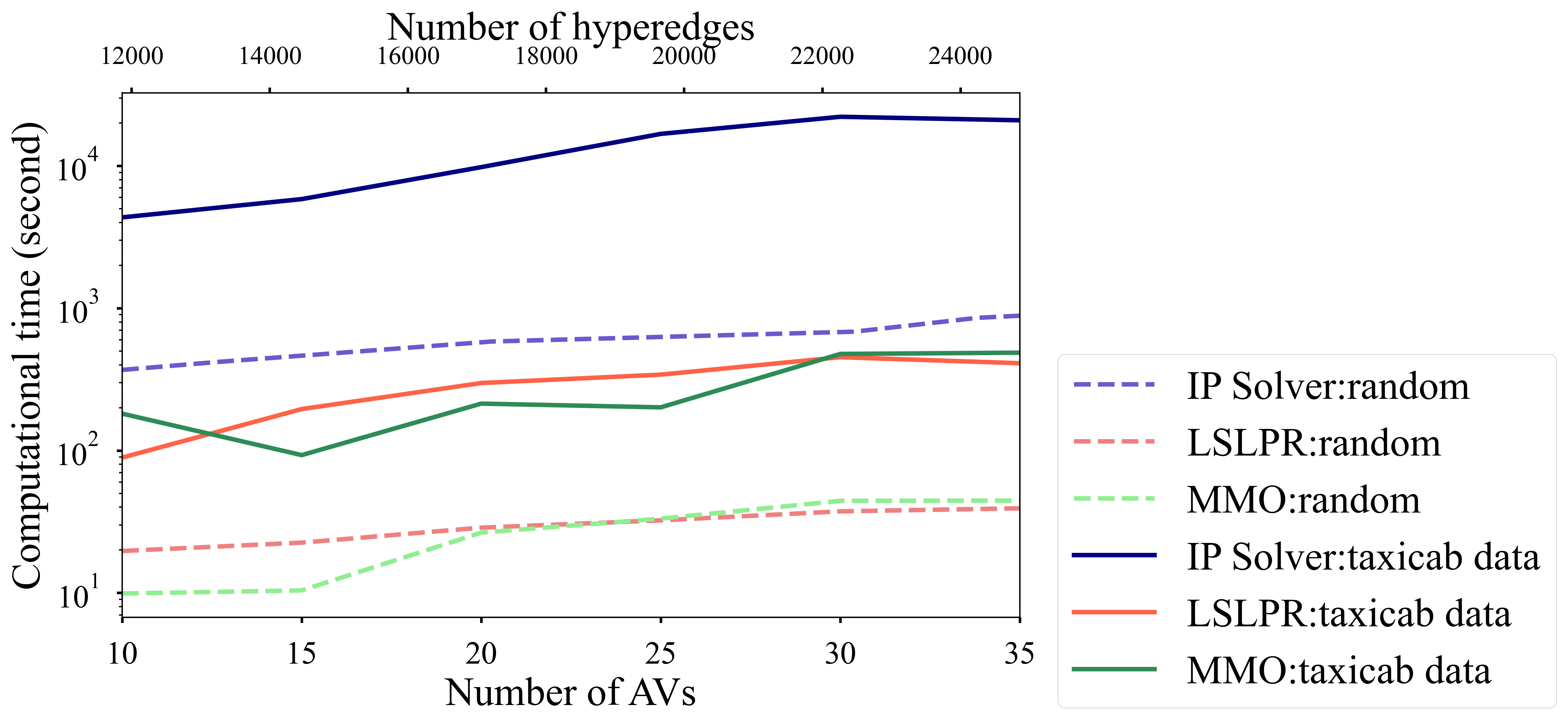} }  \\
     \subfloat[Optimality gap comparison]{ \includegraphics[width = 0.7\textwidth]{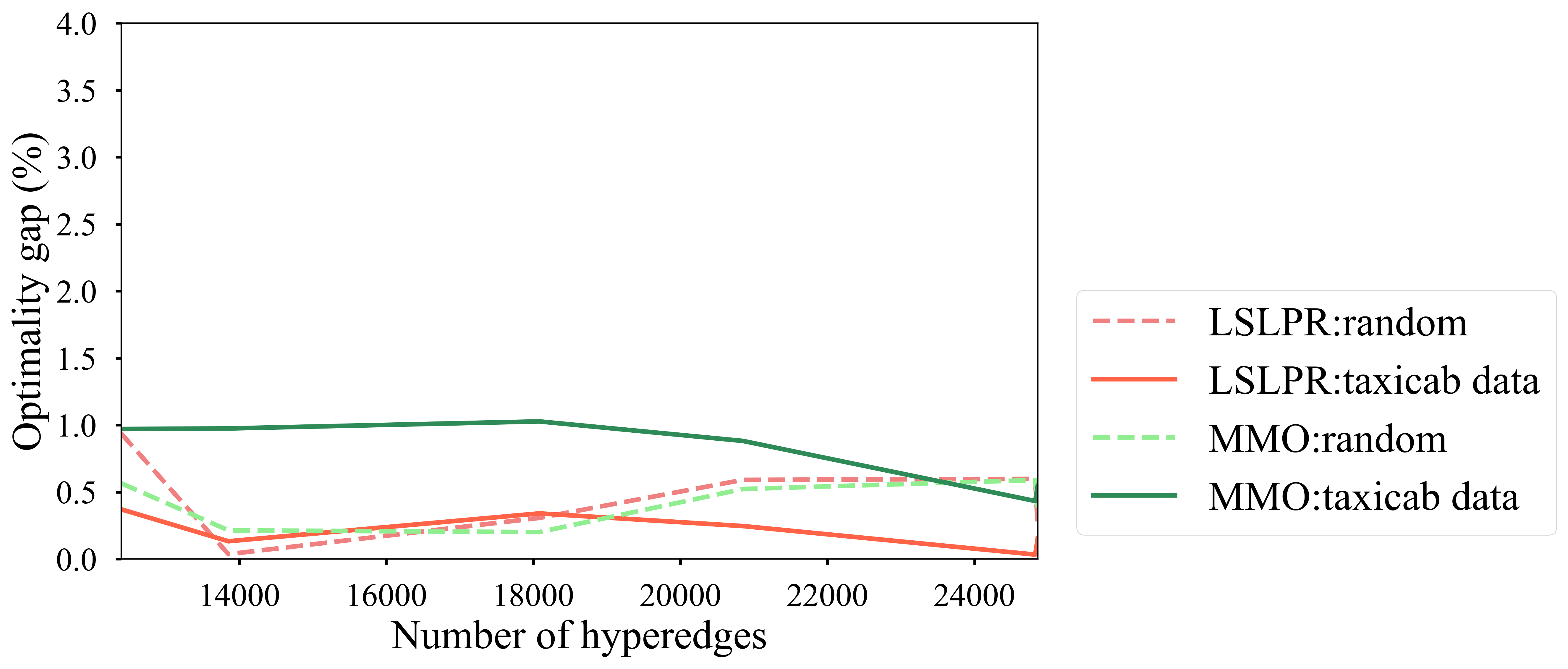} }
     
    \caption{Impact of input distribution on computation time and optimality gap.}
    \label{fig:9}
\end{figure}

\noindent \textbf{Hyperedge value distribution.} \quad
The first set of sensitivity analyses aims to check algorithms' performance degrades with different supply and demand distributions. 
By replacing the empirical hyperedge values with randomly generated hyperedge values, this analysis examines the robustness of these algorithms. 
Figure \ref{fig:9} shows the runtime and optimality gaps with uniformly generated hyperedge values (see two distributions in Figure \ref{fig:7-3}). 
This uniform distribution represents that the customers' trust in AV technology is the dominating factor in the hyperedge value such that 
$v_{e} = \sum_{j\in t} u_j  +  \sum_{j\in t} \tilde{u}_{ij}  - c(i, t) \approx \sum_{j\in t} \tilde{u}_{ij} $ for each $e\in E(\xi)$ and the joint utility function follows a uniform distribution. 

The runtime of random hyperedge values is smaller than those of real-world data, and the optimality gaps stay low across most instances. 
This is mainly because the empirical hyperedge values are more concentrated around specific values (i.e., the average trip length). 
Hence it is more difficult to find the best vehicle from the augmented set. 
In this case, the local-search-based algorithms are more efficient with more uniformly distributed values.

\noindent \textbf{Vehicle capacity.} \quad
AVs are the mass transport carriers in the mixed autonomy numerical experiment, whose capacity \blue{is up to ten requests and each request may contain more than one passenger.}
CVs have a fixed capacity of three throughout the experiments. 
Recall that $p$ bounds the vehicle capacity, Table \ref{tab:3} shows how the vehicle capacity affects the approximation ratios.
Observe that the vehicle capacity is not the bottleneck of the performance. 
Since the number of hyperedges increases with the increasing vehicle capacity, the IP benchmark's computational time increases.
However, we observe that the number of hyperedges plateaus above a certain capacity.
This is due to the method of construction of the shareability graph, as described in Section \ref{sec:5.1}.
In order for a high capacity trip to exist, all subset trips must also exist.
This leads to a combinatorially decreasing number of hyperedges with large trip set sizes, unless an even larger set of compatible trips exist.
At the density of requests chosen in this experiment in the AV zones, we do not see these large sets of compatible trips.
The optimality gaps of both approximation algorithms are not evidently affected by the AV capacity.

\begin{table}[!htb]
\centering
 \caption{Impact of vehicle capacity on computation time and optimality gap} \label{tab:3}
{\small
\begin{tabular}{c|c|c|c|c|c}
\toprule
\begin{tabular}[c]{@{}c@{}}Number of \\ AV locations\end{tabular} & \multicolumn{1}{l|}{\begin{tabular}[c]{@{}l@{}}AV capacity\\ $C_{AV}$\end{tabular}} & \begin{tabular}[c]{@{}c@{}}Number of \\ hyperedges\end{tabular} & \multicolumn{1}{l|}{\begin{tabular}[c]{@{}l@{}}Runtime of \\ IP (second)\end{tabular}} & \multicolumn{1}{l|}{\begin{tabular}[c]{@{}l@{}}Optimality \\ gap of LSLPR (\%)\end{tabular}} & \multicolumn{1}{l}{\begin{tabular}[c]{@{}l@{}}Optimality \\ gap of MMO (\%)\end{tabular}} \\ \midrule
35                                                                & 2                                                                              & 8121                                                            & 220                                                                                     & 0.62                                                                                         & 1.12                                                                                      \\ \hline
35                                                                & 4                                                                              & 11396                                                           & 288                                                                                     & 0.80                                                                                         & 0.88                                                                                      \\ \hline
35                                                                & 6                                                                              & 12048                                                           & 299                                                                                     & 0.99                                                                                         & 0.43                                                                                      \\ \hline
35                                                                & 8                                                                              & 12068                                                           & 298                                                                                     & 1.00                                                                                         & 0.03                                                                                      \\ \hline
35                                                                & 10                                                                             & 12070                                                           & 301                                                                                     & 1.00                                                                                         & 0.02                                                                                      \\ 
\bottomrule
\end{tabular}
 }
\end{table}

\noindent \textbf{Sample size.} \quad 
Although SAA guarantees a uniform convergence to the optimal value, it is not clear how the number samples affect the computation times and the optimality gaps. 
The results of running the same experiments with increasing sample size are summarized in Table \ref{tab:4}.
This table reports the runtime with unlimited computational resources (i.e., evaluating all scenarios in parallel). The cases with finite multithreading can refer to Figure \ref{fig:7}.
 
\begin{table}[!htb]
\centering
\caption{Impact of sample size on computation time and optimality gap}
    \label{tab:4}
    {\small
\begin{tabular}{c|c|c|c|c|c|c}
\toprule
\begin{tabular}[c]{@{}c@{}}Number of \\ samples\end{tabular} & \begin{tabular}[c]{@{}c@{}}Number of \\ AV locations\end{tabular} & AV capacity & \begin{tabular}[c]{@{}c@{}}Avg number \\ of hyperedges\end{tabular} & \begin{tabular}[c]{@{}c@{}}Runtime of\\ IP (second)\end{tabular} & \begin{tabular}[c]{@{}c@{}}Runtime of \\ LSLPR \\ (second)\end{tabular} & \begin{tabular}[c]{@{}c@{}}Optimality \\ gap (\%)\end{tabular} \\ \midrule
10                                                           & 20                                                                &5           & 3100                                                                & 14                                                               & 2                                                                       & 3.23                                                           \\ \hline
25                                                           & 20                                                                & 5           & 3100                                                                & 104                                                              & 5                                                                       & 1.63                                                           \\ \hline
50                                                           & 20                                                                & 5           & 3100                                                                & 445                                                              & 12                                                                      & 2.03                                                           \\ \hline
75                                                           & 20                                                                & 5           & 3100                                                                & 907                                                              & 15                                                                      & 2.01                                                           \\ \hline
100                                                          & 20                                                                & 5           & 3100                                                                & 2088                                                             & 25                                                                      & 0.91                                                           \\ \hline
150                                                          & 20                                                                & 5           & 3100                                                                & 4076                                                             & 38                                                                      & 3.15                                                           \\ \hline
200                                                          & 20                                                                & 5           & 3100                                                                & 7485                                                             & 109                                                                     & 1.01                                                           \\ \bottomrule
\end{tabular}
 }
\end{table}

\subsection{\blue{Numerical Results for Mid-Capacity SRAMF}}

\begin{figure}[!htb]
    \centering
    \subfloat[\blue{CV routes}]{\includegraphics[width=0.3\textwidth ]{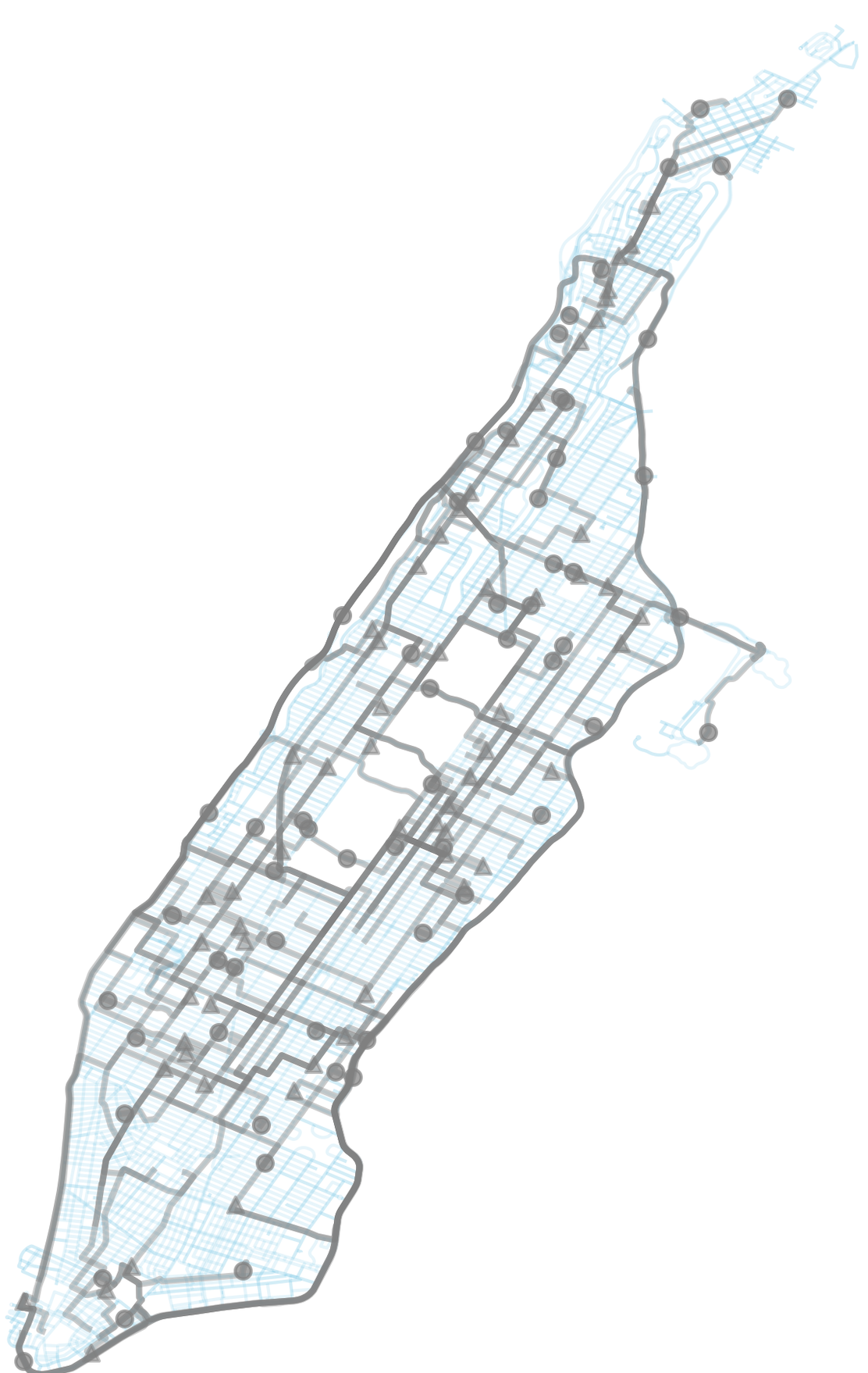} }
    \quad 
    \subfloat[\blue{AV routes}]{\includegraphics[width=0.45\textwidth]{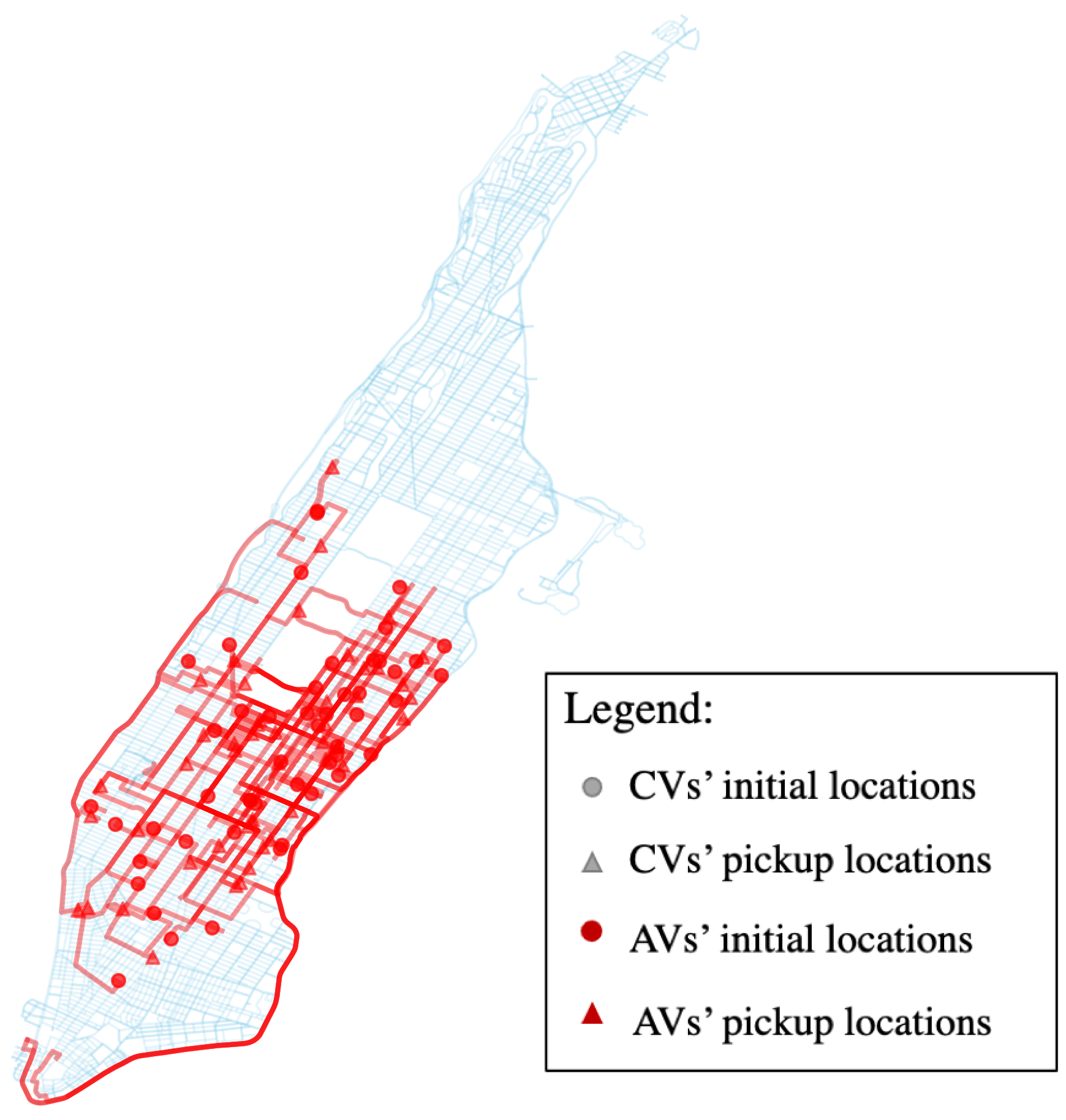} }

    \caption{\blue{Optimal trip assignment and routes in mixed autonomy mid-capacity scenario.} }
    \label{fig:map_ev}
\end{figure}

\blue{In Setting 2, a relatively large fleet of electric AVs ($|S_A| = 114$) and CVs ($|S_B| = 60$ serve the travel demand in tandem over the entire city area.
The results are summarized in Table \ref{tabler1-2}. We also show the optimal assignment and routes in Figure \ref{fig:map_ev}}.

\begin{table}[!htb]
\begin{center}
\caption{\blue{Summary of Numerical Results for Mid-Capacity SRAMF}} \label{tabler1-2}
\blue{{\small 
\resizebox{\textwidth}{!}{
\begin{tabular}{c|c|c|c|c|c|c|ccc|ccc}
\toprule
\multicolumn{1}{l|}{\multirow{2}{*}{\begin{tabular}[c]{@{}l@{}}Demand-\\ supply \\ ratio\end{tabular}}} & \multirow{2}{*}{S\_A} & \multirow{2}{*}{S\_B} & \multirow{2}{*}{K} & \multicolumn{1}{l|}{\multirow{2}{*}{\begin{tabular}[c]{@{}l@{}}Total \# IP\\ variables \\ (hyperedges \\ $\times$ samples)\end{tabular}}} & \multicolumn{1}{l|}{\multirow{2}{*}{\begin{tabular}[c]{@{}l@{}}\# of\\ samples\end{tabular}}} & Benchmark                                                                              & \multicolumn{3}{c|}{LSLPR}                                                                                                                                                                                                                                                                      & \multicolumn{3}{c}{MMO}                                                                                                                                                                                                                                                 \\ \cline{7-13} 
\multicolumn{1}{l|}{}                                                                                   &                       &     &                  & \multicolumn{1}{l|}{}                                                                                & \multicolumn{1}{l|}{}                                                                         & \multicolumn{1}{l|}{\begin{tabular}[c]{@{}l@{}}IP \\ runtime \\ (second)\end{tabular}} & \multicolumn{1}{l|}{\begin{tabular}[c]{@{}l@{}}LSLPR\\ runtime \\ (36-\\ thread)\\ (second)\end{tabular}} & \multicolumn{1}{l|}{\begin{tabular}[c]{@{}l@{}}LSLPR\\ runtime \\ (max-\\ thread)\\ (second)\end{tabular}} & \multicolumn{1}{l|}{\begin{tabular}[c]{@{}l@{}}Opt.\\ Gap\end{tabular}} & \multicolumn{1}{l|}{\begin{tabular}[c]{@{}l@{}}MMO \\ (36-\\ thread)\\ (second)\end{tabular}} & \multicolumn{1}{l|}{\begin{tabular}[c]{@{}l@{}}MMO \\ (max-\\ thread)\\ (second)\end{tabular}} & \multicolumn{1}{l}{\begin{tabular}[c]{@{}l@{}}Opt.\\  Gap\end{tabular}} \\ \hline
22.7 & 114 & 60 & 30& 8468768 & 20 & 2656 & \multicolumn{1}{c|}{718} & \multicolumn{1}{c|}{213} & 0.011\% & \multicolumn{1}{c|}{263} & \multicolumn{1}{c|}{87.6} & 0.11\% \\ \hline
34.1 & 114 & 60 & 30& 18124448 & 20 & 5622 & \multicolumn{1}{c|}{1238} & \multicolumn{1}{c|}{336} & 0.015\% & \multicolumn{1}{c|}{573} & \multicolumn{1}{c|}{191} & 0.08\% \\ \hline
45.5 & 114 & 60 & 30& 31516528 & 20 & 8934 & \multicolumn{1}{c|}{3814} & \multicolumn{1}{c|}{2151} & 0.001\% & \multicolumn{1}{c|}{1050} & \multicolumn{1}{c|}{350} & 0.05\% \\ \hline
22.7 & 114 & 60 & 30& 12572544 & 30 & 5322 & \multicolumn{1}{c|}{672} & \multicolumn{1}{c|}{267} & 0.037\% & \multicolumn{1}{c|}{381} & \multicolumn{1}{c|}{127} & 0.48\% \\ \hline
34.1 & 114 & 60 & 30& 27495912 & 30 & 11559 & \multicolumn{1}{c|}{983} & \multicolumn{1}{c|}{446} & 0.036\% & \multicolumn{1}{c|}{838} & \multicolumn{1}{c|}{279} & 0.26\% \\ \hline
45.5 & 114 & 60 & 30& 47274792 & 30 & 19763 & \multicolumn{1}{c|}{1913} & \multicolumn{1}{c|}{1504} & 0.013\% & \multicolumn{1}{c|}{1457} & \multicolumn{1}{c|}{485} & 0.14\% \\ \hline
22.7 & 114 & 60 & 60& 8195056 & 20 & 2608 & \multicolumn{1}{c|}{2665} & \multicolumn{1}{c|}{9.26} & 0.007\% & \multicolumn{1}{c|}{469} & \multicolumn{1}{c|}{227} & 0.20\% \\ \hline
34.1 & 114 & 60 & 60& 17905712 & 20 & 6277 & \multicolumn{1}{c|}{4898} & \multicolumn{1}{c|}{55.9} & 0.004\% & \multicolumn{1}{c|}{1022} & \multicolumn{1}{c|}{520} & 0.11\% \\ \hline
45.5 & 114 & 60 & 60& 31516528 & 20 & 11080 & \multicolumn{1}{c|}{28051} & \multicolumn{1}{c|}{72.0} & 0.004\% & \multicolumn{1}{c|}{1998} & \multicolumn{1}{c|}{1620} & 0.076\% \\ \hline
22.7 & 114 & 60 & 60& 12867024 & 30 & 5599 & \multicolumn{1}{c|}{3811} & \multicolumn{1}{c|}{56.7} & 0.005\% & \multicolumn{1}{c|}{649} & \multicolumn{1}{c|}{544} & 0.39\% \\ \hline
34.1 & 114 & 60 & 60& 27147312 & 30 & 11220 & \multicolumn{1}{c|}{34360} & \multicolumn{1}{c|}{79.0} & 0.001\% & \multicolumn{1}{c|}{1516} & \multicolumn{1}{c|}{1210} & 0.14\% \\ \hline
45.5 & 114 & 60 & 60& 47274792 & 30 & DNF & \multicolumn{1}{c|}{DNF} & \multicolumn{1}{c|}{DNF} & - & \multicolumn{1}{c|}{2962} & \multicolumn{1}{c|}{1407} & - \\ \bottomrule
\end{tabular}
}
}}
\end{center}

\hspace{0.1in}\text{\blue{{\scriptsize - 36-thread and max-thread are the number of parallel programs; they are not feasible for the IP benchmark. }  }}

\hspace{0.1in}\text{\blue{{\scriptsize - The max-thread runtimes assume enough threads to evaluate all potential swaps or drivers at once: }  }}

\hspace{0.2in}\text{\blue{{\scriptsize In LSLRP, max threads = 3420; in MMO, max threads = 114).}  }}
\end{table}

\subsubsection{Computation times. }
\blue{At low sample sizes, the runtime of solving the IP is fairly competitive with the approximation algorithms. 
Gurobi is a powerful and heavily optimized solver, so this is not very surprising. 
However, the advantage of MMO and LSLRP at higher sample sizes is clear because they can leverage parallel computation resources. }

\blue{The runtime of LSLRP varies widely as the swap order greatly affects the runtime. 
In the $K=60$ case, we occasionally observe very long runtimes for the LSLRP algorithm. 
This is primarily due to the number of swaps that the algorithm has to evaluate in order to guarantee that it has found a solution.
If the algorithm randomly finds a good swap, finding another one that further improves the objective value becomes increasingly difficult. 
This can lead to the algorithm evaluating a combinatorial number of swaps but the quality of approximation almost attains the optimality.}


\subsection{Limitations}
\blue{
Since this work focuses on offline algorithms, the limitations of the current experiments include: 
\begin{enumerate}
	\item 	The system does not allow alternative pickups and dropoffs in trip planning, i.e., the total number of requests is no greater than vehicle capacity in each trip clique. 
	\item The penalties of balking trips or carryover supply or demand are not directly considered, but can be incorporated into the hyperedge value as discussed in Section \ref{sec:3-1}.
\end{enumerate}
}


The computational time increases linearly with the sample size.
The optimality gaps are small across different sample sizes. 
The algorithm can use a large sample size when the stochasticity in the platform is of major interest and has access to abundant computational resources.

\section{Conclusion} \label{sec:6}

\blue{SRAMF is a generic formulation for operating shared mobility platforms with blended workforces or mixed autonomy traffic. 
This paper proposes mid-capacity and high-capacity approximation algorithms for joint fleet sizing and trip assignment in ride-pooling platforms.
Widely used real-time ride-pooling frameworks \citep{alonso2017demand,simonetto2019real} can integrate these approximation algorithms to accelerate the vehicle dispatching process and improve the quality of service with mixed fleets.}
We give provable guarantees for their worst-case performance and test their average performance in \blue{numerical experiments}.

To close this paper, we point out several promising future research avenues to address the following limitations. 
First, we focus on the \blue{SRAMF} problem in the face of a simple uncertainty structure (i.e., two-stage decision). 
Since the demand forecast can be time-varying and the vehicle repositioning decisions should adapt to revealing scenarios, extending our framework to a multi-stage setting would be interesting.
\blue{Second, the construction of hypergraph remains a computational challenge in the worst-case, hence developing efficient reformulation for updated batches of potential matchings may provide computational advantages.}
Finally, the current trip assignment does not consider cancellation and re-assignment after dispatching vehicles to passengers. 
Considering these factors in practice may improve the stability of ride-pooling algorithms.

\ACKNOWLEDGMENT{%
The work described in this paper was partly supported by research grants from the National Science Foundation (CMMI-1854684; CMMI-1904575; CMMI-1940766;  CCF-2006778) and Ford-UM Faculty Research Alliance. 
}

%
%
%


\begin{APPENDICES}
 \section{Summary of \blue{Notation}} \label{sec:app1}
%
\begin{longtable}{l|l}
        \caption{\blue{Summary of notation and acronyms}} \label{table1} \tabularnewline
		\toprule
		\multicolumn{1}{c|}{{\bf Notation} } & \multicolumn{1}{c}{ {\bf Description} } \\ \hline
		$S_A$, $S_B$     &    Augmented set and basis set of vehicles    \\ \hline
		$\xi$      &    Randomly generated scenario  \\ \hline
		$\ell \in [N]$        & \blue{Index for sampled scenarios and the total number of samples} \\ \hline 
		$D(\xi)$ &   \blue{Set of demand in scenario $\xi$} \\ \hline
		$E(\xi)$ & Set of hyperedges in scenario $\xi$  \\ \hline 
		\blue{$G$} & \blue{Shareability graph, a hypergraph consists of supply and demand vertices and hyperedges } \\ \hline 
		\blue{$e$} & 
		\blue{Each hyperedge $e = \{i, J\}_{i \in S, J\subseteq D }$ is a potential trip where vehicle $i$ serves requests $J$} \\ \hline
		\blue{$u_j$} & \blue{The expected profit of request $j$ } \\ \hline
		\blue{$\tilde{u}_{ij}$} & \blue{Utility gained from matching request $j$ with preferred vehicle type $i$ } \\ \hline
		\blue{$c(i, t)$} & \blue{Travel cost for vehicle $i$ to serve trip $t$ } \\ \hline
		\blue{$\alpha$} & \blue{Approximation ratio}    \\ \hline
		$n$        &    Total number of vehicles such that  $|S_A| = n_A$ and $|S_B| = n_B$    \\ \hline 
		$K$        &    Maximum number of vehicles allowed from the augmented set   \\ \hline
		$C_i$      &    Capacity of vehicle $i \in S_A\cup S_B$      \\ \hline
		\blue{$w_i$} &  \blue{Number of passengers in request $j$} \\ \hline   
		$p$        &    Maximum capacity of hyperedge, \blue{$p = \max_{i\in S_A\cup S_B} \{1 + C_i\}$} \\ \hline 
		\blue{$j$}        &    \blue{Index for travel demand $j\in D(\xi)$ } \\ \hline
		$w_j$      &    Size of travel demand  $j\in D(\xi)$   \\ \hline
		$E_{i, \xi}$ & Set of hyperedges contains vertex $i \in S_B \cup S_R$ in  scenario $\xi$ \\ \hline 
		$t$        &    Trip is a set of demand following \blue{the shortest pickup-and-then-dropoff order} \\ \hline
		$v_e$    &      Value of hyperedge  $e \in E(\xi)$ \\ \hline 
		$nb(e)$   &     Neighboring hyperedges $e'\in E(\xi)$ intersecting with $e$ \\ \hline
		$x_{e}$    &    Decision variable for hyperedge $e$, $x_{e}\in \{0,1\}$ \\ \hline 
		\blue{$\bar{x}_{e}$ }   &  \blue{Decision variable for fractional assignment, $x_{e}\in [0,1]$} \\ \hline 
		$y_{i}$   &     Decision variable for vehicle $i\in [S_A]$, $y_{i} \in \{0, 1\} $ \\ \hline 
		$v^*(\cdot)$     &     Optimal value of the exact GAP   \\
		\hline 
		$Q(y, \xi)$ &   Optimal value of the assignment in scenario $\xi$ \\ \hline
		\blue{$v_{\max}$} & \blue{Maximal hyperedge value for all $e \in E$ }  \\ \hline
		\blue{$v_{\min}$} & \blue{Minimal hyperedge value for all $e \in E$ such that $v_e >0$ }  \\ \hline
		$\mathcal{I}$ & Independent set as a union of hyperedges satisfying the set-packing constraint  \\ \hline 
		$S_O$    &      Optimal choice of vehicles for SRAMF  $S_O \subset S_A$ \\ \hline 
		$S_R$    &      Choice of vehicles from the algorithm \blue{ $S_R \subset S_A$}  \\ \hline
		\blue{$\mathcal{L}$} & \blue{The bijection between $S_R$ and $S_O$}  \\ \hline
		$\hat{v}(\cdot)$ & The objective value of fractional assignment  \\ \hline 
		\blue{$\pmb{z}$} & \blue{Optimal LP solutions to $\hat{v}(S_O)$} \\ \hline
		$F_i$ & Hyperedges intersect with vehicle $i$ \\ \hline 
		$H_d$ &  Hyperedges intersect with demand $d$ \\ \hline
		\blue{$\bot$} & \blue{Dummy hyperedge in LSLPR} \\ \hline
		$\Delta_d(e, f)$ & Decomposition mapping between hyperedge $e$ and $f$ \\  \hline 
		$U_{i_1 i_2}$ & Marginal value function with $i_1\in S_R$ and $i_2 \in S_O$ \\ \hline
		\blue{$\hat{v}_{ON}$} & \blue{Objective value of the online algorithm}   \\ \hline 
		$u_{g, \xi}$ & Dual variable in the MMO algorithm for $g\in e$ and scenario $\xi$ \\ \hline 
		\blue{$\Gamma_e$} & \blue{$\Gamma_e = \sum_g u_{g, \xi}$ as the left side of dual constraints} \\ \hline
		$\pmb{c}_e$ & Row of cost coefficient in the dual covering problem with entry $c_e(g, \xi)$ \\ \hline 
		\blue{$M$} & \blue{The augmented set is partitioned into $M$ subsets}  \\ \hline
		$\epsilon$  &   Error tolerance (for stopping criteria)  \\ \hline
		$\delta$  &   Error tolerance (for sample average approximation) or constant in MMO update subroutine  \\ \hline
		\blue{$OPT$} & \blue{Optimal value of the SRAMF problem} \\ \hline 
		\blue{$ALG$} & \blue{Objective value of solving SRAMF by approximation algorithms } \\
		\midrule
		\multicolumn{1}{c|}{{\bf Acronym} } & \multicolumn{1}{c}{ {\bf Description} } \\ \hline
		CV/AV & Conventional/automated vehicle
		\\	\hline
		GAP & General assignment problem \\ \hline
		LP & Linear  program  \\	\hline
		IP & Integer  program  \\	\hline
		LSLPR & Local-search linear-program-relaxation algorithm \\ \hline 
		MMO & Max-min online algorithm \\ \hline  
		SAA & Sample average approximation \\  \hline
		\blue{SRAMF}  &  \blue{Stochastic ride-pooling assignment with mixed fleets}\\ \hline
		\blue{VRP}  &  \blue{Vehicle routing problem}\\\hline
		\blue{DVRP}  &  \blue{Dynamic vehicle routing problem}\\
		\bottomrule
\end{longtable}

\section{ \blue{Performance Analysis of Construction of Shareability Graphs}} \label{sec:app2}

\blue{The main idea of recent ride-pooling assignment papers \citep{santi2014quantifying,alonso2017demand,simonetto2019real} is to separate the problem into two parts: 1) constructing the shareability graph and compatible requests and vehicles, and 2) optimally assigning those trips to vehicles by solving GAP. 
This paper primarily focuses on algorithms and approximation bounds for the stochastic extension to the second part, but we acknowledge the importance and difficulty of the first task and describe them in detail below for completeness. } 

\subsection{\blue{Procedure for Constructing Shareability Graphs}} \label{app:b1}
$D(\xi)$ is a set of all trip requests revealed in scenario $\xi$, \blue{and this section omits $\xi$ when there is no confusion because the hyperedges for scenarios are generated separately}.
\blue{We consider a number of parameters to be given by the customer or externally dictated to the platform (based on desired service parameters). These include, for each customer $j$, the maximum waiting time, $\omega_j$, and allowable delay, $r_j$.
\begin{itemize}
    \item \textit{(Constraint I) } Travel time from vehicle location to pick-up of customer $j$ in order must be less than $\omega_j$.
    \item \textit{(Constraint II)}  Travel time from origin to destination of customer $j$ in order must be less than $r_j$.
\end{itemize}}
\blue{Additionally, as defined in the setting, the hyperedge weight consists of three parts: value of trip requests $\sum_{j} u_j$, preference of the vehicle type $\tilde{u}_{sj}$, and travel cost of delivering all trip requests in a single trip. 
We take a three-request clique $(j_1, j_2, j_3)$ as an example. 
Let $t_k = \{O_{j_1},O_{j_2},\hdots,D_{j_2},D_{j_3}\}$ be a specific ordered sequence of origins and destinations and $SP(t_k)$ be the shortest path route connecting them. 
Let $T_e = \cup_k t_k $. 
Note that this is slightly less demanding than finding all feasible Hamiltonian paths if we enforce that in all trips, the origins must be picked up before any destination is visited. 
Let $c(s, e) = \min_{t_k\in T_e} v(t_k)$ where $c(t_k)$ is the cost of serving all requests following the shortest-path $SP(t_k)$.
We define a set function $f(e)$ that takes a hyperedge consisting of a vehicle $s$ and a potential combination of trips, $T_e$, as below:
\begin{align*}
    f(e) = 
    \begin{cases}
        0  & \text{ if } \forall t_k \in T_e, \, SP(t_k) \text{ violates constraints I and II}\\ 
        c(s, e)  &  \text{ otherwise}
    \end{cases}.
\end{align*}
}
\blue{The bottleneck of computation time is still finding vehicle routes that satisfy the given constraints by solving a constrained VRP problem, which is NP-hard. 
Therefore, all heuristic methods can only minimize this bottleneck as much as possible by reducing the number of combinations to check at each step. 
For example, \cite{ke2021data} suggested a reformulation for finding $c(s, e)$ to avoid enumerating all possible paths.}

\blue{We combine multiple heuristic methods in literature to construct the shareability graph.  
First, we need to identify the valid single customer trips for a given vehicle $s$.
Let $D_s$ be the demand that can be served by vehicle $s$ in a single trip within the allowable pickup time.
We may further reduce the number of trips by planning on a spatiotemporal graph and examining compatible trips' cliques. 
By testing trips in order of increasing size and only considering a trip if all subsets of trips (where one request is removed from the trip) are feasible, we reduce the number of candidate trips by orders of magnitude. 
This heuristic generates the shareability graph in Figure \ref{fig:1} in which a set of requests is tested for trip compatibility only if every subset of that set of requests is also compatible.}
\blue{
\begin{lemma}{(\cite{alonso2017demand})}
 A trip associated with the hyperedge $e$ is feasible for vehicle $s$ only if, for all $j\sim e, j \in D_s$, hyperedges (subtrips) $e' = e\backslash \{j\}$ are feasible. 
\end{lemma}
}
\blue{The heuristic reduces the candidate hyperedge sets by leveraging the topological relationship between matchable trips of size $k$ and $k+1$ (see Figure \ref{fig:b1}), without eliminating potentially feasible trips.
The hypergraph can then be constructed in order of increasing capacity to minimize the number of request sets tested.
Additionally, we adopt the following rules to further reduce the number of candidate trips:
\begin{enumerate}
    \item Since only hyperedges with nonnegative edge weights are of interest, we remove all the trips from the candidate set subject to $f(e) <= 0$.
    \item  If a vehicle v is not feasible for trip $t_k$ at time $\tau$, it will not be feasible for $t_k$ at any time  $\tau' > \tau$ \citep{liu2020proactive}.
\end{enumerate}
}


\blue{Let $C_k(D)$ be the set of combinations of size $k$ of the elements of the set $D$. This process is summarized as follows:
\begin{center}
	\begin{minipage}{\linewidth}
		\begin{algorithm}[H]
		\caption{Construction of Shareability Graph}
		\KwData{Vehicle locations and requests (request time, pick-up, drop-off, preferred vehicle type, acceptable delay)}
		\KwResult{Set of hyperedges, $E$, each containing a vehicle, $s$, and a set of compatible requests for that vehicle to serve in one trip. Hyperedge values are $v_e$ for all $e\in E$. }
		Initialize $E = \varnothing$\\
		\For{$s \in S_A \cup S_B$}{
		    Identify candidate passengers\\
		    $D_s^{1} \leftarrow \{e \in D \mid f((s,e)) >0\}$\\
		    Add hyperedges of size one\\
		    $E_k \leftarrow \bigcup_{e \in D_s^{1}} (s,e)$ \\
		    \For{$k = 2,\hdots,c$}{
		        \For{Demand set $d \in C_k(D_s^{k-1})$}{
		            Add trips of size k if all subsets exist and value greater than 0\\
		            \If{$(s,e') \in E_{k-1} \forall i \in C_{k-1}(d)$ and $f((s,e')) >0$}{
		                $E_k \leftarrow E_k \cup (s,d)$ \\
		                $D_s^k \leftarrow D_s^k \cup d$
		            }
		        }
		    }
		}
		$E = \bigcup_{k=1}^{p-1} E_k$ \\
		Return hyperedges $E$ and their values $v_e$
		\end{algorithm}
	\end{minipage}
\end{center}
}

\begin{figure}[!htb]
    \centering
    \includegraphics[width = 0.5 \textwidth]{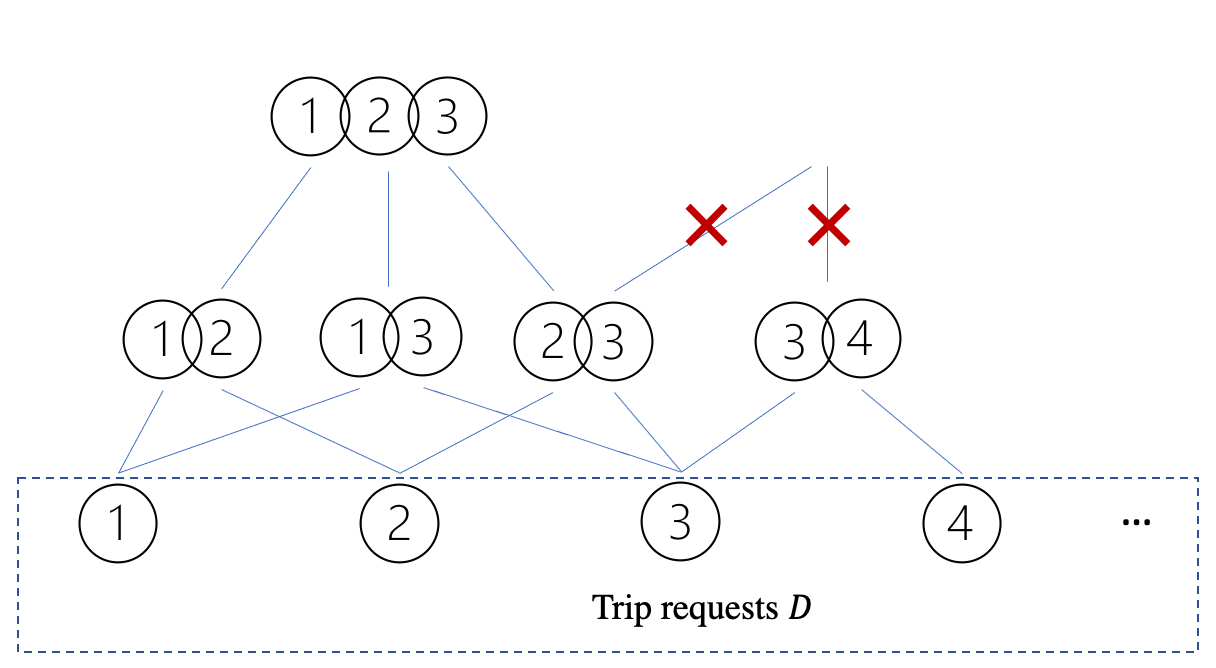}
    \caption{\blue{Topological relationship between cliques of matchable requests. In this example, (2,3,4) is not a valid combination of requests because the (2,4) combination was not valid.} }
    \label{fig:b1}
\end{figure}


\subsection{\blue{Performance and Complexity Analysis of the Hypergraph Construction Procedure}}  \label{app:b2}
\subsubsection{Optimality analysis.}

\blue{The two-step ride-pooling assignment that first constructs the hypergraph and then solves GAP obtains the exact optimal solution of solving the joint VRP and enjoys the computational advantage for large fleets. 
Since this work focuses on stochastic assignment, the optimality analysis does not consider the errors of computing hyperedge values.
The following results from \cite{alonso2017demand} provide positive guarantees for returning a feasible set of hyperedges in the shareability graph: without enumerating all trip combinations:}

\begin{enumerate}
    \item \blue{$v^*$ from solving GAP on the shareability graph obtains the optimal value for ride-pooling for an arbitrary batch of supply and demand.} 
    \item \blue{The construction of the shareability graph is anytime optimal, i.e., given additional computational resources, the set of hyperedges is only expanded to allow for improved matching.}
\end{enumerate}


\blue{The second property guarantees using a capacity bound, the threshold of which is derived below, as an early stopping criterion in generating the hypergraph will still provide satisfactory results.
Solving GAP on this reduced shareability graph can guarantee anytime optimality such that the output is near-optimal for the original problem with high probability.}

 

\subsubsection{Computational complexity analysis.}
\blue{We consider a fixed sample of demand $D$ and vehicles $S$ in this section as the hyperedges of each scenario can be generated in parallel.
The realized demand has size $|D| = d$.}

 \begin{lemma}
 \blue{In the worst case, where all demand is compatible and can be served by all vehicles, the runtime is $O(|S| d^{p-1})$}.
 \end{lemma}
 
 \blue{While in the worst-case runtime is large, this scenario only arises when all trips are compatible for ride-pooling, which is unlikely in practice.
 Therefore, we consider the Erdos-Renyi model in which an arbitrary pair of demand is matchable (i.e., satisfy the conditions above) with probability $q$.  
 Empirical studies showed that $q$ was often a small number ($<0.1$) over a large area \citep{ke2021data}.}

  \begin{lemma}[\cite{bollobas1976cliques}]
 \blue{The expected number of cliques of size $k$ is $ {d \choose k} q^{{k \choose 2}} $.}
\end{lemma}

 
 \blue{For example, with $d = 1000$ and $q=0.1$, the expected number of cliques of size $3$ (each vehicle deliver at most two requests in a single trip) is $500$.
 Often we observe the size of complete cliques of compatible trips to be less than 10, our maximum tested capacity, and the total number of hyperedges is manageable. }
 \begin{lemma}[\cite{matula1976largest}]
 \blue{As $d\rightarrow \infty$, the maximal clique size $\rho$ takes on one of at most two values around $\frac{2\log d}{log 1/q}$ with probability tending to one, i.e. with $b=1/q$, $\lfloor 2\log_{b} d \rfloor < \rho < \lceil 2\log_{b} d\rceil$.}
\end{lemma}
 
 \blue{Therefore, we only need to consider hyperedges with size less than $p^* = \min\{p, \rho + 1\}$ (i.e., the height of the cliques' graph in Figure \ref{fig:b1}).
 We have the following theorem for the runtime of constructing shareability graphs. }

\begin{theorem}
  \blue{In the average case that the demand and supply profiles satisfying the random geometric graph conditions, the runtime is  $O(|S|d^{p^*-1} )$. }
\end{theorem}

\proof{Proof:}
\blue{
The expected number of hyperedges connected to vehicle $s$,  $E_{s, \max}$, is bounded by
\begin{align*}
    E_{s, \max} &= {d \choose 1} 1^2 + {d \choose 2} 2^2 q + \dots + {d \choose (p^*-1)} (p^*-1)^2 q^{(p^*-2)} \\
    & \leq e d +  \left( \frac{e d}{2}\right)^2 2^2 q + \dots + \left(\frac{e d}{p^* - 1} \right)^{p^*-1} (p^*-1)^2 q^{(p^*-2)} \\
    &=   ed +  \frac{1}{q} \left[\left( \frac{eqd}{2}\right)^2 2^2  + \dots + \left(\frac{eqd}{p^* - 1} \right)^{p^*-1} (p^*-1)^2 \right] = O(d^{p^*-1}),
\end{align*}
where $e$ is the Euler's number.
}
\endproof

 
 \section{Supplementary Results for \blue{Approximation Algorithms} }  \label{sec:app3}
 \subsection{Proof for \blue{Sample Average Approximation in SRAMF} }
 \proof{Proof:}
 We denote the optimal value of the SRAMF problem \eqref{eq1} as $v^*$ and the optimal value of problem for the objective from Algorithm \ref{alg2} as $\hat{v}(S_O)$. 
 Let $\delta$ be the upper bound of the optimality gap $v^* - \hat{v}(S_O)$.
 \blue{We assume that $\mathbb{E}[Q(y,\xi)] = \Omega(m^{-2})$ for any $\xi$ where $m$ is a given constant.}  
 The main task is to show that:
 \begin{enumerate}
 	\item $ \mathbb{E}[\hat{v}(S_O)] = \Omega(m^2) $
 	with the sample size $N = \frac{m^4}{\delta^2}$. 
 	\item $ Pr(\hat{v}(S_O) \notin  [ (1-\delta) v^*, (1+\delta) v^* ]) \leq \exp(-\frac{\delta^2}{2} \mathbb{E}[\hat{v}(S_O)])$.
 \end{enumerate}
 
Let $D(\xi) < D$ for all $\xi$.
For any selection of vehicles in $S_A$ denoted by $y\in Y$, $\mathbb{E}[Q(y, \xi)^2] < \infty$, because we can choose $K$ vertices in $S_A$ with maximum number of  $D$ edges. 
The upper bound of of hyperedge value $v_{ij}$ is $v_{\max}$. 
Thus we have $\mathbb{E}[Q(y, \xi)^2] < K^2 |v_{\max}|^2 D^2< \infty$. 
 Without loss of generality, we draw the following observations from the standard stochastic programming literature \citep{pagnoncelli2009sample}:
 \begin{enumerate} 
 	\item $\hat{v}(S_O) \to v^*$ as $N\to \infty$;
 	\item $ \mathbb{E}[\hat{v}(S_O)] \geq v^*$.
 \end{enumerate}
 
 Since $N$ samples are i.i.d., we can use the Chernoff bound on the measure:
 \begin{align*}
 	Pr(\hat{v}(S_O) \notin  [ (1-\delta) v^*, (1+\delta) v^* ]) \leq \exp(-\frac{\delta^2}{2} \mathbb{E}[\hat{v}(S_O)]).
 \end{align*}
 
 Setting $N = m^4 /\delta^2$ and using the assumption that $\mathbb{E}[Q(y,\xi)] = \Omega(m^{-2})$, by Jensen's inequality, we have:
 $$\delta^2 \mathbb{E}[\hat{v}(S_O)]   \geq \delta^2 N \cdot E[Q(y,\xi)],$$
 i.e., $\delta^2 \cdot \mathbb{E}[\hat{v}(S_O)] = \Omega(m^{2})$. We have
 \begin{align*}
 	Pr\Big(\hat{v}(S_O) \notin  [ (1-\delta) v^*, (1+\delta) v^* ] \Big) \leq \exp(-\Omega(m^2)),
 \end{align*}
 which achieves the second task as
 \begin{align*}
 	& \blue{Pr\left( (\frac{1}{2} - \epsilon) \hat{v}(S_O)< (\frac{1}{2} - \epsilon)(1-\delta) v^* \right) + Pr\left( (\frac{1}{2} - \epsilon) 
 	\hat{v}(S_O)> (\frac{1}{2} - \epsilon)(1+\delta) v^*\right) }
 	\\
 	\leq & Pr\left(\hat{v}(S_R) < (\frac{1}{2} - \epsilon)(1-\delta) v^* \right)   +  Pr\left(\hat{v}(S_R)> (\frac{1}{2} - \epsilon)(1+\delta) v^*\right) 
 	\\
 	\leq& \exp(-\Omega(m^2)).
 \end{align*}
 The first inequality is because $\frac{1}{p^2} \hat{v}(S_O) \leq \hat{v}(S_R) \leq \hat{v}(S_O) $. This concludes the approximation ratio for LSLPR algorithm for the stochastic counterpart of the ride-pooling problem.  
 \hfill $\square$
 
 \endproof
 
 \subsection{Proof for Lemma \ref{lemma4-1}}

\blue{We use a network flow formulation to prove the existence of the mapping $\Delta_d : H'_d\times H'_d \rightarrow \mathbb{R}_+$. 
Consider a bipartite graph with nodes $L=\{\ell_e : e\in H'_d\}$ and $R=\{r_f : f\in H'_d\}$, and arcs $L\times R$. 
There is an additional source node $s$, and arcs from $s$ to each $L$-node and arcs from each $R$-node to $s$. 
Every arc $(i,j)$ in this network has a \textit{lower bound} $\alpha(i,j)$ and an \textit{upper bound} $\beta(i,j)$. 
The goal is to find a \textit{circulation} $z$ such that $\alpha(i,j) \le z(i,j)\le \beta(i,j)$ for all arcs $(i,j)$. 
Recall that a circulation is an assignment of non-negative values to the arcs of the network so that the in-flow equals the out-flow at every node. 
The lower/upper bounds are set as follows. }

\blue{
\begin{enumerate}
    \item For each arc  $(i,j)\in L\times R$, we have $\alpha(i,j)=0$ and $\beta(i,j)=\infty$.
    \item For each arc $(s,\ell_e)$ where $e\in H'_d$, we have $\alpha(s,\ell_e)=\beta(s,\ell_e)=x_e$. 
    \item For each arc $(r_f,s)$ where $f\in H'_d$, we have $\alpha(r_f,s)=\beta(r_f,s)=z_f$. 
\end{enumerate}
    Recall that $\x$ and $\z$ are the LP solutions corresponding to $\hat{v}(S_R)$ and $\hat{v}(S_O)$.}
    
\blue{Given \textit{any} circulation $z$, we define $\Delta_d(e,f) = z(\ell_e, r_f)$ for all $e,f\in H'_d$. Then, it is easy to see that all 3 conditions in Lemma~\ref{lem:1-swap} are satisfied.}

\blue{It just remains to prove the existence of some circulation. By  Hoffman's circulation theorem \citep{hoffman2003inequalities}, there is a circulation if and only if
\begin{equation}
    \label{eq:hoffman}
\alpha\left( \delta^-(T)\right) \le \beta\left( \delta^+(T)\right),\qquad \forall T \mbox{ subset of nodes}.\end{equation}
Above, $\delta^-(T)$ denotes all arcs from a node outside $T$ to a node inside $T$; similarly, $\delta^+(T)$ denotes all arcs from a node inside $T$ to a node outside $T$. This condition can be verified using the following cases:}
\begin{itemize}
    \item  \blue{$T\cap L\ne\varnothing$ and $T\cap R\ne R$. In this case, there is some arc from $L\times R$ in $\delta^+(T)$, so the RHS in \eqref{eq:hoffman} is $\infty$, which is clearly satisfies the condition.
    \item \blue{$T\cap L= \varnothing$.} If source $s\not\in T$ then $\alpha(\delta^-(T))=0$; so \eqref{eq:hoffman} is clearly true. If source $s\in T$ then $\beta(\delta^+(T))\ge \sum_{e\in H'_d} x_e = 1$ as all of $L$ lies outside $T$, and clearly $\alpha(\delta^-(T))\le 1$; so  \eqref{eq:hoffman} holds. }    
    \item  \blue{ $T\cap R=R$. If $s\in T$ then 
 $\alpha(\delta^-(T))=0$; so \eqref{eq:hoffman} is clearly true. If source $s\not \in T$ then $\beta(\delta^+(T))\ge \sum_{f\in H'_d} z_f = 1$ as all of $R$ lies inside $T$, and clearly $\alpha(\delta^-(T))\le 1$; so  \eqref{eq:hoffman} holds. }
 \end{itemize}

\subsection{Supplementary Results for MMO Algorithm} \label{sec:app4} Recall that $\hat{v}(S_R)=\sum_\xi \hat{v}(S_R,\xi)$ where $\hat{v}(S_R,\xi)$ is defined as the LP in \eqref{eq2.5}. So, we can write $\hat{v}(S_R)$ as the following LP:
\begin{align}
        \hat{v}(S_R) = \maximize_{x}& \frac{1}{N} \sum_\xi \sum_{e \in E(\xi)} v_{e} x^\xi_{e}    \\
        s.t. & \sum_{e\in E(\xi): j\in e } x^\xi_{e} \leq 1     & \forall j \in D(\xi)  \quad \forall \xi \nonumber \\
        & \sum_{e\in E(\xi): i\in e } x^\xi_{e} \leq 1     & \forall i \in S_A \cup S_B  \quad \forall \xi   \nonumber \\
        &  x^\xi_{e} = 0    & \forall e \sim S_A\setminus S_R  \quad \forall \xi  \nonumber   \\
        &  x^\xi_{e} \geq 0       & \forall e \in E(\xi)  \quad \forall \xi . \nonumber
\end{align}

For any vehicle $i$ and scenario $\xi$, set $F_{i,\xi}\subseteq E(\xi)$ denotes all the hyperedges incident to $i$ in scenario $\xi$.   Note that all variables $x^\xi_e$ with $e\sim S_A\setminus S_R $ are set to zero. So, it suffices to consider the LP with variables   $x^\xi_e$  for $e\in F_{i,\xi}$ and $i\in S_B\cup S_R$. 

\blue{We now consider the dual of the above LP (which has the same optimal value by strong duality). 
Let $G=S_A\cup S_B\cup \left( \cup_\xi D(\xi)\right)$ denote a combined groundset consisting of \textit{all} vehicles and demands from all scenarios. The dual variables are $u_{g, \xi}$ for all $g\in G$ and scenarios $\xi$.}
The dual LP is: 
\begin{align*} \label{eqapp4}
    \hat{v}(S_R) \,\,= &\,\,\minimize_u \sum_{\xi} \sum_{ g \in G  }    u_{g, \xi } \\ 
 &      s.t.\,\, \sum_{g\in e}  u_{g, \xi }  \geq \frac{v_e}{N} , & \forall e \in F_{i,\xi},\quad \forall \xi, \quad  \forall i\in S_R \cup S_B  \nonumber \\
    & \qquad \pmb{u} \geq 0. &  \nonumber
\end{align*}

\end{APPENDICES}


\bibliographystyle{informs2014trsc} 
\bibliography{references} 

\end{document}